\documentclass[11pt]{article}
\usepackage{lscape}
 \usepackage[ansinew]{inputenc}
 \usepackage[dvips]{graphicx}
\usepackage[psamsfonts]{amssymb} \usepackage[all]{xy}  \usepackage[all]{xy}
 \usepackage[psamsfonts]{eucal}
 \usepackage[psamsfonts]{eucal}
\usepackage{amssymb, amsmath}
\usepackage{amsthm}
\usepackage{geometry}
\usepackage{enumerate}
\usepackage{float}
\newcommand{\be}{\begin{equation}}
\newcommand{\ee}{\end{equation}}

\newcommand{\bd}{\begin{displaymath}}
\newcommand{\ed}{\end{displaymath}}
\newcommand{\ben}{\begin{enumerate}}
\newcommand{\een}{\end{enumerate}}
\newtheorem{Theorem}{Theorem}

\newtheorem{Prop}{Proposition}

\newtheorem{lem}{Lemma}

\title{Enriched Lie algebras in topology, I}
 \author{Yves F\'elix and Steve Halperin}
\begin{document}

\maketitle

To each path connected space $X$   the Sullivan theory of minimal models  associates a   commutative differential graded algebra, its minimal model $(\land V,d)$, and with it a graded Lie algebra $L_X$ that is the homotopy Lie algebra of its geometric realization $<\land V,d>$. When $X$ is a simply connected space with finite Betti numbers, then $L_X$ is isomorphic to the   Lie algebra of $X$, $ \pi_*(\Omega X)\otimes \mathbb Q$, equipped with the Whitehead bracket.

In the general case, $L_X$ is   a complete enriched Lie algebra: An \emph{enriched Lie algebra} is a graded Lie algebra $L= L_{\geq 0}$ together with a family of ideals $I_\alpha$ indexed by a partially ordered set with $\alpha \geq \beta$ if and only if $I_\alpha \subset I_\beta$, and such that each quotient $L/I_\alpha$ is a finite dimensional nilpotent Lie algebra. The enriched Lie algebra $(L, \{I_\alpha\})$ is \emph{complete} if $L= \varprojlim_\alpha L/I_\alpha$. 

Enriched Lie algebras are the natural extension of graded Lie algebras for path connected spaces, because each complete enriched Lie algebra $L$ is the homotopy Lie algebra of a path connected space $X$.  

This text is the first part of a general study of complete enriched Lie algebras.  Applications will be contained in a second part.

\newpage
\tableofcontents

\newpage
\part{Enriched Lie algebras}

\section{Introduction} By an \emph{enriched Lie algebra} we mean a graded Lie algebra, $L=L_{\geq 0}$, equipped with a family of ideals ${\mathcal I}= \{I_\alpha\}$ satisfying the following conditions
\begin{enumerate}
\item[(i)] Each $L_\alpha = L/I_\alpha$ is a finite dimensional nilpotent Lie algebra. In particular, for some $n_\alpha$, $L^{n_\alpha}= 0$, $L^k$ denoting the linear span of iterated Lie brackets of length $k$.
\item[(ii)] The index set, $\{\alpha\}$ is a directed set under the partial order given by
$$\alpha \geq \beta \Longleftrightarrow I_\alpha \subset I_\beta:$$
i.e. for each $\alpha, \beta$ there is a $\gamma$ such that $I_\gamma \subset I_\alpha \cap I_\beta$.
\item[(iii)] $\cap_\alpha I_\alpha = 0$.
\end{enumerate}
The corresponding quotient maps will be denoted 
$$\rho_\alpha : L\to L/I_\alpha:=L_\alpha.$$
Moreover, we shall use ${\mathcal I}$ to denote both the family $\{I_\alpha\}$ of ideals, and the index set $\{\alpha\}$.

\vspace{3mm}\noindent {\bf Example.} In any graded Lie algebra $L= L_{\geq 0}$ the ideals $J_\alpha \subset L$ for which $L/J_\alpha$ is finite dimensional and nilpotent form a directed system as above. The intersection of these ideals consists of the elements $x\in L$ which, for every $n$, can be expressed as a linear combination of commutators of length $n$. If this intersection is zero, the ideals $J_\alpha$ define an enriched structure, and so this condition is necessary and sufficient for $L$ to admit an enriched structure.
However, in principle, $L$ may admit two distinct such structures.

Finally, as we shall show in \S 8, a minimal Sullivan algebra induces a natural enriched structure in its homotopy Lie algebra.

\vspace{3mm}
The \emph{completion} of an enriched Lie algebra $(L, {\mathcal I})$ is the inclusion
$$\lambda_L: L\to \varprojlim_\alpha L/I_\alpha := \overline{L},$$
and $L$ is \emph{complete} if $\lambda$ is an isomorphism. In particular, the enriched structure in the homotopy Lie algebra of a minimal Sullivan algebra is complete. Now, in general, the inverse limit provides surjections $\overline{\rho_\alpha}: \overline{L} \to L/I_\alpha$, and their kernels make $\overline{L}$ into a complete enriched Lie algebra, $(\overline{L}, \overline{\mathcal I})$. 

Note that if $\beta>\alpha$ then $\rho_\alpha$ factors to yield a surjection $\rho_{\alpha\beta} : L_\beta\to L_\alpha$. Thus $\overline{L}$
 consists of the \emph{coherent} families $(x_\alpha)$, $\alpha \in \mathcal I$, i.e. those satisfying $\rho_{\alpha\beta}x_\beta = x_\alpha$. In particular, 
\begin{eqnarray}
\label{i1}
L\cap \mbox{ker}\, \overline{\rho_\alpha}= \mbox{ker}\, \rho_\alpha.
\end{eqnarray}

\vspace{3mm}\noindent {\bf Definition}. Suppose $R\subset E$ and $T\subset F$ are subspaces of enriched Lie algebras $(E, {\mathcal J})$ and $(F, {\mathcal I})$. A linear map $f : R\to T$ is \emph{coherent} if deg$\, f= 0$ and for each $\alpha \in {\mathcal I}$ there is a $\beta\in {\mathcal J}$ such that 
$$J_\beta \cap R \subset f^{-1}(I_\alpha \cap T).$$

In particular, a \emph{morphism} of enriched Lie algebras is a coherent homomorphism $\varphi : (E, {\mathcal J}) \to (F, {\mathcal I})$ of graded Lie algebras.
 In other words, for each $\alpha \in {\mathcal I}$ there is a  $\beta(\alpha) \in \mathcal G$ and a morphism $\varphi_\alpha : E_{\beta (\alpha)}\to F_\alpha$ such that $\rho_\alpha \circ \varphi = \varphi_\alpha \circ \rho_{\beta (\alpha)}$. Thus a morphism of complete enriched Lie algebras $\varphi: E\to F$ is a morphism of limits
$$E= \varprojlim_\beta E_\beta \to \varprojlim_\alpha E_{\beta (\alpha)} \to \varprojlim_\alpha F_\alpha = F.$$
It is therefore straightforward to verify that a morphism $\varphi : (E, {\mathcal G}) \to (F, {\mathcal I})$ extends uniquely to a morphism $\overline{\varphi}: (\overline{E}, \overline{\mathcal G})\to (\overline{F}, \overline{\mathcal I})$ satisfying
$$\overline{\varphi} \circ \lambda_E = \lambda_F\circ \varphi.$$

Next, if ${\mathcal I}$ and ${\mathcal I}'$ are families of ideals in $L$ satisfying the conditions above, we say that ${\mathcal I}$ and ${\mathcal I}'$ are \emph{equivalent} (${\mathcal I}\sim {\mathcal I}'$) if the identities, $(L, {\mathcal I})\to (L, {\mathcal I}')$ and $(L, {\mathcal I}')\to (L, {\mathcal I})$  are morphisms. This is an equivalence relation, and if $\varphi : (F, {\mathcal G})\to (E, {\mathcal I})$ is a morphism, then $\varphi : (F, {\mathcal G}') \to (E, {\mathcal I}') $ is also a morphism whenever ${\mathcal G}\sim {\mathcal G}'$ and ${\mathcal I}\sim {\mathcal I}'$. In particular, if ${\mathcal I}$ and ${\mathcal I}'$ are equivalent families of ideals, then the identity induces an isomorphism
$$\varprojlim_{I_\alpha \in \mathcal I} L/I_\alpha \stackrel{\cong}{\longrightarrow} \varprojlim_{J_\beta\in {\mathcal I}'} L/J_\beta.$$
Thus the completion, $\lambda_L: L\to \overline{L}$ is independent of the choice of ${\mathcal I}$ in its equivalence class.

\vspace{3mm}\noindent {\bf Examples.} (1). Let $(L,\{I_\alpha\})$ be an enriched Lie algebra and let $(I_\beta)$ be a subfamily such that for each $\alpha $ there is a $\beta$ with $I_\beta\subset I_\alpha$. Then $(I_\alpha)$ and $(I_\beta)$ are equivalent families of ideals.

(2). \emph{Full enriched Lie algebras}. An enriched Lie algebra $(L, {\mathcal I})= \{I_\alpha\})$ is \emph{full} if whenever an ideal $I$ contains some $I_\beta$ then $I= I_\alpha$ for some $\alpha$. If $(L, {\mathcal J})= \{J_\beta\})$ is any enriched Lie algebra then ${\mathcal J} $ extends to a unique full enriched family ${\mathcal I}$ consisting of all the ideals containing some $J_\beta$. Evidently ${\mathcal I}\sim {\mathcal J}$. Moreover if ${\mathcal I}$ and ${\mathcal J}$ are equivalent full enriched structure, then they coincide.

(3). Suppose that $L$ is the projective limit of a directed system of nilpotent finite dimensional Lie algebras $Q_\alpha$, $L=\varprojlim_\alpha Q_\alpha$. Denote $\rho_\alpha : L\to Q_\alpha$ the corresponding projections. Then $(L, \{\mbox{ker}\, \rho_\alpha\})$ is a complete enriched Lie algebra with $L_\alpha = L/\mbox{ker}\, \rho_\alpha= \mbox{Im}\, \rho_\alpha \subset Q_\alpha$, $$L = \varprojlim_\alpha L_\alpha.$$

(4). {\sl Free Lie algebras.}  The Lie algebra freely generated by a subspaces $S$ is denoted by $\mathbb L_S$ and, in particular, is the direct sum of the subspaces
$$\mathbb L_S(k) = [S, \dots , S] \hspace{3mm} (k \, \mbox{factors}.)$$
Its \emph{classical completion}, $\widehat{\mathbb L}_S$ is then the inverse limit $\varprojlim_n \mathbb L_S/\oplus_{k>n} \mathbb L_S(k).$\\
When $S$ is a graded vector space of finite type the ideals $\oplus_{k>n}  \mathbb L_S(k)$ define an enriched in $\mathbb L_S$, and $\widehat{\mathbb L}_S = \overline{\mathbb L_S}$, the corresponding completion.

(5). When the space $[L,L]$ of commutators in a graded Lie algebra $L$ has finite codimension then (see Proposition \ref{p4}) its completion $\overline{L}$ has a unique structure as an enriched Lie algebra, namely that provided by its lower central series. 

\vspace{3mm}\noindent {\bf Remark.} Let $f: (L, \{J_\beta\})\to (L', \{I_\alpha\})$ be a coherent morphism of complete enriched Lie algebras. For each $\alpha$, let $\beta (\alpha)$ an index such that $J_{\beta (\alpha)}\subset f^{-1}(I_\alpha)$ and form the index set $K= \{(\alpha, \beta)\, \vert \beta\geq \beta (\alpha)\}$. Then for each $(\alpha, \beta)\in K$, let $I_{(\alpha, \beta)} = I_\alpha$ and $J_{(\alpha, \beta)}= J_\beta$. This gives for each $(\alpha, \beta)\in K$ a morphism $f_{(\alpha, \beta)} = L_{(\alpha, \beta)} \to L'_{(\alpha, \beta)}$, and $f = \varprojlim_{(\alpha, \beta)} f_{(\alpha, \beta)}$. 

\vspace{3mm}
Now we define the category, ${\mathcal C}$, of enriched Lie algebras as follows:
\begin{enumerate}
\item[$\bullet$] The objects of ${\mathcal C}$ are the pairs $(L, \{{\mathcal I}\})$ in which $L=L_{\geq 0}$ is a graded Lie algebra and $\{{\mathcal I}\}$ is an equivalence class of ideals satisfying the conditions above.
\item[$\bullet$] The morphism of ${\mathcal C}$ are the morphisms $\varphi : (F,{\mathcal G})\to (E, {\mathcal I})$.
\end{enumerate}

In particular $L\leadsto \overline{L}$ is a functor from ${\mathcal C}$ to the sub category of complete enriched Lie algebras, the inclusion $L\to \overline{L}$ is a morphism, and any morphism $(E, {\mathcal J})\to (F, {\mathcal I})$ extends uniquely to a morphism $\overline{E}\to \overline{F}$.

\vspace{3mm} A key aspect of an enriched Lie algebra $(L, {\mathcal I})$ is that ${\mathcal I}$ makes the Lie algebra accessible to finite dimensional arguments, even when $L$ is not a graded vector space of finite type. In particular, enriched Lie algebras behave well with respect to inverse limits, as follows from the following Theorem (\cite{Bou}):

 \emph{Suppose $0\to A_\alpha\to B_\alpha\to C_\alpha \to 0$ is a directed system of short exact sequences of graded vector spaces of finite type. Then
\begin{eqnarray}
\label{i2}
0\to \varprojlim_\alpha A_\alpha \to \varprojlim_\alpha B_\alpha \to \varprojlim_\alpha C_\alpha \to 0
\end{eqnarray}
is also short exact. In particular, for chain complexes of finite type, homology commutes with inverse limits.}

\vspace{3mm} Complete enriched Lie algebras arise naturally in Sullivan's approach to rational homotopy theory :   the homotopy Lie algebra of a Sullivan model is naturally a complete enriched Lie algebra and this is an essential aspect when the theory is extended to all connected CW complexes. Thus our objective here  is  to develop the properties of enriched Lie algebras, in particular with a view to future topological applications. In particular, while many of the basic properties of ${\mathcal C}$ are parallel those for the category of all graded Lie algebras, the enrichment gives the objects in ${\mathcal C}$ a distinct "homotopy flavour". In particular each enriched Lie algebra $(L, {\mathcal I})$ determines a quadratic Sullivan algebra $\land Z$ and a simplicial set $\langle \land Z\rangle$.

\section{Subspaces, sub Lie algebras, ideals, and universal enveloping algebras} 

\subsection{General properties}

\emph{Throughout this entire subsection $(L, {\mathcal I})$ denotes a fixed complete enriched Lie algebra, with quotient maps $\rho_\alpha: L\to L/I_\alpha:=L_\alpha$.}

Any subspace $S\subset L$ determines the inclusions,
$$S\to \varprojlim_\alpha S/S\cap I_\alpha = \varprojlim_\alpha \rho_\alpha S \subset L,$$
and the subspace $\varprojlim_\alpha \rho_\alpha S\subset L$ is called the \emph{closure} of $S$ in $L$. In particular, \emph{the closure of a sub Lie algebra $E\subset L$ coincides with its completion, $\overline{E}$,} and so we use the same notation and write
$$\overline{S}= \varprojlim_\alpha \rho_\alpha S.$$
Since $\rho_\beta S\to \rho_\alpha S$ is surjective for $\beta \geq \alpha$ it follows that
$$\rho_\alpha \overline{S}= \rho_\alpha S,\hspace{5mm} \alpha \in {\mathcal I}.$$

\vspace{3mm}\noindent {\bf Remark.} Suppose $(x_k)_{k\geq 1}$ is a sequence of elements in $S$. If for each $\alpha \in {\mathcal I}$ there is some $r(\alpha)$ such that
$$\sum_{k=r(\alpha) + 1}^m \, x_k \in I_\alpha, \hspace{1cm} m>r(\alpha),$$
then the elements $\rho_\alpha (\sum_{k=1}^{r(\alpha)} x_k) \in \rho_\alpha (S)$ define a single element $y\in \overline{S}$,  
 and we write $y = \sum_k x_k.$

\begin{lem}
\label{l1} If $S\subset T\subset L$ are subspaces, then
$$\overline{T}/\overline{S} = \varprojlim_\alpha \rho_\alpha T/\rho_\alpha S.$$
\end{lem}

\vspace{3mm}\noindent {\sl proof:} Since $\rho_\alpha S\subset \rho_\alpha T$,
$$0\to \rho_\alpha S\to \rho_\alpha T\to \rho_\alpha T/\rho_\alpha S\to 0$$
is a short exact sequence of finite dimensional spaces. Thus the Lemma follows from (\ref{i2}).

 \hfill$\square$

\begin{lem}
\label{l2}  
\begin{enumerate}
\item[(i)] A finite sum of closed subspaces of $L$ is closed.  
\item[(ii)] An arbitrary intersection of closed subspaces of $L$ is closed.
\item[(iii)] For any $S\subset L$, $\overline{S}= \cap_\alpha (S+I_\alpha)$.
\item[(iv)] If $S\subset L$   is a graded space of finite type, then $S$ and $[S,L]$ are closed.

\end{enumerate}\end{lem}

\vspace{3mm}\noindent {\sl proof.} (i) It is sufficient to show that if $S$ and $T$ are closed then $S+T$ is closed. By Lemma \ref{l1},
$$\renewcommand{\arraystretch}{1.7}
\begin{array}{ll}
\overline{S+T}/\overline{S} &= \displaystyle\varprojlim_\alpha \rho_\alpha (S+T) / \rho_\alpha (S) = \displaystyle\varprojlim_\alpha \displaystyle\frac{\rho_\alpha (S) + \rho_\alpha (T)}{\rho_\alpha (S)}
\\
& = \displaystyle\varprojlim_\alpha \displaystyle\frac{\rho_\alpha (T)}{\rho_\alpha (S) \cap \rho_\alpha (T)}\\
& = \displaystyle \overline{T} / \displaystyle \varprojlim_\alpha \rho_\alpha (S) \cap \rho_\alpha (T)
\end{array}
\renewcommand{\arraystretch}{1}$$
In particular, it follows that the inclusion of $\overline{S}$ and $\overline{T}$ in $\overline{S+T}$ define a surjection $\overline{S}+ \overline{T}\to \overline{S+T}$. Since $S$ and $T$ are closed, $S+T\to \overline{S+T}$
is an isomorphism.

\vspace{3mm} (ii) Suppose $x\in \overline{\cap_\sigma S_\sigma}$, where each $S_\sigma$ is a closed subspace of $L$. Then $$\rho_\alpha x\in \rho_\alpha (\cap S_\sigma) \subset \cap_\sigma \rho_\alpha (S_\sigma).$$
Therefore $\rho_\alpha x \in \rho_\alpha (S_\sigma)$ for each $\sigma$. Thus $x = (\rho_\alpha x)\in \varprojlim_\alpha \rho_\alpha (S_\sigma) = \overline{S_\sigma} = S_\sigma$. Thus $x\in \cap_\sigma S_\sigma$ and $\cap_\sigma S_\sigma = \overline{\cap_\sigma S_\sigma}$. 

\vspace{3mm} (iii) If $x\in \overline{S}$ then $\rho_\alpha x\in \rho_\alpha (S)$. Thus for some $x_\alpha \in S$, 
$$\rho_\alpha (x-x_\alpha) = 0$$
and so $x-x_\alpha \in I_\alpha$. Therefore $x\in S+I_\alpha$ for all $\alpha$.

In the reverse direction, suppose $x\in \cap_\alpha (S+I_\alpha)$. Then
$$\rho_\alpha x\in \rho_\alpha S.$$
Thus the coherent family $x = (\rho_\alpha x) \in \overline{S}$. 

\vspace{3mm} (iv) In view of (i) it is sufficient to prove that for any $x\in L$, $[x,L]$ is closed. Let $C_\alpha = \{y_\alpha \in L_\alpha \,\vert\, [\rho_\alpha x, y_\alpha]= 0\}.$ It is immediate from the definition that $\{C_\alpha\}$ is an inverse system, and we set
$$C = \varprojlim_\alpha C_\alpha.$$

Now we show that $C = \{y\in L\, \vert\, [x,y]= 0\}.$  Indeed, if $y\in C$, then for all $\alpha$,
$$\rho_\alpha [x,y] = [\rho_\alpha x, \rho_\alpha y] \in [\rho_\alpha x, C_\alpha] = 0,$$
and hence $ [x,y]= 0$. On the other hand, if $[x,y]= 0$ then $[\rho_\alpha x, \rho_\alpha y]= \rho_\alpha [x,y]= 0$ and so each $\rho_\alpha y\in C_\alpha$. Thus $y\in C$. 

Finally, apply Lemma \ref{l1} to obtain the commutative diagram
$$\xymatrix{L/C = \varprojlim_\alpha L_\alpha \, / \, \varprojlim_\alpha C_\alpha\ar[d]_\cong^{[x,-]} \ar[rrr]^-= &&&\varprojlim_\alpha L_\alpha/C_\alpha\ar[d]_\cong^{[x,-]}\\
[x,L] \ar[rrr]^-= &&&\varprojlim_\alpha [\rho_\alpha x, L_\alpha]}$$
Since $\varprojlim_\alpha [\rho_\alpha x, L_\alpha]= \overline{[x,L]}$, (iv) is established. \hfill$\square$

\vspace{3mm}\noindent {\bf Definition.} (i). If $E\subset L$ is a sub Lie algebra then the corresponding sub enriched Lie algebra $(E, {\mathcal I}_E)$ is defined by
$${\mathcal I}_E = \{E\cap I_\alpha\}.$$

(ii). If $I\subset L$ is a closed ideal then the corresponding quotient enriched Lie algebra, $(L/I, {\mathcal I}_{L/I})$, defined by
$${\mathcal I}_{L/I} = \{ (I+I_\alpha)/I\},$$
is complete.

Straightforward arguments from the definitions then establish the 

\vspace{3mm}\noindent {\bf Remarks.}
\begin{enumerate}
\item[1.] If $E\subset F\subset L$ are sub Lie algebras, then $(E, {\mathcal I}_E)\to (F, {\mathcal I}_F)$ is a morphism. Moreover, the closure and completion of $E$ coincide as complete enriched Lie algebras. 
\item[2.]   If $I$ is an ideal in $L$ then   $\overline{I}$ is an ideal.
\item[3.] If $I$ is a closed ideal in $L$ then $(L/I, {\mathcal I}_{L/I})$ is complete. (This requires Lemma \ref{l1}.)
\item[4.] Suppose $(E, {\mathcal G})\stackrel{\varphi}{\to} (L, {\mathcal I})$ is a morphism of complete enriched Lie algebras. Then ker$\,\varphi$ is closed in $E$, and so $(\mbox{ker}\, \varphi, {\mathcal I}_{\mbox{\scriptsize ker}\, \varphi})$ is a complete enriched Lie algebra. The induced map $E/\mbox{ker}\, \varphi \to L$ is a morphism of enriched Lie algebras.
\item[5.] If $L = \varinjlim_\sigma L(\sigma)$ is the direct limit of closed sub Lie algebras $L(\sigma)$ then for each $\alpha$ there is some $\sigma$ for which $I_\alpha = L(\sigma)$, and
$$C(L,-) = \varinjlim_\sigma C(L(\sigma), -).$$

\end{enumerate}

\subsection{Weight decompositions}

\vspace{3mm}\noindent   A \emph{weight decomposition} in a graded Lie algebra is a decomposition
$$L= \oplus_{k\geq 1} L(k)$$
in which $[L(k), L(\ell)] \subset L(k+\ell)$. A \emph{weighted subspace} of $L$ is a subspace $S\subset L$ such that $S = \oplus_k S\cap L(k)$. We denote $S(k) = S\cap L(k)$.

A \emph{weighted enriched Lie algebra} is a weighted Lie algebra with a defining set of ideals of the form $I_\alpha = \oplus_k I_\alpha (k)$. 

\vspace{2mm} \emph{Henceforth in this example, $L= \oplus_{k\geq 1}L(k)$ denotes a fixed weighted   enriched Lie algebra with a defining set of weighted ideals $I_\alpha$.}

\vspace{2mm} There follows the

\begin{Prop} \label{p1}
With the hypotheses and notation above,  
\begin{enumerate}
\item[(i)] 
$\overline{L} = \prod_{k} \overline{L(k)}.$
\item[(ii)] If $S\subset L$ is a weighted subspace then $\overline{S} = \prod_k \overline{S(k)}$. In particular, $S$ is closed if and only if each $S(k)$ is closed, in which case any subspace of the form $\prod_i S(k_i)$ is closed.
\item[(iii)] If each $L(k)$ is finite dimensional, then  $\overline{L} = \varprojlim_r L/\, \prod_{k>r} L(k).$
\end{enumerate}\end{Prop}

 \vspace{3mm}\noindent {\sl proof.} (i) Denote as usual by $\rho_\alpha : L\to   L_\alpha= L/I_\alpha$ the projection on the finite dimensional Lie algebra $L_\alpha$. Then
 $$L_\alpha = \oplus_k \rho_\alpha (L(k)) = \prod_k \rho_\alpha (L(k)),$$
and
 $$L=\varprojlim_\alpha L_\alpha = \prod_k \lim_\alpha \rho_\alpha (L(k)) = \prod_k \varprojlim_\alpha  L(k)/I_\alpha (k) = \prod_k \overline{L(k)}.$$
 
 (ii) Similar proof than the one for (i), and (iii) is a direct consequence. 
 
 \hfill$\square$

\subsection{Universal enveloping algebras}
 
\vspace{1mm} Recall that the classical completion of the universal enveloping algebra $UE$ of a graded Lie algebra $E$ is the inverse limit
 $$\widehat{UE} = \varprojlim_n UE/J^n,$$
 $J^n$ denoting the $n^{th}$ power of the augmentation ideal. In particular
 $$\widehat{UE} = \mathbb Q \oplus \widehat{J},$$
 and $\widehat{J} := \varprojlim_n J/J^n$ is the augmentation ideal for $\widehat{UE}$.

 \vspace{3mm}\noindent {\bf Definition.} 
The \emph{Sullivan completion} $\overline{UL}$ of an enriched Lie algebra is defined  by
 $$\overline{UL}:= \varprojlim_\alpha \widehat{UL_\alpha} = \varprojlim_{n,\alpha} UL_\alpha/J_\alpha^n= \varprojlim_\gamma \overline{UL}/K_\gamma,$$
 where $K_\gamma$ runs over all finite codimensional ideals of $UL$. 
 
 \vspace{2mm}Passing to inverse limits shows that the inclusions $L_\alpha \to L_\alpha/J_\alpha^n$ define an inclusion
 $$\overline{L}\to \overline{UL}.$$ Note also that   $$ \overline{UL} = \mathbb Q \oplus \overline{J},$$
 where $\overline{J} = \varprojlim_{\alpha, n} J_\alpha/J_\alpha^n$; $\overline{J}$ is the augmentation ideal for $\overline{UL}$. It is also immediate that 
  a morphism $\varphi : E\to L$ of enriched Lie algebras extends to a morphism $\overline{U\varphi}: \overline{UE}\to \overline{UL}$.

 Finally, in analogy with the filtration of $UL$ by the ideals $J^n$ we filter $\overline{UL}$ by the ideals
 $$\overline{J}^{(n)}:= \varprojlim_{k,\alpha} J_\alpha^n/J_\alpha^{k+n} \subset \overline{J}= \overline{J}^{(1)}.$$
 Then $\overline{UL}$ is complete with respect to this filtration:
 \begin{eqnarray}\label{i3}
 \overline{UL}= \varprojlim_n \overline{UL}/\overline{J}^{(n)}.
 \end{eqnarray}
 In fact, (\ref{i2}) yields
 $$
 \renewcommand{\arraystretch}{1.4}
 \begin{array}{ll}
 \overline{UL} = \varprojlim_{n,\alpha} UL_\alpha /J_\alpha^n & = \varprojlim_k \varprojlim_{n,\alpha} (UL_\alpha/J_\alpha^{n+k})/ (J_\alpha^n/J_\alpha^{n+k})\\
 & = \varprojlim_n \varprojlim_{k,\alpha} (UL_\alpha/J_\alpha^{n+k})/ \varprojlim_{k,\alpha} (J_\alpha^n/J_\alpha^{n+k})\\
 & = \varprojlim_n \overline{UL}/\overline{J}^{(n)}.
 \end{array}
 \renewcommand{\arraystretch}{1}
$$ 
 
\vspace{3mm}\noindent {\bf Examples.} 
\begin{enumerate}
\item When dim$\, L/[L,L]<\infty$, then $\overline{UL}= \widehat{UL}$.
\item Let $L= \widehat{\mathbb L}(V)$ be the completion of the free graded Lie algebra on a finite dimension graded vector space $V$. Then $\overline{UL}$ is isomorphic to $\prod_n T^nV$. Moreover the short exact sequence of $\overline{UL}$-modules
$$0\to \overline{UL}\otimes sV\stackrel{d}{\to} \overline{UL} \stackrel{\varepsilon}{\to} \mathbb Q\to 0,$$
with $d(sx)=x$ shows that $\overline{UL}$ is an algebra of global dimension one.
\end{enumerate}

\section{Lower central series}
\emph{Throughout this entire section $(L, {\mathcal I})$ denotes a fixed complete enriched Lie algebra, with quotient maps $\rho_\alpha: L\to L/I_\alpha:=L_\alpha$.} 

\vspace{2mm} A key invariant for any graded Lie algebra $E$ is its lower central series
$$E= E^1\supset \dots \supset E^k\supset \dots$$
of ideals in which $E^k$ is the linear span of iterated commutators of length $k$ of elements in $E$. 
The \emph{classical completion} of $E$ is the inverse limit $\widehat{E}=\varprojlim_n E/E^n$, and $E$ is \emph{pronilpotent} if $E\stackrel{\cong}{\to} \widehat{E}$.

\vspace{3mm}The analogue of the lower central series for   $L$ is the sequence of ideals
$$L=L^{(1)}\supset \dots \supset L^{(k)}\supset \dots$$
defined by $L^{(k)}= \overline{L^k}$. It is immediate that $L\leadsto L^{(k)}$ is a functor. Moreover, since $L^k\subset L^{(k)}$, \emph{each $L/L^{(k)}$ is a nilpotent Lie algebra : $\left(L/L^{(k)}\right)^k = 0$.}

\begin{lem}
\label{l3}  Let $L$ be a complete enriched Lie algebra. \begin{enumerate}
\item[(i)] $ L^{(k)} = \varprojlim_\ell L^{(k)}/L^{(k+\ell)}.$ In particular
$$L = \varprojlim_k L/L^{(k)}.$$
\item[(ii)] $[L^{(k)}, L^{(\ell)}] \subset L^{(k+\ell)}.$
\item[(iii)] $L$ is a retract of $\widehat{L} =\varprojlim_n L/L^n$.
\end{enumerate}
\end{lem}

\vspace{1mm}\noindent {\sl proof.} (i) Since each $L_\alpha$ is nilpotent, for some $\ell = \ell (\alpha)$, $L_\alpha^{k+\ell}= 0$. Thus by Lemma \ref{l1},
$$\varprojlim_\ell \overline{L^k}/\overline{L^{k+\ell}} = \varprojlim_{\ell, \alpha} L_\alpha^k/L_\alpha^{k+\ell}= \varprojlim_\alpha L_\alpha^k = \overline{L^k}.$$

(ii) From (2) we obtain $\rho_\alpha (L^{(k)})= \rho_\alpha (\lim_{\beta \geq \alpha} \rho_\beta (L^k))= \rho_\alpha (L^k) = L_\alpha^k$. It follows that
$$\rho_\alpha [L^{(k)}, L^{(\ell)}] \subset [L_\alpha^k, L_\alpha^\ell] \subset L_\alpha^{k+\ell},$$
and so $[L^{(k)}, L^{(\ell)}] \subset \overline{[L^{(k)}, L^{(\ell)}]} \subset L^{(k+\ell)}.$
 
 (iii) For each $n$ and $\alpha$, the surjection $L\to L_\alpha/L_\alpha^n$ (induced from $\rho_\alpha$) factors as
 $$L\to L/L^n \to L_\alpha/L_\alpha^n.$$
 Thus we obtain
 $$L\to \varprojlim_n L/L^n \to \varprojlim_{n,\alpha} L_\alpha /L_\alpha^n = \varprojlim_\alpha L_\alpha = L,$$
 which decomposes id$_L$ as $L\stackrel{\varphi}{\to} \widehat{L}\stackrel{\psi}{\to} L.$

 \hfill$\square$

\vspace{3mm}\noindent {\bf Example.} The completion of the free Lie algebra on an infinite number of generators may not be pronilpotent.
 
\vspace{2mm} Let $W$ be the vector space whose basis is given by the elements $v_i$ and $w_i$, $i\geq 1$, deg$\, v_i=$ deg$\, w_i=1$. We equip $W$ with   trivial multiplication $W\cdot W = 0$ and we denote by $(\land V,d)$ the minimal model of $(\mathbb Q \oplus W,0)$. So $V = V^1= \oplus_{n\geq 0} V_n$ with $V_0= V\cap \mbox{ker}\, d$ and for $n>0$, $d : V_n \to (\land^2 V)_{n-1}$. In particular $V_0= W$. Then, for $i\geq 1$ we denote by $c_{i1}$   an element of $V_1$ with $dc_{i1}= v_iw_i$, and by induction for $n>1$. Let   $c_{in}$ be an element of $V_n$ with $d(c_{in})= v_ic_{i,n-1}$.

Let $L$ be the homotopy Lie algebra of $(\land V,d)$. We denote by $x_i$ and $y_i$ the elements of $L$ defined by 
$$<v_i, sx_j >= \delta_{ij}, \hspace{3mm} <w_i,sx_j> = 0, \hspace{3mm}\mbox{and } <V_{\geq 1}, sx_j>= 0,$$
 $$<w_i, sy_j >= \delta_{ij}, \hspace{3mm} <v_i,sy_j> = 0, \hspace{3mm}\mbox{and } <V_{\geq 1}, sy_j>= 0.$$
 It follows that 
 $$<c_{j1}, s[x_i,y_i]> = -<v_j\land w_j, sx_i, sy_i>= -\delta_{ij}.$$
 In the same way,
  $$<c_{jn}, s\,\mbox{ad}^p_{x_i}(y_i)> = -<v_j\land c_{j,n-1}, sx_i, s\,\mbox{ad}^{p-1}_{x_i}(y_i)>= (-1)^n\delta_{ij}\delta_{pn}.$$

Now let $U\subset V$ be the subspace formed by the elements $z$ such that $<z, sx_i>= <z, sy_j>= <z, s \,\mbox{ad}^p_{x_k}(y_k)>= 0$ for all $i,j,k,p$. Then if $z\in V$ is such that $dz= v_i\land  c_{j,n} + \sum u_k\land z_k$ with $i\neq j$ and $u_k\in U$ then $z\in U$. The space $U$ is a direct summand of the vector space $S$ generated by the $v_i, w_j$ and $c_{kn}$.  To see that, let $z\in V_r$. Since $W = \cup_n W(n)$ where $W(n)$ is the vector space generated by the $v_i$ and $w_i$, $i\leq n$, $\land V= \cup_n \land V(n)$ where $\land V(n)$ is the minimal model of $(\mathbb Q \oplus W(n), 0)$. Therefore $v\in V(s)$ for some $s$ and the elements $x_i, y_i, \mbox{ad}^p_{x_i}(y_i)$ vanish on $v$ when $i>s$. Since $v\in V^r$ when $p\neq r$, for any $i$, the element $\mbox{ad}^p_{x_i}(y_i)$ vanish on $v$. Therefore the element
$$z_s = \sum_i <z, sx_i> v_i + \sum_j  <z, sy_j> w_j + \sum_{kp} <z, s\, \mbox{ad}^p_{x_k}(y_k)> c_{kp}$$
is a well defined element in $S$ and by construction $ z-z_s\in U$.

 Let $E$ be the free Lie algebra on the generators $x_i,y_i$, $i\geq 1$.    Then $L$ is the completion of $E$.
In particular the series  $$\omega = \sum_{n\geq 1} \mbox{ad}^n_{x_n}(y_n) = [x_1, y_1] + [x_2,[x_2,y_2]] + \dots$$
 is a well defined element in $L$. 
 
 We first suppose that $\omega$ belongs to $L^2$ and will arrive to a contradiction. This will imply that $\omega$ is an indecomposable element in $L$. So suppose $\omega = \sum_{i=1}^N [\omega_i, \omega_i']$. Denote by $K$ the intersection $K = \cap_{i=1}^N \mbox{ker}\, \omega_i \, \cap \, \cap_{i=1}^N \mbox{ker}\, \omega'_i$. Then $K$ is a finite codimensional subspace of $V$. Since the $v_i$ generate a subspace of infinite dimension, some linear combination $a=\sum_{i=p}^m \alpha_i v_i$,  $\alpha_i\in \mathbb Q$,  belongs to $K$, and we can suppose that the first term $\alpha_p= 1$. We consider then the sequence of elements $a_r\in V_r$, $r\geq 1$, defined by $da_1= ay_p$ and for $r\geq 2$, $da_r = a\cdot a_{r-1}$.
 Then $$<a_r, s\omega> = - \sum_{i=1}^N <da_r, s\omega_i, s\omega'_i > = 0.$$
 
 On the other hand, $a_1 = c_{p1}+ u_1$ with $u_1\in U$, and by induction $a_k = c_{pk}+ u_k$ with $u_k\in U$. Now we compute $<a_p, s\omega>$ using the decomposition of $\omega$ as a series. Since $a_p\in V_p$, for each iterated Lie bracket $\alpha$ in the $x_i, y_i$ of length different of $p$, we have $<a_p, \alpha>= 0$. Therefore
 $$<a_p, s\omega> = <a_p, s \,\mbox{ad}^p_{x_p}(y_p)> = <c_{pp}, s\,\mbox{ad}^p_{x_p}(y_p)> = -1,$$
 which gives the required contradiction.
 
 Now write $\alpha_n = \sum_{k=1}^n \mbox{ad}^k_{x_k}(y_k)$. Then $\alpha_{n+1}-\alpha_n \in L^{n+1}$ and the sequence $(\alpha_n)$ is a coherent sequence of elements in the tower 
 $$\xymatrix{\dots \ar[rr]&& L/L^{n+1} \ar[rr]^{q_{n+1,n}} && L/L^n \ar[rr]^{q_{n, n-1}} && L/L^{n-1} \ar[rr] && \dots}$$ This sequence defines an element $\alpha \in \widehat{L}= \varprojlim_n L/L^n$. 
 
We use the notation of Lemma \ref{l3} for the retraction $L\stackrel{\varphi}{\to} \widehat{L} \stackrel{\psi}{\to} L$. Since the natural projections $q_n : L/L^n \to L/L^{n}$ maps $\alpha_n$ to $\alpha_n$, at the limit we have $\psi (\alpha)= \omega$. Now if $\varphi$ is an isomorphism, because $\psi\circ \varphi$ is the identity, we would have $\varphi (\omega)= \alpha$, which is not the case because $\omega$ being indecomposable, in $L/L^2$ we have  $  p_2(\omega) = \omega \neq p_2(\alpha)= \alpha_1$.

\vspace{3mm}\noindent {\bf Remark.} While $L$ may not be pronilpotent, Lemma \ref{l3} identifies the sequence $L\supset \dots \supset L^{(k)}\supset \dots$ as an $N$-suite as defined by Lazard in \cite{La}.

 \vspace{3mm}\noindent {\bf Example.}  Let $E$ be an enriched Lie algebra. Then $\overline{E}$ and $\widehat{E}$ can be very different.
Let $E$ be the abelian Lie algebra with basis the countably infinite family $x_1, x_2, \dots$. We equip $E$ with the enriched structure given by the ideals $I_n$ generated by the $x_i$, $i\geq n$. The completion $\overline{E}$ is the vector space of series $\sum_{i\geq 1} \alpha_i x_i$, $\alpha_i\in\mathbb Q$. Therefore the injection $\widehat{E}\to \overline{E}$ is not an isomorphism, even when $E$ is nilpotent.

\begin{lem}
\label{l4} Suppose $E\subset L$ is a sub Lie algebra.
\begin{enumerate}
\item[(i)] If $E+ L^{(2)} = L$, then $\overline{E}= L$.
\item[(ii)]   $\overline{E^n} = \overline{E}^{(n)}$, $n\geq 2$. \end{enumerate}\end{lem}

\vspace{3mm}\noindent {\sl proof.} (i) By hypothesis $\rho_\alpha L = \rho_\alpha E+ \rho_\alpha L^{(2)}= \rho_\alpha E + [\rho_\alpha L, \rho_\alpha L]$. Since $\rho_\alpha L$ is nilpotent this gives $\rho_\alpha E = \rho_\alpha L$ and $\overline{E}= L$.

(ii) Since $E\subset \overline{E}$, $\rho_\alpha (\overline{E})= \rho_\alpha (E)$. Therefore
$$\renewcommand{\arraystretch}{1.4}
\begin{array}{ll} \overline{E^n} & = \varprojlim_\alpha \rho_\alpha (E^n) = \varprojlim_\alpha (\rho_\alpha (E))^n\\ & = \varprojlim_\alpha (\rho_\alpha (\overline{E})^n = \varprojlim_\alpha \rho_\alpha (\overline{E}^n)= \overline{E}^{(n)}.\end{array}\renewcommand{\arraystretch}{1}$$ \hfill$\square$

\begin{Prop}
\label{p2} Suppose $L$ is a complete enriched Lie algebra. If dim$\, L/L^{(2)}<\infty$, then
\begin{enumerate}
\item[(i)] For any integer $k$ there is an $\alpha$ such that the projection $L/L^{(k)}\to L_\alpha /L_\alpha^k$ is an isomorphism.
\item[(ii)] $L^k= L^{(k)}, \hspace{3mm}\mbox{ for } k\geq 1$.
\end{enumerate}
\end{Prop}

\vspace{3mm}\noindent {\sl proof.} (i).  Since $L/L^{(2)}$ is finite dimensional,  by Lemma \ref{l1} there is an $\alpha_0$, such that for each $\beta\geq \alpha_0 $  the projection $p_\beta : L/L^{(2)}\to L_\beta/L_\beta^2$ is an isomorphism. This shows that for each $\alpha $ the dimension of $L_\alpha/L_\alpha^2$ is bounded by the integer $N = \mbox{dim}\, L_{\alpha_0}/L_{\alpha_0}^2$. It follows that for any $r\geq 1$, and any $\alpha$ the dimension of $L_\alpha^r/L_\alpha^{r+1}$ is bounded by $N^r$. Thus it follows from Lemma \ref{l1} that   $L^{(r)}/L^{(r+1)}= \varprojlim_\alpha L_\alpha^r/L_\alpha^{r+1}$ is finite dimensional. Therefore for any integer $k$ there is an index $\alpha_k$ such that the projection 
$$L/L^{(k)}\to L_{\alpha_k}/L_{\alpha_k}^k$$
is an isomorphism.

(ii) 
It follows from Lemma \ref{l1} that for some $\alpha$ the Lie bracket in $L$  
\begin{eqnarray}\label{i4}
L/L^{(2)} \otimes L^{(r)}/L^{(r+1)} \to L^{(r+1)}/L^{(r+2)}, \hspace{1cm} r\geq 1,
\end{eqnarray}
 may be identified with the corresponding surjection
$$L_\alpha/L_\alpha^2\otimes L_\alpha^r/L_\alpha^{r+1} \to L_\alpha^{r+1} / L_\alpha^{r+2}.$$
This implies that the Lie bracket (\ref{i4}) is   surjective.

Now let $x_1, \dots , x_m\in L$ represent a basis of, $L/L^{(2)}$. Then for any
$y\in L^{(r+1)}$ it follows   that there are elements $y_i\in L^{(r)}$ such that
$$y-\sum_{i=1}^r [x_i, y_i]\in L^{(r+2)}.$$
Given $x\in L^{(n)}$, this yields an inductive construction  of elements $y_i(\ell)\in L^{(\ell-1)}$, $\ell\geq k$, such that
$$x-\sum_{i=1}^m\sum_{\ell = k-1}^r [x_i, y_i(\ell)]\in L^{(r)}.$$
Set $y_i=\sum_\ell y_i(\ell)$. Then by construction
$$x= \sum_{i=1}^m [x_i, y_i].$$

This shows that $L^{(k)}\subset [L, L^{(k-1)}]$. Induction on $k$ gives
$$L^{(k)}\subset L^k,$$
and the reverse inclusion is obvious.

\hfill$\square$

\vspace{3mm}\noindent {\bf Corollary 1.}  Let $E$ be an enriched Lie algebra with dim$\, E/E^2<\infty$. Then dim$\,\overline{E}/\overline{E^k}<\infty$, $k\geq 1$, and 
the lower central series for $\overline{E}$ satisfies
$$\overline{E}^{(k)}= \overline{E}^k = \overline{E^k}.$$ 

\vspace{3mm}\noindent {\sl proof.} There is an 
$\alpha_0$ such that for each $\alpha \geq \alpha_0$, the projection $E_\alpha /E_\alpha^2\to E_{\alpha_0}/E_{\alpha_0}^2$ is an isomorphism. Thus by Lemma \ref{l1} $\overline{E}/\overline{E}^{(2)}$  is isomorphic to $  E_{\alpha_0}/E_{\alpha_0}^2$, and so is   finite dimensional. Proposition \ref{p2} gives then an isomorphism $\overline{E}/\overline{E}^{(2)}= \overline{E}/\overline{E}^2=E_{\alpha_0}/E_{\alpha_0}^2$ In particular, $E+\overline{E}^{(2)}= \overline{E}$. 
The result is then a direct consequence of Proposition \ref{p2}(ii) and  Lemma \ref{l4}(ii).  

\hfill$\square$

\vspace{3mm}\noindent {\bf Corollary 2.}  An complete enriched Lie algebra $L$ is the direct limit of closed   sub Lie algebras, $E$, satisfying dim$\, E/E^2<\infty$.

\vspace{3mm}\noindent {\sl proof.} Any Lie algebra is the direct limit of its finitely generated Lie algebras, $F$. Thus $L$ is the direct limit of the completions, $\overline{F}$. By Proposition \ref{p2} and Corollary 1, each $F$   satisfies dim$\, \overline{F}/\overline{F}^2<\infty$. 

\hfill $\square$

\vspace{3mm}Proposition \ref{p2} has the following analogue:

\begin{Prop}
\label{p3} If $E\subset L$ is a sub Lie algebra and $E/E^2$ is a graded vector space of finite type then 
\begin{enumerate}
\item[(i)] Each $\overline{E}/\overline{E^k}$ has finite type.
\item[(ii)] $\overline{E}^k = \overline{E^k}=\overline{E}^{(k)}.$
\end{enumerate}
\end{Prop}

\vspace{3mm}\noindent {\sl proof:} (i)  In the proof of Proposition \ref{p2}(i), replace $L/L^k$ by any $E_n/(E^k)_n$ to obtain that $(\overline{E})_n / (\overline{E}^k)_n$ is finite dimensional. 

\vspace{2mm} (ii)  This follows directly by the same proof as in Corollary 3..

\hfill$\square$

\begin{Prop}
\label{p4} 
 Suppose $L$ is a complete enriched Lie algebra such that
 $L/L^2 \mbox{ has finite type,}$ 
then $L$ has a unique (up to equivalence) structure as an enriched Lie algebra $(L, {\mathcal I})$, namely that given by
$${\mathcal I}= \{I_k\}_{k\geq 1} \hspace{5mm}\mbox{with } I_k = L^k+ L_{\geq k}.$$
\end{Prop}

\vspace{3mm}\noindent {\sl proof:} Suppose a second family ${\mathcal J}= \{J_\beta\}$ of ideals in $L$ also makes $L$ into an enriched Lie algebra. Since $L/J_\beta$ is finite dimensional and nilpotent, it is also immediate that $J_\beta \supset I_k$ for some $k$.

In the reverse direction, apply Proposition \ref{p2}(i) to conclude that for fixed $k$ and $n$ there is some $\beta = \beta (k,n)$ such that 
$$\rho_\beta : L_n /(L^k)_n \stackrel{\cong}{\longrightarrow} (L_\beta)_n /(L_\beta^k)_n.$$
It follows that 
$$(L^k)_n = \rho_\beta^{-1}(L_\beta^k)_n \supset (J_\beta)_n.$$
Choose $\gamma \in {\mathcal J}$ so that 
$$J_\gamma \subset \cap_{n=1}^k J_{\beta (k,n)}.$$
Then $$(L^k)_{\leq n} \supset (J_\gamma)_{\leq n},$$
and so $$I_k\supset J_\gamma.$$ \hfill$\square$

 \vspace{3mm}\noindent {\bf Corollary.} Let $L$ be a complete enriched Lie algebra and $E\subset L$ be a subalgebra such that $E/E^2$ has finite type. Then $E$ is closed if and only if $E$ is pronilpotent.
 
 \vspace{3mm}\noindent {\sl proof.} By Proposition \ref{p4}, the structure systems $\rho_\alpha E$ and $E/E^n$ are equivalent. This implies the result.\hfill$\square$

\vspace{3mm}\noindent {\bf Example.} Suppose a graded Lie algebra $E$ satisfies
$$\cap_k E^k= 0\,, \hspace{3mm}\mbox{dim}\, E/E^2<\infty, \hspace{3mm}\mbox{and dim}\, E= \infty.$$
Let $z\in \overline{E}$ satisfy $E\oplus \mathbb Q z \subset \overline{E}$ and let $F$ be the sub Lie algebra generated by $E$ and $\mathbb Q z$. Then $E\subset F\subset \overline{E}$, and it follows that $\overline{E}= \overline{F}$. However the dimension of $F/F^2$ depends on the choice of $z$. 

For instance, suppose $E$ is the completion of the free Lie algebra $\mathbb L(x,y)$ generated by $x,y$, and with the unique enriched structure described in Proposition \ref{p4}. Let $F$ be the sub Lie algebra generated by $x$, $y$ and $z= e^{ad\, x}(y)$. Denote by $t$ the operator $\mbox{ad}_x$. The series $t$ and $e^t$ are algebraically independent in $\mathbb Q[[t]]$, i.e., if $P(a,b)$ is a polynomial in 2 variables $a$ and $b$ and $P(t,e^t)= 0$ then $P$ is identically zero. It follows that if we have a polynomial relation of the form
 $$\sum_n \alpha_n \,\mbox{ad}_x^n (y) + \sum_n \beta_n \,\mbox{ad}^n_x\, e^{\mbox{ad}_x}(y) = 0$$
 then $(\sum \alpha_nt^n + \sum\beta_n t^n e^t) (y)= 0$, and $\alpha_n=\beta_n= 0$ for all $n$.   Therefore the classes of $x,y$ and $z$ are linearly independent in $F/F^2$, and $F/F^3$ has   dimension 3.  
 
 On the other hand, with the same $E$, if $z= \sum_{n\geq 0} \mbox{ad}^n_x(y)$, then $z-y= [x,z]$ and $F/F^2 = E/E^2$.

\section{The quadratic model of $(L, {\mathcal I})$}

First recall that a \emph{quadratic Sullivan algebra} is a Sullivan algebra $(\land V,d)$ in which $d : V\to \land^2V$. This endows $H(\land V)$ with the decomposition
$$H(\land V) = \oplus_k H^{[k]}(\land V)$$
in which $H^{[k]}(\land V)$ is the subspace represented by $\land^kV\cap \mbox{ker}\, d$. Moreover, $V = \varinjlim_n V_n$ where $V_0 = V\cap \mbox{ker}\, d$ and $V_{n+1} = V\cap d^{-1}(\land^2V_n)$. In particular a morphism $\varphi : \land V\to \land W$ of quadratic Sullivan algebras restricts to linear maps 
$$\varphi_n : V_n \to W_n.$$

On the other hand, recall that
the classical functor $L \mapsto C_*(L)$ from graded Lie algebras to cocommutative chain coalgebras is given by $C_*(L)= \land sL$ with differential determined by the condition
$$\partial  (sx\land sy) = (-1)^{deg\, x+1} s[x,y].$$
The dual, $C^*(L)= \mbox{Hom}(C_*(L), \mathbb Q)$ is, with the differential forms on a manifold, one of the earliest examples of a commutative differential graded algebra (cdga). In particular, if dim$\, L<\infty$ then
$
C^*(L) = \land (sL)^\vee,
$ and the differential is determined by the condition
\begin{eqnarray}\label{i5}
<dv, sx,sy> = (-1)^{deg\, y +1} <v, s[x,y]>, \hspace{1cm} v\in L^\vee.
\end{eqnarray}

Now suppose $(L, {\mathcal I})$ is an enriched graded Lie algebra. Then $L= \varprojlim_{\alpha \in {\mathcal I}} L_\alpha$, and since $L_\alpha$ is finite dimensional, 
 $C^*(L_\alpha)= \land (sL_\alpha)^\vee$. Set $V_\alpha = (sL_\alpha)^\vee$ and note that then $sL_\alpha = V_\alpha^\vee$. In particular,
 $$(\land V, d) := \varinjlim_\alpha \land V_\alpha = \varinjlim_\alpha C^*(L_\alpha) $$
 is a cdga, and $d: V\to \land^2V$. Moreover, because each $L_\alpha$ is nilpotent, each $\land V_\alpha$ is a (quadratic) Sullivan algebra, and therefore  $(\land V,d)$ is a quadratic Sullivan algebra. In view of (\ref{i5}), which then holds for all $v\in V$ and $x,y\in L$, \emph{$L$ is the homotopy Lie algebra of $\land V$}, as defined for example in \cite[\S 2.1]{RHTII}.

\vspace{3mm}\noindent {\bf Definition.} $\land V$ is the quadratic model of $(L, {\mathcal I})$.

\vspace{3mm}\noindent {\bf Remark.} Note that the natural morphism $\land V= \varinjlim C^*(L_\alpha) \to C^*(L)$ will not, in general, be a quasi-isomorphism. For instance, if $V= V^1$ has zero differential and countably infinite dimension then $\land V$ is countable but $L$ and $C^*(L)$ are uncountable, and since $L$ is abelian, $C^*(L)= H(C^*(L))$.

\vspace{3mm} With the notation at the start of the Introduction note that a morphism of complete enriched Lie algebras, $\psi : E\to L$, induces morphisms
$$E = \varprojlim_\beta E_\beta \longrightarrow \varprojlim_\alpha E_{\beta (\alpha)} \longrightarrow \varprojlim_\alpha L_\alpha = L.$$
The morphism $E_{\beta (\alpha)}\to E_\alpha$ dualizes to a morphism of quadratic models, and the composites define a unique morphism
$$\varphi : \land V_L\to \land V_E$$
such that $\psi = (\varphi\vert_{V_L})^\vee$. This identifies the correspondence $(L, {\mathcal I}) \mapsto \land V_L$ as a functor from the category of complete enriched Lie algebras to the category ${\mathcal C}'$ of quadratic Sullivan algebras in which the morphisms
$\land V\to \land W$ map $V$ into $W$.

More generally, and in the same way, a coherent linear map $f: E\to L$ induces a linear map $\varphi : V_L\to V_E$ for which $\varphi^\vee = f$.

\begin{Prop}
\label{p5} The correspondence $(L, {\mathcal I})\leadsto \land V$ is a contravariant isomorphism from the category   of complete enriched Lie algebras to the category ${\mathcal C}'$ of quadratic Sullivan algebras.
\end{Prop}

\vspace{3mm}\noindent {\sl proof.}  Each quadratic Sullivan algebra, $\land V$, determines an enriched Lie algebra $(L, {\mathcal I})$. Namely, $L$ is the homotopy Lie algebra of $\land V$, given by
$$L= (sV)^\vee \hspace{5mm}\mbox{and } <v,s[x,y]> = (-1)^{deg\, y+1} <dv, sx\land sy>.$$
Moreover $V= \varinjlim_\beta V_\beta$ where the $V_\beta$ are the finite dimensional subspaces for which $\land V_\beta$ is preserved by $d$. The inclusions $\land V_\beta \to \land V$ dualize to surjections $\rho_\beta : L\to L_\beta$ between the homotopy Lie algebras, and these induce an isomorphism
$$L\stackrel{\cong}{\longrightarrow} \varprojlim_\beta L_\beta.$$
This then endows $L$ with the enriched structure $ (L, {\mathcal I})$ with ${\mathcal I}= \{\mbox{ker}\, \rho_\beta\}$. If $\varphi: \land V\to \land W$ is a morphism then $(\varphi\vert_V)^\vee$ is the corresponding morphism $L_W\to L_V$. \hfill$\square$

 \vspace{3mm}\noindent {\bf Remark.} Proposition \ref{p5} couples the categories of complete enriched Lie algebras, $L$, and quadratic Sullivan algebras, $\land V$, as pairs $(L, \land V)$ where $\land V$ is the quadratic Sullivan model of $L$ and, equivalently, $L$ is the homotopy Lie algebra of $\land V$.

\vspace{3mm}\noindent {\bf Example.}  \emph{Inverse limits}

Suppose $\{L(\sigma), \varphi_{\sigma, \tau}: L(\tau)\to L(\sigma)\}$ is an inverse system of morphisms of complete enriched Lie algebras. Then
$$L:= \varprojlim_\sigma L(\sigma)$$
is naturally a complete enriched Lie algebra.

In fact under the correspondence of Proposition \ref{p5}, this inverse system is the dual of the directed system $\{\land V(\sigma)\}$ of the quadratic Sullivan models of the $L(\sigma)$. But then $\land V := \varinjlim_\sigma \land V(\sigma)$ is a quadratic Sullivan algebra, and $L$ is the homotopy lie algebra of $\land V$. It is immediate that $L$ is the inverse limit in the category of complete enriched Lie algebras.

\vspace{3mm}  The correspondence of Proposition \ref{p5} is reflected in the next Lemma.

\begin{lem}
\label{l5} The degree 1 identification $L\stackrel{\cong}{\to} V^\vee$ induces isomorphisms
$$L/L^{(n+2)} \stackrel{\cong}{\longrightarrow} V_n^\vee, \hspace{1cm} n\geq 0,$$
and therefore identifies $L^{(n+2)} = \{x\in L\, \vert\, <V_p, sx>= 0\}.$\end{lem}

\vspace{3mm}\noindent {\sl proof.}  This is straightforward when $L$ and $V$ are replaced by $L_\alpha$ and $V_\alpha$. In the general case by Lemma \ref{l1} we have the following sequence of isomorphisms
$$
\renewcommand{\arraystretch}{1.6}
\begin{array}{ll} L/L^{(n+2)} &= \overline{L}/\overline{L^{n+2}}= \displaystyle\varprojlim_\alpha \rho_\alpha (L)/ \rho_\alpha (L^{n+2})= \displaystyle\varprojlim_\alpha L_\alpha / L_\alpha^{n+2}\\
& = \displaystyle \varprojlim_\alpha  {(V_\alpha)_n}^\vee = \left( \displaystyle\varinjlim_\alpha  (V_\alpha)_n\right)^\vee = V_n^\vee.
\end{array}
\renewcommand{\arraystretch}{1}
$$
\hfill$\square$

\vspace{3mm} \noindent {\bf Remark.}  A morphism $\psi : E\to L$ of complete enriched Lie algebras induces morphisms $\psi(n) : E/E^{(n)}\to L/L^{(n)}$, $n\geq 1$. If $\varphi : \land V_L\to \land V_E$ is the corresponding morphism of quadratic models then $\psi(n+2)$ is dual to to the linear map $\varphi_n : (V_L)_n \to (V_E)_n$.

\vspace{3mm} Finally, note that a surjection $\xi : V\to P$ dualizes to an inclusion $\xi^\vee : (sP)^\vee \to L$.

\begin{lem}
\label{l6}\begin{enumerate} 
\item[(i)] If $S\subset L$,  let $\xi : V\to P= V/K$ be the surjection given by
$$K = \{ v\in V\, \vert\, <v, sS>= 0\}.$$ Then 
$$\overline{S}= \mbox{Im}\, \xi^\vee.$$
\item[(ii)] In particular, the closed subspaces of $L$ are precisely the subspaces of the form $\mbox{Im}\, \xi^\vee$ as $\xi$ runs through the surjections $V\to P$. 
 
\item[(iii)] Any closed subspace $S\subset L$ has a closed direct summand.
\end{enumerate}
\end{lem}

\vspace{3mm}\noindent {\sl proof.} (i) Use $\xi^\vee$ to identify $(sP)^\vee$ with a subspace of $L$. Since $L = \varprojlim_\alpha L_\alpha$ and $L_\alpha = V_\alpha^\vee$, for $x\in L$ and $v\in V_\alpha$  we have
$$<v, sx> = <v, s\rho_\alpha x>.$$
 Thus
$$V_\alpha \cap K = \{v\in V_\alpha \, \vert \, <v,s\rho_\alpha S>= 0\}.$$
Set $P_\alpha = V_\alpha /V_\alpha \cap K$. Then, because $V_\alpha$ and $L_\alpha$ are finite dimensional, this gives
$$\rho_\alpha S = \{x\in L_\alpha \, \vert\, <V_\alpha \cap K, sx> = 0\} = (sP_\alpha)^\vee.$$
Thus 
$$\overline{S}= \varprojlim_\alpha \rho_\alpha S = \varprojlim_\alpha (sP_\alpha)^\vee = (\varinjlim_\alpha sP_\alpha)^\vee = (sP)^\vee.$$

(ii) If follows exactly as in (i) that each $(sP)^\vee \subset L$ is closed.  

(iii) The inclusion $S\to L$ is the dual of a surjection $V\to V/K$. Dualizing the inclusion $K\to V$ provides a surjection $L \to (sK)^\vee$ onto a closed direct summand of $S$. 
\hfill$\square$

\subsection{Closed subalgebras and ideals}
 
\vspace{1mm} Suppose $E\subset L$ is a closed subalgebra of a complete enriched Lie algebra. Then the inclusion is the dual of a surjection $\xymatrix{\land V \ar[r]^\rho & \land Z}$ of the respective quadratic models, and by Lemma \ref{l6}, 
 $$E = \{x\in L\,\vert\, <\mbox{ker}\, \rho_{\vert V}, sx>= 0\}.$$
 Moreover (\cite[Lemma 10.4]{RHTII}), $d: \mbox{ker}\, \rho_{\vert V}\to \mbox{ker}\, \rho_{\vert V}\, {\scriptstyle \land}\, V$. Recall also that every such surjection induces an inclusion of a closed sub algebra of $L$.
 
 Now suppose $I\subset L$ is a closed ideal. Let $\rho : \land V\to \land Z$ be the corresponding surjection, and denote $W = \mbox{ker}\, \rho_{\vert V}$. In this case (\cite[Lemma 10.4]{RHTII}) the condition that $I$ be an ideal is equivalent to the condition $d: W\to \land^2W$. Thus $\land W$ is a sub quadratic Sullivan algebra and $I$ decomposes $\land V$ as the Sullivan extension
 $$\xymatrix{\land W\ar[r]^-\lambda & \land W\otimes \land Z = \land V \ar[r]^-\rho & \land V\otimes_{\land W}\mathbb Q = \land Z,}$$
 which dualizes to
 $$L/I \leftarrow L\leftarrow I.$$
 In particular this identifies $L/I$ as the homotopy Lie algebra of $\land W$.

  \begin{Prop}
\label{p1n} 
Let $f : L\to E$ be a morphism of complete  enriched Lie algebras. If $L/L^{(2)}\to E/E^{(2)}$ is surjective, then $f$ is surjective.
\end{Prop}

\vspace{3mm}\noindent {\sl proof.}  Denote by $\varphi : \land V_E\to \land V_L$ the corresponding morphism given by Proposition \ref{p5}. By Lemma \ref{l5}, $(V_E\cap \mbox{ker}\, d)^\vee \cong L/L^{(2)}$ and $(V_L\cap \mbox{ker}\, d)^\vee \cong E/E^{(2)}$. It follows that $\varphi : V_E\cap \mbox{ker}\, d\to V_L\cap \mbox{ker}\, d$ is injective. We suppose that for some integer $n$, $\varphi : (V_E)_n\to (V_L)_n$ is injective. Then $\varphi : \land (V_E)_n\to \land (V_L)_n$ is injective. Suppose that $v\in (V_E)_{n+1}$ is in the kernel of $\varphi$, then $\varphi (dv)= 0$, and so $dv= 0$. Therefore $v\in (V_E)_0$ which implies that $v=0$. It follows that the dual map $f$ is surjective.

\hfill$\square$

\begin{Prop}
\label{p4n}
If $f : L\to L'$ is a morphism of complete enriched Lie algebras, then Ker$\, f$ is a closed ideal in $L$ and $Im\, f$ is a closed subalgebra of $L'$.
\end{Prop}

\vspace{3mm}\noindent {\sl proof.} Let $\varphi : \land V'\to \land V$ be the morphism associated to $f$. We denote by $E$ the ideal generated by $V'\cap \mbox{Ker}\, \varphi$ in $\land V'$ and we denote by $\land W = \land V'/E$ the quotient quadratic Sullivan algebra. Then $\varphi$ factorizes through an injection $\varphi' : \land W\to \land V$. We extend then $\varphi'$  to an isomorphism
$\psi : \land W\otimes \land T\stackrel{\cong}{\to} \land V$ from a $\Lambda$-extension. Then $\land W$ is the quadratic model of Im$\, f$ and the quotient cdga $(\land T, \overline{d})$ is the quadratic model of Ker$\, f$. 

The surjection $\land V'\to \land W$ shows that Im$\, f$ is a closed subalgebra of $L'$, and the morphism $\land W\otimes \land T\to \land T$ that $I$ is a closed ideal in $L$. \hfill$\square$

\vspace{3mm}\noindent {\bf Corollary.} 
  Let $f : L\to W$ a surjective coherent morphism between complete enriched abelian Lie algebras. Then $f$ admits a coherent section, $\sigma$.

\vspace{3mm}\noindent {\sl proof.}  By Proposition \ref{p4n}, Ker$\, f$ is a closed subspace. Now by Lemma \ref{l6}(iii), Ker$\, f\subset L$ admits a closed direct summand $S$. It follows that $f_{\vert S} : S\to $ Im$\, f=W$ is an isomorphism, and we define $\sigma = (f_{\vert S})^{-1}$. \hfill$\square$

 \vspace{3mm}\noindent {\bf Example. Weighted Lie algebras.}  Let $E= \oplus_k  E(k)$ be a weighted enriched Lie algebra with defining ideals $I_\alpha = \oplus_k I_\alpha (k)$. The weighting then induces another gradation $V = \oplus_k V(k)$ in the generating space of the corresponding quadratic model:
 $$V(k) = \varinjlim_\alpha s\left( E(k)/I_\alpha(k)\right)^\vee.$$
 It is immediate that the differential preserves the induced gradation in $\land V$.

\vspace{3mm}\noindent {\bf Remark.}   ZFC - set theory alone does not permit us to extend Proposition \ref{p5} to all graded Lie algebras, since   it is consistent with the ZFC-axioms that the same graded Lie algebra can support two inequivalent enriched structures.

In fact, let $L=L_0$ be a vector space and suppose (consistent with the ZFC axioms) that there are isomorphisms 
$$V^\vee \cong sL\cong W^\vee$$
in which card$\,V\neq$ card$\, W$.
Let $\land V$ and $\land W$ be the quadratic Sullivan algebras with zero differential. Then $L$, regarded as an abelian Lie algebra is the homotopy Lie algebra of both $\land V$ and $\land W$, but since $\land V \not\cong \land W$ the corresponding enriched structures in $L$ are not equivalent.

\section{Representations}
\emph{Throughout this entire section $(L, {\mathcal I})$ denotes a fixed complete enriched Lie algebra, with quotient maps $\rho_\alpha: L\to L/I_\alpha:=L_\alpha$  and quadratic model $\land V$.} 

\vspace{3mm}\noindent {\bf Definition.} (i) An \emph{elementary $L$-module} is a finite dimensional nilpotent $L$-module, $Q$, for which some $I_\alpha \cdot Q = 0$.

(ii) An \emph{enriched $L$-module} is an $L$-module, $N$, together with a decomposition,
$$N = \varprojlim_{\tau} N_\tau,$$
of $N$ as an inverse limit of elementary $L$-modules.

(iii) A \emph{coherent $L$-module} is an $L$-module, $M$, together with a decomposition,
$$M=\varinjlim_\sigma M_\sigma,$$
of $M$ as a direct limit of elementary $L$-modules.

(iv) A \emph{morphism} of enriched (resp. coherent) $L$-modules is a morphism of $L$-modules which is the inverse limit (resp. direct limit) of morphisms of elementary $L$-modules. 

\vspace{3mm} Note that an enriched $L$-module $N$ is the dual of the coherent $L$-module $M$ defined by $M= \varinjlim_\sigma N_\sigma^\vee$. 

\vspace{3mm}\noindent {\bf Example.}
Let $L$ be the free Lie algebra on one generator $x$ in degree $0$. The space $M$ of finite sequences of rational numbers $(a_1, \dots , a_n)$ equipped with the $L$-structure defined by $x\cdot (a_1, \dots , a_n)= (a_2, \dots a_{n})$ is a coherent $L$-module. On the other hand the space $N$ of infinite sequences $(a_1, a_2, \dots )$ with the $L$-structure defined by $ (a_1, a_2, \dots)\cdot x = (0, a_1, a_2, \dots )$ is an enriched $L$-module.

\vspace{3mm}\noindent {\bf Remarks.}

1. When dim$\, L/[L,L]<\infty$, any finite dimensional nilpotent $L$-module is elementary. This however is not true in general. In fact, if $L= L_0$ is an infinite dimensional abelian Lie algebra with basis $\{x_i\}$, the subspaces $I(k)$ spanned by the $x_i$, $i\geq k$, define an enriched structure in $L$. Then, consider the finite dimensional nilpotent $L$-module $M = \mathbb Q a \oplus \mathbb Q b$   defined by 
$$x_i\cdot a = b \hspace{3mm}\mbox{and } x_i\cdot b = 0.$$
Since dim$\, L_0=\infty$, $M$ cannnot be an elementary $L$-module.

2. The dual, $M^\vee$, of a coherent $L$-module $M= \varinjlim_\sigma M_\sigma$ inherits the enriched $L$-module structure given by $M^\vee = \varprojlim_\sigma M^\vee_\sigma$. Moreover, each enriched $L$-module is the dual of a unique coherent $L$-module.

3. The dual, $\varphi^\vee$, of a morphism, $\varphi$, of coherent $L$-modules is a morphism of the dual enriched $L$-modules. Moreover, every morphism of the dual modules is the dual of a unique morphism of the corresponding coherent modules.

4. Since for each elementary $L$-module $Q$ there is an ideal $I_\alpha$ with $I_\alpha \cdot Q= 0$, the module $Q$ is naturally a $\widehat{UL_\alpha}$-module for some $\alpha$. Thus the representations of $L$ in enriched and coherent $L$-modules naturally extend to representations of the algebra $\overline{UL}$. 

5. Let $N$ be a right enriched $L$-module. The space  of \emph{decomposable elements} of $N$ is by definition $
\overline{N\cdot L} = \varprojlim_\alpha  N_\alpha \cdot L$.
For instance each complete Lie algebra $L$ is a right enriched module over itself for the adjoint action and $\overline{N\cdot L} = L^{(2)}$. 
 
 \subsection{Augmentations} An elementary $L$-module $N_\tau$ determines the family $I_\alpha\subset {\mathcal I}$ of ideals satisfying $N_\tau \cdot I_\alpha  = 0$. However, the   trivial $L$-modules, $\widetilde{N}_\tau:= N_\tau \otimes_{\widehat{UL_\alpha}} \mathbb Q= N_\tau/N_\tau \cdot L$ are independent of the choice of $\alpha$. Thus $\widetilde{N}= \varprojlim_\tau \widetilde{N_\tau}$ is a trivial enriched $L$-module and the induced morphism
 $$\varepsilon_N : N\to \widetilde{N}$$
 is a morphism of enriched $L$-modules: $\varepsilon_N$ is the \emph{augmentation} for $N$. By construction, $\varepsilon_N : N\to \widetilde{N}$ depends naturally on $N$. 
 
 Let $M$ and $\widetilde{M}$ be the coherent $L$-modules defined by $M^\vee = N$ and $\widetilde{M}^\vee = \widetilde{N}$. Then   $\widetilde{M} = \{m\in M\, \vert\, L\cdot m= 0\}.$

\subsection{Hom and $\otimes$}

Suppose $N$, $M$, and $M'$ are respectively a right enriched $L$-module and two left coherent $L$-modules. Then,
$$\mbox{Hom}(M,N) = \varprojlim_{\sigma, \tau} \mbox{Hom}(M_\sigma, N_\tau)$$
and
$$M\otimes M' = \varinjlim_{\sigma, \sigma'} M_\sigma \otimes M'_{\sigma'}$$
identify $\mbox{Hom}(M,N)$ and $M\otimes M'$ respectively as a right enriched and a left coherent $L$-module.

\subsection{$\overline{UL}$}

The surjections $\overline{UL}\to \widehat{UL_\alpha}/J_\alpha^n$ identify $\overline{UL}$ as a right enriched $L$-module with the representation given by right multiplication.

\subsection{Quadratic $\land V$-modules and holonomy representations}

Quadratic $\land V$-modules are examples of the semifree $\land V$-modules recalled in the Appendix:

\vspace{3mm}\noindent {\bf Definition.} (i) A \emph{quadratic} $\land V$-module is a $\land V$-module of the form $\land V\otimes M$ in which $d: 1\otimes M\to V\otimes M$ and $M= \cup_{k\geq 0} M_k$, where
$$M_0 = M\cap d^{-1}(V) \hspace{3mm}\mbox{and} M_{k+1} = M \cap d^{-1}(V\otimes M_k).$$

(ii) The \emph{holonomy representation} for a quadratic $\land V$-module $\land V\otimes M$ is the representation of $L$ in $M$ given by
$$x\cdot m = -\sum <v_i, sx> m_i,$$
where $d(1\otimes m) = \sum v_i\otimes m_i.$

\vspace{3mm} It is immediate that the holonomy operation makes $M$ into a coherent left $L$-module. In the reverse direction, suppose $M= \varinjlim_\sigma M_\sigma$ is a coherent left $L$-module. The classical Catan-Eilenberg-Serre construction then has the form
$$C^*(L_{\alpha (\sigma)}, M_\sigma) = \land (sL_{\alpha (\sigma)})^\vee \otimes M_\sigma = \land V_{\alpha (\sigma)} \otimes M_\sigma$$
in which $d : M_\sigma \to V_{\alpha (\sigma)}\otimes M_\sigma$. Passing to direct limits constructs the quadratic $\land V$-module $\land V\otimes M = \varinjlim_\sigma \land V_{\alpha (\sigma)}\otimes M_\sigma$.
This establishes a bijection between quadratic $\land V$-modules and coherent left $L$-modules.

Observe as well that since quadratic $\land V$-modules are semifree  a morphism from a quadratic $\land V$-module lifts through any surjective quasi-isomorphism of $\land V$-modules \cite[\S 6]{FHTI}. In particular a quasi-isomorphism $\land V\otimes M \stackrel{\simeq}{\to} \land V\otimes M'$ of quadratic $\land V$-modules has a homotopy inverse. It follows that it induces an isomorphism $M\stackrel{\cong}{\to} M'$ of coherent left $L$-modules.

Finally, as recalled in \S 4, associated with $L$ is the differential coalgebra $(\land sL, \partial)$ in which $\partial (sx\land sy) = (-1)^{deg\, sx} s[x,y]$. More generally, (eg. \cite[Chap 2]{RHTII}) associated with $M$ is the differential $(\land sL, \partial)$-comodule, $(\land sL\otimes M, \partial)$, characterized by
$$\partial (sx\otimes m)= (-1)^{deg\, sx} x\cdot m.$$
As recalled in \cite[Chap 10]{RHTII}, there are natural isomorphisms
$$\mbox{Tor}_p^{UL}(\mathbb Q, M) \cong H_{[p]}(\land sL\otimes M), \hspace{5mm} p\geq 0,$$
where $H_{[p]}$ denotes the subspace of $H(\land sL\otimes M)$ represented by cycles in $\land^psL\otimes M$. 

\subsection{Acyclic closures} 

Recall (\cite[Chap 3]{RHTII}) that the acyclic closure   is the special case, $\land V\otimes \land U$,  of a $\Lambda$-extension constructed inductively as follows. Write $V = \cup_n V_n$ with $V_0= V\cap \mbox{ker}\, d$ and $V_{n+1}= V\cap d^{-1}(\land^2V_n)$. Then $U = \cup_{n\geq 0} U_n$ and there is a degree $1$ isomorphism $p: U\stackrel{\cong}{\to} V$ restricting to isomorphisms $U_n \stackrel{\cong}{\to} V_n$. Thus this identifies $L= U^\vee$. Finally, the differential is determined by the conditions
$$du=pu, \hspace{3mm} u\in U_0\hspace{3mm}\mbox{and } (d-p) : U_{n+1}\to V_n\otimes \land U_n.$$
In particular the augmentation $\varepsilon_U : \land U\to \mathbb Q$, $U\to 0$ together with the unique augmentation in $\land V$ define a quasi-isomorphism,
$$\land V\otimes \land U \stackrel{\simeq}{\longrightarrow} \mathbb Q.$$

This construction identifies $\land V\otimes \land U$ as a quadratic $\land V$-module. In this case the holonomy representation is a representation $\theta : L\to \mbox{Der}(\land U)$ of $L$ by derivations in $\land U$, and if $x\in L$ and $u\in U_n$ then 
\begin{eqnarray}\label{i6}
\theta (x)u \,\,+ <pu, x>\in V_{n-1}\otimes \land U_{n-1}.
\end{eqnarray}

Now   denote by $\eta_L$ the morphism of enriched $L$-modules,
$$\eta_L : \overline{UL}\stackrel{ }{\longrightarrow} (\land U)^\vee,$$
defined by $$\eta_L(a)(\Phi) = \varepsilon_U(a\cdot \Phi), \hspace{5mm} a\in \overline{UL}, \Phi\in \land U.$$
By \cite[Theorem 6.1]{RHTII}, $\eta_\alpha : \widehat{UL_\alpha}\to (\land U_\alpha)^\vee$ is an isomorphism, and since $\eta_L$ is a morphism of inverse limits, we have

\begin{Prop}
\label{eta}
The morphism $$\eta_L: \overline{UL}\stackrel{\cong}{\longrightarrow} (\land U)^\vee $$
is an isomorphism of $\overline{UL}$-modules. 
\end{Prop}

\subsection{Adjoint representations}

The right adjoint representation of $L_\alpha$ in $L_\alpha$ makes $L_\alpha$ into an $L_\alpha$-module. This exhibits $L= \varprojlim_\alpha L_\alpha$ with the right adjoint representation of $L$ as a right enriched $L$-module. Since $sL_\alpha = V_\alpha^\vee$ the corresponding  left coherent $L$-module is a representation of $\overline{UL}$ in $V$. It is given explicitly by
$$<x\cdot v, sy> = -<dv, sx,sy>, \hspace{1cm} v\in V, x,x\in L.$$

\begin{lem}
\label{l7} If $S\subset L$ is any subspace then $S\cdot \overline{UL}$ is an ideal. If $S$ is a graded space of finite type, then $S\cdot \overline{UL}$ is closed.\end{lem}

 \vspace{3mm}\noindent {\sl proof.} Observe that
 $$[x\cdot \Phi, y]= x\cdot \Phi y\,, \hspace{1cm} x,y\in L,\, \Phi\in \overline{UL},$$
 and so $S\cdot \overline{UL}$ is an ideal. Since a graded space is closed if and only if each component in a given degree is closed, if $S$ has finite type it is sufficient to show that each $S^k\cdot \overline{UL}$ is closed. Since a finite sum of closed subspaces is closed we have only to show that $x\cdot \overline{UL}$ is closed for each $x\in L$.
 
 Denote $\rho_\alpha x$ by $x_\alpha$. Then the closure of $x\cdot \overline{UL}$ is the inverse limit
 $$\varprojlim_\alpha \rho_\alpha (x\cdot \overline{UL}) = \varprojlim_\alpha x_\alpha \cdot \widehat{UL_\alpha} = \mathbb Q x_\alpha \oplus \varprojlim_\alpha x_\alpha \cdot \widehat{J_\alpha} = \mathbb Q x_\alpha \oplus \varprojlim_{n,\alpha} (x_\alpha \cdot J_\alpha)/(x_\alpha\cdot J_\alpha^n).$$
 On the other hand, because each $J_\alpha/J_\alpha^n$ is finite dimensional, the surjections
 $J_\alpha /J_\alpha^n \to (x_\alpha \cdot J_\alpha)/x_\alpha \cdot J_\alpha^n)$ induce a surjection
 $$\overline{J} \to \varprojlim_{n,\alpha} x_\alpha\cdot J_\alpha /(x_\alpha \cdot J_\alpha^n).$$
 This factors through the map $\overline{J}\to x\cdot \overline{J}$, and therefore shows that
 $$x\cdot \overline{J} \to \varprojlim_\alpha \rho_\alpha (x\cdot \overline{J})$$
 is surjective. But this is the inclusion of $x\cdot \overline{J}$ in its closure, and so $x\cdot \overline{J}$ is closed.

 \hfill$\square$

\subsection{Profree $L$-modules}

\vspace{3mm}\noindent {\bf Definition.} A \emph{profree} $L$-module is an enriched $L$-module of the form Hom($S, \overline{UL}$) where $S$ is a trivial coherent $L$-module.

\vspace{3mm}\noindent {\bf Remarks.}

1. It is immediate from the definition that the augmentation in $\mbox{Hom}(S, \overline{UL})$ is the morphism 
$$\mbox{Hom}(S, \overline{UL}) \to S^\vee$$
given by $\gamma \mapsto \varepsilon\circ \gamma$, where $\varepsilon$ is the augmentation $\overline{UL}\to\mathbb Q$ with kernel $\overline{J}$.

2. If $S$ is a graded vector space of finite type and $S = S^{\geq 0}$ then it is immediate that
$$\mbox{Hom}(S, \overline{UL}) = S^\vee \otimes \overline{UL}$$
is a free $\overline{UL}$-module.

\vspace{3mm}\emph{Now fix a left coherent $L$-module, $M$, } and define a left coherent $L$-module $S$ by setting 
 $$S = \{a\in M\,\vert\, L\cdot a = 0\}.$$
 It follows from    $\S$ 5.1 that $M^\vee \to S^\vee$ is the augmentation for $S^\vee$. 
 
 \vspace{3mm} As usual we denote by $\land V$ the quadratic Sullivan model of $L$ and  by $(\land V\otimes M,d)$ the holonomy representation defined in \S 5.4. We remark that by definition, $d(1\otimes S)= 0$.

\begin{Prop}
\label{p6} With the hypotheses and notation above the following conditions are equivalent:
\begin{enumerate}
\item[(i)] The inclusion $S\to 1\otimes S\subset \land V\otimes M$ induces an isomorphism
$$S\stackrel{\cong}{\longrightarrow} H(\land V\otimes M).$$
\item[(ii)] $M\cong \land U\otimes S$ as left coherent  $L$-modules.
\item[(iii)] $M^\vee $ is a profree $L$-module with augmentation $M^\vee \to S^\vee$.
\end{enumerate}
\end{Prop}

\vspace{3mm}\noindent {\sl proof.} Choose any linear retraction $\rho : M\to S$. Tensored with the augmentation of $\land V$ this provides a morphism $\land V\otimes M\to S$ of $\land V$-modules. Lifting this through the quasi-isomorphism $\land V\otimes \land U\otimes S\to S$ provided by the augmentation of $\land V\otimes \land U$, we obtain a morphism,
$$\psi : \land V\otimes M\to \land V\otimes \land U\otimes S,$$
of $\land V$-modules. Then division by $\land^{\geq 1}V$ yields a morphism $\varphi : M\to \land U\otimes S$ of coherent left $L$-modules.

To complete the proof we show first that (i) $\Leftrightarrow$ (ii). In fact it is immediate that (ii) $\Rightarrow $ (i). On the other hand, if (i) holds then $\land V\otimes M\to S$ is a quasi-isomorphism, and hence  so is $\psi$. It follows that $\varphi$ is an isomorphism.

Finally, to show (ii) $\Leftrightarrow$ (iii) note that since $(\land U)^\vee = \overline{UL}$, (ii) $\Rightarrow$ (iii) is immediate from the definition of profree. In the reverse direction, suppose $M^\vee  $ is profree with augmentation $M^\vee \to S^\vee$. Then, since augmentations are preserved by morphisms, there is a commutative diagram
$$\xymatrix{
M^\vee \ar[rr]^\chi_\cong\ar[rd] && \mbox{Hom}(S, \overline{UL})\ar[ld]\\
&S^\vee}$$
of right enriched $L$-modules.  

On the other hand, the dual of the morphism $\varphi$ is a morphism $\varphi^\vee: (\land U\otimes S)^\vee \to M^\vee $ of   enriched $L$-modules. Since the representation of $L$ in $\land U\otimes S$ is just the holonomy representation in $\land U$ and since $(\land U)^\vee\cong \overline{UL}$, $\varphi^\vee$ is a morphism from Hom$(S, \overline{UL})$. Thus the diagram
$$
\xymatrix{
\mbox{Hom}(S,\overline{UL}) \ar[rr]^{\varphi^\vee}\ar[rrd] && M^\vee \ar[d] \ar[rr]^\chi_\cong && \mbox{Hom}(S, \overline{UL})\ar[lld]\\
&& S^\vee
}$$
commutes.

To establish (ii) we will show that $\varphi$ is an isomorphism, and of course it is sufficient to show that $\varphi^\vee$ is an isomorphism or, equivalently, that $\gamma := \chi \circ \varphi^\vee$ is an isomorphism. But the commutative diagram above shows that
$$\gamma - id : \mbox{Hom}(S, \overline{UL})\to \mbox{Hom}(S, \overline{J}).$$
Since $\gamma $ is a morphism of enriched $L$-modules it follows that
$$(\gamma - id): \mbox{Hom}(S, \overline{J}^{(n)}) \to \mbox{Hom}(S, \overline{J}^{(n+1)}).$$
Finally, because $\overline{UL}= \varprojlim_n \overline{UL}/\overline{J}^{(n)}$, it follows that
$$\mbox{Hom}(S, \overline{UL}) = \varprojlim_n \mbox{Hom}(S, \overline{UL}/\overline{J}^{(n)})$$
and therefore $\gamma$ is an isomorphism.
\hfill$\square$

\subsection{Closed ideals}

Suppose $I\subset L$ is a closed ideal and denote by  
$$\xymatrix{\land W \ar[r]^-\lambda & \land W\otimes \land Z = \land V \ar[r]^-\rho & \land V\otimes_{\land W}\mathbb Q = \land Z,}$$
the $\Lambda$-extension in which $\lambda$ and $\rho$ dualize to the short exact sequence
 $L/I \leftarrow L\leftarrow I.$ 

Next let $(\land W\otimes \land U_W, d_W)$ and $(\land Z\otimes \land U_Z, d_Z)$ be the respective acyclic closures. Since $\land W\otimes \land U_W$ and $\land W\otimes \land Z$ are $\Lambda$-extensions we can form the $\Lambda$-extension 
$$\land W \otimes \land U_W\otimes \land Z := (\land W\otimes \land U_W)\otimes_{\land W} (\land W\otimes \land Z).$$
Now it follows from \cite[Proposition 1]{Sdepth} that the acyclic closure of $\land V$ is a $\Lambda$-extension of $\land W\otimes \land U_W\otimes \land Z$ of the form
$$(\land W\otimes \land U_W\otimes \land Z,d) \to (\land W\otimes \land U_W\otimes \land Z\otimes \land U_Z,d)$$
in which
\begin{eqnarray}
\label{i7} (d-d_Z): \land U_Z \to W\land (\land W\otimes \land U_W\otimes \land Z\otimes \land U_Z).  \end{eqnarray}
Since $\Lambda$-extensions are semifree modules with respect to the underlying quadratic Sullivan algebras it follows that $\land U_Z$ is the union of subspaces
$M(0) \subset \dots \subset M(k)\subset M(k+1) \subset \dots $ satisfying
\begin{eqnarray}
\label{i8} d  : M(0)\to \land W\otimes \land U_W\otimes \land Z\hspace{3mm}\mbox{and } d  : M(k+1)\to (\land W\otimes \land U_W\otimes \land Z)\otimes M(k). \end{eqnarray}

On the other hand,    since $[I,L]\subset I$, $I$ is an $L$-submodule with respect to the right adjoint representation. Moreover, this representation factors to yield a representation of $L/I$ in $I/I^2$. 

\begin{Prop}
\label{p7} With the hypotheses above:
\begin{enumerate}
\item[(i)] $I$ is an enriched sub-module of $L$ for the adjoint representation.
\item[(ii)] The adjoint representation of $L/I$ in $I/I^2$ factors over the surjection $I/I^2\to I/I^{(2)}$ to identify $I/I^{(2)}$ as an enriched $L/I$-module.   
\end{enumerate}
\end{Prop}

\vspace{3mm}\noindent {\sl proof.} (i). Each $\rho_\alpha I \subset L_\alpha$ is an $L_\alpha$-module: $[\rho_\alpha I, L_\alpha]\subset \rho_\alpha I$. This exhibits $I= \varprojlim_\alpha \rho_\alpha I$ as an enriched $L$-sub module of $L$.

(ii). By definition, $I^{(2)} = \varprojlim_\alpha [\rho_\alpha I, \rho_\alpha I]$ is also a closed ideal. Now the right adjoint representation of $L_\alpha$ factors to give a representation of $L_\alpha/\rho_\alpha I$ in $\rho_\alpha I/[\rho_\alpha I, \rho_\alpha I]$. But by \S 2,
$$\varprojlim_\alpha L_\alpha /\rho_\alpha I = L/I\hspace{5mm}\mbox{and } \varprojlim_\alpha \rho_\alpha/[\rho_\alpha I, \rho_\alpha I]= I/I^{(2)}.$$

\hfill$\square$

\vspace{3mm} Now, since $\land V = \land W\otimes \land Z$ is a $\Lambda$-extension of $\land W$ it is a semifree $\land W$-module. In particular from the equation $d^2= 0$ n $\land V$ it follows \cite[Chap. 4]{RHTII} that the holonomy representation is a representation of $L_W$ in $(\land Z, d_Z)$.

Denote by $\overline{\theta}$ the induced representation in $H(\land Z, d_Z)$, and note that it follows immediately from the construction that $\overline{\theta}$ restricts to a representation in the subspaces $H^{[p]}(\land Z)$ represented by the cycles in $\land^pZ$.
 Moreover, filtering $\land W\otimes \land Z$ by the ideals $\land^{\geq r}W\otimes \land Z$   induces a spectral sequence whose $E_1$-term $\land W\otimes H(\land Z)$ is a semifree $\land W$-module. It is immediate from the construction that $\overline{\theta}$ is the resulting holonomy representation.
 
 \hfill$\square$

\begin{Prop}
\label{p8} Let $I$ be a closed ideal in a profree Lie algebra $L$ with associated quadratic algebra $\land Z$. Then  there is a natural isomorphism
$$H^{[1]}(\land Z)^\vee \cong s\, I/I^{(2)}$$
of enriched $L/I$-modules, where $L/I$ acts in $H^{[1]}(\land Z)^\vee$   by the dual of the holonomy representation and in $s\, I/I^{(2)}$  by the right adjoint representation. \end{Prop}

\vspace{3mm}\noindent {\sl proof.} First observe that Lemma \ref{l5} provides an isomorphism
$$H^{[1]}(\land Z)^\vee = (Z\cap \mbox{ker}\, \overline{d})^\vee \cong s\, I/I^{(2)}.$$
It remains to show that this is an isomorphism of enriched $L/I$-modules. A limit argument reduces this to the case dim$\, V<\infty$, where it is sufficient to show that it is an isomorphism of $L$-modules.

 Let $x\in I$, $y\in L$ and $z\in \mbox{ker}\, \overline{d}$. From (\ref{i5}) we have
$$<d(1\otimes z), sx,sy> = (-1)^{deg\, y+1} <z, s[x,y]>= (-1)^{deg\, x} <s[x,y], z>.$$
On the other hand, since $<V, sx>= 0$, we   have from (\ref{i6}) that
$$\renewcommand{\arraystretch}{1.3}
\begin{array}{ll}
<d(1\otimes z), sx,sy> & = \sum <v_i, sy>\, <z_i, sx>\\
& = - <\overline{\theta}(y)z, sx> = (-1)^{deg\, x} <sx, \overline{\theta}(y)z>. 
\end{array}
\renewcommand{\arraystretch}{1}$$
 \hfill$\square$

\newpage
\part{Profree Lie algebras}

Free graded Lie algebras freely generated by a space $V$, and which we denote by  $\mathbb L_V$, play a key role in Lie algebra theory. Profree Lie algebras are the analogue of free Lie algebras in the category of complete enriched Lie algebras. More precisely, let $L$ be a complete enriched Lie algebra, and recall from Lemma \ref{l6} that $L$ admits direct decompositions $L= L^{(2)}\oplus T$ in which $T$ is a closed subspace. Moreover, as we shall show in Proposition \ref{p9}, if $L$ is profree then $T$ freely generates a free sub  Lie algebra of $L$.

\vspace{3mm}\noindent {\bf Definition} \begin{enumerate}
\item[1.] The complete enriched Lie algebra, $L$, is \emph{profree} if, for some decomposition $L= L^{(2)}\oplus T$ with $T$ closed, any coherent linear map
$$f : T\to E$$
into a complete enriched Lie algebra extends to a morphism
$$\varphi : L\to E.$$
\item[2.] A decomposition $L= L^{(2)}\oplus T$ with $T$ closed is \emph{extendable} if it satisfies this condition.
\end{enumerate}

\vspace{3mm}\noindent {\bf Remarks} \begin{enumerate}
\item[1.] If $T$ is a closed direct summand of $L^{(2)}$ in a complete enriched Lie algebra $L$ then (Lemma \ref{l4}) $L$ is the closure of the sub Lie algebra generated by $T$. It follows that $\varphi$ is uniquely determined by $f$.
\item[2.] Suppose $\land V$ is the quadratic model of a complete enriched Lie algebra $L$. As in \S 4 we write
$$V_0:= V\cap \mbox{ker}\, d.$$
If $L= L^{(2)}\oplus T$ with $T$ closed then by Lemma \ref{l5}, division by $L^{(2)}$ induces an isomorphism
$$T\stackrel{\cong}{\longrightarrow} V_0^\vee.$$
\end{enumerate}

\section{Characterization of profree Lie algebras} 

\begin{Theorem}
\label{t1}
A complete enriched Lie algebra, $L$, is profree if and only if its quadratic model $\land V$ satisfies
$$H(\land V)= \mathbb Q \oplus (V\cap \mbox{ker}\, d).$$
In this case each direct decomposition $L= L^{(2)}\oplus T$ with $T$ closed is extendable.
\end{Theorem}

\vspace{3mm}

The proof of Theorem \ref{t1} will be carried out in the next two subsections.  

\subsection{Quadratic Sullivan algebras} 

For any quadratic Sullivan algebra, $\land V$, we write
$$H(\land V) = \oplus_k H^{[k]}(\land V),$$
where $H^{[k]}(\land V)$ is the image in $H(\land V)$ of the space of cycles in $\land^kV$. 
Our objective here is to prove (Lemma \ref{l9}) that
$$H^{[2]}(\land V)= 0 \Longleftrightarrow H(\land V) = \mathbb Q \oplus (V\cap \mbox{ker}\, d).$$

\begin{lem} 
\label{l8} Suppose $\psi: E\to L$ is a morphism of complete enriched Lie algebras.
\begin{enumerate}
\item[(i)] $\psi $ is surjective if and only if $\psi (2): E/E^{(2)}\to L/L^{(2)}$ is surjective.
\item[(ii)] Suppose the quadratic model, $\land V$, of $L$ satisfies $H^{[2]}(\land V)= 0$. Then the following assertions are equivalent:
\begin{enumerate}
\item[(a)] $\psi (2)$ is an isomorphism,
\item[(b)] $\psi(n)$ is an isomorphism, $n\geq 2$,
\item[(c)] $\psi$ is an isomorphism.
\end{enumerate}
\end{enumerate}
\end{lem}

\vspace{3mm}\noindent {\sl proof.} First recall that $\psi$ is the dual of   $\varphi :  V\to  W$ where $\varphi : \land V\to \land W$ is the corresponding morphism between the quadratic models of $L$ and $E$. In view of Lemma \ref{l5}, $\psi(n+2): E/E^{(n+2)}\to L/L^{(n+2)}$ is the dual of $\varphi_n: V_n\to W_n$.  

(i) We have only to show that if $\varphi_0$ is injective then $\varphi$ is injective. Assume by induction that $\varphi_n$ is injective. if $v\in V_{n+1}$ and $\varphi v= 0$ then since $dv\in \land ^2V_n$ and $\varphi_n$ is injective it follows that $dv= 0$. Thus $v\in V_0$ and since $\varphi_0$ is injective $v= 0$.

(ii) Suppose (a) is satisfied and assume by induction that $\varphi_n : V_n \stackrel{\cong}{\to} W_n$. If $w\in W_{n+1}$ then $dw$ is a cycle in $\land^2W_n$ and so $dw= \varphi\Phi$ where $\Phi$ is a cycle in $\land^2V_n$. Because $H^{[2]}(\land V)= 0$, it follows that $\Phi = dv$ for some $v\in V_{n+1}$. Thus $d(w-\varphi v)= 0$ and so $w-\varphi(v)\in W_0$. By hypothesis, $w-\varphi v = \varphi v_0$ for some $v_0\in V_0$, and so $W_{n+1} \subset \varphi (V_{n+1})$. On the other hand, by (i), $\varphi_{n+1} : V_{n+1}\to W_{n+1}$ is injective. This proves that (a) $\Rightarrow$ (b).

But then $\varphi = \varinjlim \varphi_n$ and so $\varphi$ is an isomorphism. Thus (b) $\Rightarrow$ (c). Finally, it is immediate that if $\varphi$ is an isomorphism then so is $\varphi_0$. \hfill$\square$

\begin{lem}
\label{l9}
If a morphism $\rho : \land V\to \land Z$ of quadratic Sullivan algebras restricts to a surjection $\rho_0 : V_0\to Z_0$, and if $H^{[2]}(\land V)= 0$, then
$$H^{[k]}(\land Z)= 0, \hspace{3mm} k\geq 2.$$
\end{lem}

\vspace{3mm}\noindent {\sl proof.} The proof is in two Steps.

\vspace{2mm}\noindent \emph{Step One. If, in addition, $\rho$ is surjective then $H^{[k]}(\land Z) = 0$, $k\geq 2$.}

Here we first show that $H^{[2]}(\land Z)= 0$. Since  $\rho$ is surjective, by \cite[Cor.3.4]{RHTII} $\rho$ extends to a quasi-isomorphism
$$\rho : \land V\otimes \land U \stackrel{\simeq}{\to} \land Z$$
in which $\land V\otimes \land U$ is a $\Lambda$-extension of $\land V$ and $d: \land U\to V\otimes \land U$. Since $\rho$ is surjective there is a quasi-isomorphism $\sigma : \land Z\to \land V\otimes \land U$ satisfying $\rho\circ \sigma = id$.  Because $d: \land U\to V\otimes \land U$ and $\land V\otimes \land U$ is a $\Lambda$-extension it is straightforward to verify that $\sigma$ may be constructed so that $\sigma: Z\to V\otimes \land U$. In particular, the decomposition $\land V\otimes \land U = \oplus_k \land^kV\otimes \land U$ induces a decomposition of its homology, and
$$\sigma: H^{[k]}(\land Z) \stackrel{\cong}{\to} H^{[k]} (\land V\otimes \land U), \hspace{5mm} k\geq 0.$$

On the other hand, because $\land V\otimes \land U$ is a $\Lambda$-extension, $\land U$ is the increasing union of the subspaces $(\land U)_q$ given by
$$(\land U)_0= \land U \cap d^{-1}(\land V\otimes 1)\hspace{5mm}\mbox{and } (\land U)_{q+1} = \land U \cap d^{-1}(\land V\otimes (\land U)_q).$$
Thus
$$H^{[k]}(\land V\otimes \land U) = \varinjlim_q H^{[k]}(\land V\otimes (\land U)_q).$$
But the differential in the quotients $\land V\otimes ((\land U)_{q+1}/(\land U)_q)$ is just $d\otimes id$, and it follows that $H^{[2]}(\land V\otimes \land U)= 0$.

It remains to show that $H^{[k]}(\land Z)= 0$, $k\geq 3$. Suppose $\Phi\in \land^kZ$ is a cycle. There is then a sequence $z_1, \dots , z_r$ of elements in $Z$ such that $dz_1= 0$, $dz_{i+1}\in \land^2(z_1, \dots , z_i)$, $\Phi\in \land (z_1, \dots , z_r)$ and such that division by $z_1, \dots , z_r$ maps $\Phi$ to zero. When $r=1$, $\Phi = z_1^k$ is a boundary. We use induction on $r$ and on $k$ to show that $\Phi$ is a boundary.

Observe first that division by $z_1$ gives a quadratic Sullivan algebra $(\land Z',d')$. By what we just proved, $H^{[2]}(\land Z',d')= 0$. Thus by induction on $r$, the image of $\Phi$ in $\land Z'$ is a boundary. It follows that for some $\Phi'\in \land^{k-1}Z$,
$$\Phi -d(1\otimes \Phi') = z_1\otimes \Phi'',$$
with $\Phi''\in \land Z$. In particular, in $\land Z'$, $\Phi''$ is a $d'$-cycle. By induction on $k$ in $\land Z'$, $\Phi'' = d'\Omega$ for some $\Omega \in \land ^{k-2}Z'$. Therefore $\Phi''= d\Omega + z_1\otimes \Psi$ for some $\Psi$.

If deg$\, z_1$ is odd then $z_1\otimes \Phi''= d(-z_1\otimes \Omega)$, and so $\Phi$ is a boundary. If deg$\, z_1$ is even then, since $H^{[2]}(\land Z)= 0$, we may choose $z_2$ so $dz_2= z_1^2$. Division by $z_2$ and $z_1^2$ then gives a quasi-isomorphism $\land Z\stackrel{\cong}{\longrightarrow} \left((\land z_1)/z_1^2\right) \otimes \land Z''$ and the same argument as above shows that $\Phi$ is a boundary.

\vspace{2mm}\noindent \emph{Step Two. Completion of the proof of Lemma \ref{l9}.}

  We define a sequence of surjective morphisms
$$\land V= \land V(1) \twoheadrightarrow \land V(2) \twoheadrightarrow \dots \twoheadrightarrow \land V(p)\twoheadrightarrow \dots$$
such that $\rho$ factors through each to yield morphisms
$$\rho(p) : \land V(p)\twoheadrightarrow \land Z.$$
In fact, if $\rho(p)$ has been defined let $\varphi_p: \land V(p)\twoheadrightarrow \land V(p+1)$ be obtained by division by $V(p)_0 \cap \mbox{ker} \rho(p).$

Then the kernels of the surjections $\land V \to \land V(p)$ form an increasing sequence of subspaces $K(p)\subset V$. Set $K = \cup_p K(p)$ and let
$$\varphi : \land V\to \land W$$
be the surjection obtained by dividing $V$ by $K$. Then $\varphi = \varinjlim_p \varphi_p$, and so   $\varphi$ and $\varphi_0$ are surjective. Thus by Step One   $H^{[k]}(\land W)= 0$, $k\geq 2$.

Now by construction,   $\rho$ factors as
$$\land V\stackrel{\varphi}{\to} \land W\stackrel{\gamma}{\to} \land Z.$$
Moreover, $\land W = \varinjlim_p \land V(p)$ and so $W_0 = \varinjlim_p V(p)_0$. Let $I(p)$ be the image of $V(p)_0$ in $V(p+1)_0$. Then also
$$W_0 = \varinjlim_p I(p).$$
Thus by construction, $\gamma   : W_0\stackrel{\cong}{\to} Z_0$. Now Lemma \ref{l8}(ii) asserts that $\gamma$ is an isomorphism, and hence
$$H^{[k]}(\land Z) = H^{[k]}(\land W)= 0, \hspace{3mm} k\geq 2.$$

\hfill$\square$

\subsection{Proof of Theorem \ref{t1}}

This is in two steps.

\vspace{3mm}\noindent {\emph{Step One.} If $L= L^{(2)}\oplus T$ is extendable then $H(\land V)= \mathbb Q \oplus (V\cap \mbox{ker}\, d)$.

\vspace{2mm} In this case, by adjoining 	additional variables construct an inclusion $\lambda : \land V\to \land W$ of quadratic Sullivan algebras for which 
$$V_0\stackrel{\cong}{\longrightarrow} W_0\hspace{5mm}\mbox{and } H^{[2]}(\land W)= 0.$$
Dualizing $V\to W$ gives a surjection $L\stackrel{\rho}{\leftarrow} L_W$ of complete enriched Lie algebras. Moreover, since $V_0\stackrel{\cong}{\to} W_0$, $L_W = L_W^{(2)}\oplus T_W$ where $T_W$ is closed and 
$$\rho_T : T\stackrel{\cong}{\longleftarrow} T_W.$$

Here both $\rho_T$ and its inverse, $\sigma$, are coherent and so $\sigma$ extends to a morphism $\varphi : L\to L_W$. Since $\rho\circ \varphi \vert_T = id\vert_T$ it follows that $\rho\circ \varphi = id_L$. Now dualizing, $\varphi$,  gives a morphism
$$\psi : \land V\leftarrow \land W$$
such that $\psi\circ \lambda = id$.

On the other hand, the inclusion $\lambda$ defines a $\Lambda$-extension
$$\land V\otimes \land Z \stackrel{\cong}{\longrightarrow} \land W,$$
in which, if $Z\neq 0$ then some $z\in Z$ satisfies $dz\in \land^2V$. Thus
$$d\psi z= \psi dz= dz.$$
Hence $\psi z-z$ is a cycle in $W$. This gives
$$\psi z -z\in W_0= V_0$$
and so $z\in V$. This contradicts $z\in Z$ and it follows that $Z=0$ and $\psi$ is an isomorphism inverse to $\lambda$. In particular $H^{[2]}(\land V)= 0$ and now Lemma \ref{l9} implies that $H(\land V) = \mathbb Q \oplus (V\cap \mbox{ker}\, d).$

\vspace{3mm}\noindent {\emph{Step Two.}  If $H(\land V)= \mathbb Q \oplus (V\cap \mbox{ker}\, d)$ then any decomposition $L= L^{(2)}\oplus T$ with $T$ closed is extendable.

\vspace{2mm} Suppose $f : T\to E$ is a coherent linear map into a complete enriched Lie algebra. Since $f$ is coherent by Lemma \ref{l6} it is the dual of a linear map
$$g: V_0\leftarrow W,$$
where $\land W$ is the quadratic model of $E$. Now recall from Lemma \ref{l6} that since $T$ is closed,
$$V = V_0\oplus C$$
with $C = \{v\in V\, \vert\, <v,sT>= 0\}$. We shall construct a morphism
$$\varphi : \land V\leftarrow \land W$$
of quadratic Sullivan algebras such that $\varphi w-gw\in C$. 

In fact, extend $g$ to a morphism $\land V_0\leftarrow \land W_0$ and then recall that $W = \cup_k W_k$ with $W_{k+1}= W\cap d^{-1}(\land^2W_k)$. Suppose that $\varphi$ has been constructed in $\land W_k$. If $w\in W_{k+1}$ then $\varphi dw$ is a cycle in $\land^2V$ and therefore there is a unique $v\in C$ such that $dv= \varphi dw$. Extend $\varphi$ to $W_{k+1}$ by setting 
$$\varphi w= v+ g(w).$$
The uniqueness of $v$ implies that this coincides with the construction of $\varphi$ in $W_k$, and so this inductive process constructs the morphism $\varphi$.

Finally, let $\gamma : L\to L_W$ be the morphism determined by $\varphi$. Then for $x\in T$ and $w\in W$,
$$<w,s(\gamma x)> = <\varphi w, sx> = <v+ gw, sx>.$$
Since $<C, sT>= 0$ and $v\in V$ we have
$$<w, s(\gamma x)> = <gw, sx> = <w, sf(x)>.$$
Thus $\gamma$ extends $f$.

\vspace{2mm} This completes the proof of Step Two, and with it of Theorem \ref{t1}.

\hfill$\square$

\section{Properties of profree Lie algebras}

\vspace{3mm} In analogy with properties of the  free Lie algebras, we have from Theorem \ref{t1}:

\vspace{3mm}\noindent {\bf Corollary 1.} A complete enriched Lie algebra $L$ admits a surjective morphism of enriched Lie algebras
$E \twoheadrightarrow L$ in which $E$ is profree.

\vspace{3mm}\noindent {\sl proof.} Let $\land V$ be the quadratic model of $L$, and adjoin variables to $V$ to obtain an inclusion $\land V\to \land W$ of quadratic Sullivan algebras such that $V_0\stackrel{\cong}{\to} W_0$ and $H^{[2]}(\land W)= 0$. Then dualize to obtain a surjective morphism $E\twoheadrightarrow L$ from the homotopy Lie algebra of $W$, which, by Lemma \ref{l9} and Theorem \ref{t1}, is profree. 
\hfill$\square$

\vspace{3mm}\noindent {\bf Corollary 2.} If $\psi : (L,{\mathcal I})\to (F, {\mathcal G})$ is a surjective morphism from a complete enriched Lie algebra to a profree Lie algebra then there is a morphism
$$\sigma : (F, {\mathcal G})\to (L, {\mathcal I})$$
such that $\psi\circ \sigma = id_F$.

\vspace{3mm}\noindent {\sl proof.} Let $\varphi : \land V_F\to \land V$ be the injective morphism of the corresponding quadratic models which dualizes to $\psi$.   Write $V = \varphi (V_F)\oplus S$. Inverting $\varphi$ then defines an isomorphism $\gamma$ from the sub quadratic algebra $\land \varphi(F)$ to $\land V_F$ and clearly $\gamma\circ \varphi = id$. It remains to extend $\gamma$ to $S$ so that $\gamma (dw)= d\gamma (w)$, $w\in S$.

For this set $\gamma = 0$ in $S\cap V_0$. Now suppose $\gamma$ has been defined in $S\cap V_n$ and let $\{z_i\}$ be a basis of a direct summand of $S\cap V_n$ in $S\cap V_{n+1}$. Then $\gamma (dz_i)$ is a cycle in $\land^2V_F$. Hence $\gamma (dz_i)= dw_i$ for some $w_i\in V_F$. Extend $\gamma $ by setting $\gamma z_i = w_i$. \hfill$\square$

 \vspace{3mm} Now recall from \cite[Chap 9]{RHTII} that the category of a minimal Sullivan algebra, $\land V$, is the least $p$ (or $\infty$) such that $\land V$ is a homotopy retract of $\land V/\land^{>p}V$.
 
 \begin{Prop}
 \label{pXbis}
 (i) A complete enriched Lie algebra is profree if and only if its quadratic model satisfies cat$(\land V)= 1$.  
 
 (ii) Any closed sub Lie algebra of a profree Lie algebra is profree.
 \end{Prop}
 
 \vspace{3mm}\noindent {\sl proof.} (i) The condition cat$(\land V)= 1$ for a quadratic Sullivan algebra is equivalent to the condition $\mathbb Q \oplus (V\cap \mbox{ker}\, d) \stackrel{\simeq}{\to} \land V$. Thus (i) follows from Theorem \ref{t1}.
 
 (ii) Suppose $\land V\to \land Z$ is the morphism of quadratic models corresponding to an inclusion $E\to L$ of a closed sub Lie algebra in a complete enriched Lie algebra. Then \cite[Theorem 9.3]{RHTII} gives cat$(\land Z)\leq $ cat$(\land V)$. Thus (ii) follows from (i). 
  \hfill$\square$

\begin{Prop}
\label{p9}
(i) Suppose $E= R \oplus E^2$ is a free graded Lie algebra. If dim$\, R<\infty$ and $E$ is equipped with the unique (Proposition 4) enriched structure then $\overline{E}$ is profree.

(ii) If $L= T\oplus L^{(2)}$ is a  decomposition of a profree Lie algebra in which $T$ is closed then the sub Lie algebra generated by $T$ is a free Lie algebra.
\end{Prop}

\vspace{3mm}\noindent {\sl proof.} (i) By Proposition 2, $\overline{E}= R\oplus \overline{E}^{(2)}$. Moreover any linear map from $R$ into a complete enriched Lie algebra, $F$, is coherent, and therefore extends to a  morphism $E\to F$. This therefore extends to a morphism $\overline{E}\to F$, and so by definition, $\overline{E}$ is profree.

(ii) It is sufficient to show that any finite dimensional subspace $R$ of $T$ generates a free Lie algebra. It follows from Lemma \ref{l6} that $T = R\oplus Q$ in which $Q$ is also closed. Then division by $Q$ is a surjection $T\to R$ and this is a coherent linear map.

Now let $\overline{E}$ be the completion of the free Lie algebra generated by $R$ with respect to its unique enriched structure. The surjection $T\to R$ extends to a morphism
$$L\to \overline{E}$$
which maps the sub Lie algebra generated by $R$ in $L$ onto the free Lie algebra generated by $R$ in $\overline{E}$. It follows that the sub Lie algebra generated by $R$ in $L$ is free.

 \hfill$\square$

\begin{Prop}\label{p11} Suppose $L= T\oplus L^{(2)}$ is a profree Lie algebra. Then
\begin{enumerate}
\item[(i)]   If $S\subset T$ is any closed subspace then the closure $\overline{E}$ of the sub Lie algebra $E$ generated by $S$ satisfies $$\overline{E} = S\oplus \overline{E}^{(2)},$$ and $\overline{E}^{(2)} = \overline{E}\cap L^{(2)}$.
\item[(ii)] If $0\neq E$ is a solvable sub Lie algebra of the profree Lie algebra, $L$, then $E$ is a free Lie algebra on a single generator.
\end{enumerate}\end{Prop}

\vspace{3mm}\noindent {\sl proof.} (i)   Since the Lie algebra $E_T$ generated by $T$ is a free Lie algebra and $E$ is a retract of $E_T$, $E$ is also a free Lie algebra. Choose a decomposition $T = S\oplus S'$ in which $S'$ is also closed. This gives a coherent retraction $E_T\to E$ which then extends uniquely to a retraction $\rho : L\to \overline{E}$. In particular the identity of $\overline{E}^{(2)}$ factors as
$$\overline{E}^{(2)} \to \overline{E}\cap L^{(2)} \stackrel{\rho}{\to} \overline{E}^{(2)},$$
because $\rho : L^{(2)}\to \overline{E}^{(2)}$. Since $\rho_{\vert E}$ is injective, $\rho: \overline{E}\cap L^{(2)}\to \overline{E}^{(2)}$ is injective, and it follows that
$$\overline{E}^{(2)} = \overline{E}\cap L^{(2)}.$$
This gives $\overline{E}= S\oplus \overline{E}^{(2)}$. 

(ii) Denote by $E \supset \dots \supset E^{[k]}\supset \dots$ the sequence of ideals defined by $E^{[k+1]} = [E^{[k]}, E^{[k]}]$. By hypotheses, some $E^{[n+1]}= 0$. Then, by induction on $n$, we may assume $[E,E]$ is either zero or a free Lie algebra on a single generator. It follows that the closure, $\overline{E}\subset L$ satisfies dim$\, \overline{E}^{(2)}\leq 2$, since (Lemma \ref{l4}) $\overline{E}^{(2)} = \overline{E^2}$.

But since $E\neq 0$, Proposition 13 gives that $\overline{E}= S\oplus \overline{E}^{(2)}$ and $S$ generates a free Lie algebra. This implies that $\overline{E}$ is the free Lie algebra generated by a single element. In particular dim$\, \overline{E}\leq 2$ and so $E= \overline{E}$. 

\hfill$\square$

\vspace{3mm}\noindent {\bf Remark.} Recall that any sub Lie algebra of a free Lie algebra is free (\cite{Sh}). By Proposition 12(i) the analogous statement for complete Lie algebras and closed sub algebras is also true.

 \subsection{The structure of a profree Lie algebra}
 
 \begin{Prop}
 \label{structure}
 Let $L$ be a profree Lie algebra.
 \begin{enumerate}
 \item[(i)] When dim$\, L/[L,L]<\infty$, there is a graded vector space $T$, and 
 $$L \cong \widehat{\mathbb L}(T) = \varprojlim_n \mathbb L_T/\mathbb L^n_T.$$
 \item[(ii)] In the general case, $L= \varprojlim_\sigma L_\sigma$, where the $L_\sigma$ are profree Lie algebras satisfying dim$\, L_\sigma/[L_\sigma, L_\sigma]<\infty$.
 \end{enumerate}
 \end{Prop}
 
 \vspace{3mm}\noindent {\sl proof.} (i). Since $L/[L,L]$ is finite dimensional, $L^{(k)}= L^k$ and $L\cong \widehat{\mathbb L}(sV_0)^\vee$.
 
 (ii)
 The quadratic model, $\land V$, of a profree Lie algebra, $L$, decomposes as the direct limit
 $$\land V= \varinjlim_\sigma \land V(\sigma)$$
 where the sub quadratic algebras $\land V(\sigma)$ are characterized by the two conditions
 $$\mbox{dim}\, V(\sigma)\cap \mbox{ker}\, d<\infty, \hspace{3mm}\mbox{and } V(\sigma ) \cap \mbox{ker}\, d \stackrel{\cong}{\longrightarrow} H^{\geq 1}(\land V(\sigma)).$$
 This in turn yields the isomorphism
 $$L \stackrel{\simeq}{\longrightarrow} \varprojlim_\sigma L(\sigma)$$ of the corresponding homotopy Lie algebras.
 \hfill$\square$
 
 \vspace{3mm} With the notations of the proof of Proposition \ref{structure}, denote by $\rho_\sigma : L\to L(\sigma)$ the dual of the inclusion $V(\sigma)\to V$.
 
 \begin{lem}
 \label{l10} If $S\subset L$ is any subspace then
 $$\overline{S} = \varprojlim_\sigma \overline{\rho_\sigma(S)}.$$
 \end{lem}
 
 \vspace{3mm}\noindent {\sl proof.} Set $K(\sigma) = \{ v\in V(\sigma)\, \vert \, <v, s\rho_\sigma(S)>= 0\}$, and $K= \varinjlim_\sigma K(\sigma)$. It is straightforward to check that $K = \{ v\in V\, \vert\, <v, sS>= 0\}$. It follows from Lemma \ref{l6} that
 $$
 \renewcommand{\arraystretch}{1.5}
 \begin{array}{ll}
 \overline{S} = \mbox{image}\, \left( \, (V/K)^\vee \to L\, \right)& = \varprojlim_\sigma \, \mbox{image}\, \left(\, (V(\sigma)/K(\sigma))^\vee \to L(\sigma)\right) \\ &= \varprojlim_\sigma \overline{\rho_\sigma(S)}.
 \end{array}
 \renewcommand{\arraystretch}{1}$$
 \hfill$\square$
 
 \vspace{3mm}\noindent {\bf Definition.} A vector space $S$ is called \emph{profinite} if $S$ is the projective limit $S= \varprojlim_\alpha S_\alpha$ of a family of finite dimensional vector spaces indexed by a directed set. Associate to $S$ the enriched Lie algebra 
 $$L_S = \varprojlim_{\alpha, n} \mathbb L(S_\alpha)/\mathbb L^{>n}(S_\alpha).$$
 We denote then by $I_\alpha$ the kernel of the projection $S\to S_\alpha$.
 
 \begin{Prop}
 \label{structure2}
 The enriched Lie algebra $L_S$ is a profree Lie algebra. Conversely, if $L= T\oplus L^{(2)}$ is a profree Lie algebra and $T$ is closed, then $L= L_T$.
 \end{Prop}
 
 \vspace{3mm}\noindent {\sl proof.} By construction $L_S = S\oplus L_S^{(2)}$. Now let $E= \varprojlim_\beta E_\beta$ be a complete enriched Lie algebra and $f : S\to E$ be a coherent linear map. Denote by $J_\beta$ the kernel of the projection $p_\beta: E\to E_\beta$. Since $f$ is coherent, for each $\beta$ there is $\alpha$ with $f(I_\alpha) \subset J_\beta$. Then $p_\beta\circ f$ factors through $S_\alpha$, and since $L_\beta$ is nilpotent, we get a map $\mathbb L(S_\alpha)/\mathbb L^{>n}(S_\alpha) \to E_\beta$ for some integer $n$. By composition with the projection $L_S\to \mathbb L(S_\alpha)/\mathbb L^{>n}(S_\alpha)$ this gives a map $f_\beta: L_S\to E_\beta$. Finally, together the $f_\beta$ define a map $f : L_S\to E$.
 
 The converse is Proposition \ref{structure}(ii).
 
 \hfill$\square$
 
 \vspace{3mm} Again, suppose that $\land V$ is the quadratic model of a profree Lie algebra, $L$. Recall that for some closed subspace $T\subset L$ we have $L= T\oplus L^{(2)}$, so that $T$ generates a free Lie algebra, $E$, with $E^2\subset L^{(2)}$. In particular
 $$E = \oplus_{k\geq 1} T(k),$$
 where $T(k)$ is the linear span of the iterated Lie brackets of length $k$ in elements of $T$. In general,   however the subspaces $T(k)$ may not be closed. Nevertheless we do have an analogous structure for $L$.
 
 \begin{Prop}
 \label{p12}
 If $L$ is a profree Lie algebra then, with the notation above and of Proposition 15(ii),
 \begin{enumerate}
 \item[(i)] $L^{(k)} = \varprojlim_\sigma L(\sigma)^k =\overline{T(k) }\oplus L^{(k+1)},$ $k\geq 1$.
 \item[(ii)] For each $k, \ell$, $[\overline{T(k)}, \overline{T(\ell)}] \subset \overline{T(k+\ell)}$. In particular
 $$F:= \oplus_k \overline{T(k)}$$
 is a weighted sub Lie algebra of $L$.
 \item[(iii)] $L\stackrel{\cong}{\longrightarrow}\varprojlim_n L/J(n)$, where $J(n)= \prod_{k\geq n} \overline{T(k)}$.
 \item[(iv)] $L^{(k)}/L^{(k+1)} \cong \overline{T(k)} \cong \varprojlim_\sigma L(\sigma)^k/L(\sigma)^{k+1}$. 
 \end{enumerate}
 \end{Prop}
 
 \vspace{3mm}\noindent {\sl proof.} (i) The inclusions $V(\sigma)\subset V$ dualize to surjections $\rho_\sigma : L\to L(\sigma)$ and, as observed above, since $V = \varinjlim_\sigma V(\sigma)$, these induce an isomorphism
 $$L\stackrel{\cong}{\longrightarrow} \varprojlim_\sigma L(\sigma).$$
 Moreover, for each $n$, $V_n = \varinjlim_\sigma V(\sigma)_n$, and therefore 
 for each $n$,
 $$V/V_n = \varinjlim_\sigma V(\sigma)/V(\sigma)_n.$$
It follows from Lemma \ref{l6} that $\rho_\sigma : L^{(n)}\to L(\sigma)^{(n)}$ is also surjective and that these define an isomorphism
$$L^{(n)} \stackrel{\cong}{\longrightarrow} \varprojlim_\sigma L(\sigma)^{(n)}.$$

Next observe that since $T$ is closed the inclusion of $T$ in $L$ is (Lemma \ref{l6}) dual to a surjection of $V$ onto a direct summand $V'$ of $V\cap \mbox{ker}\, d$ in $V$. This then restricts to surjections of $V(\sigma)$ onto direct summand of $V(\sigma)\cap \mbox{ker}\, d$ in $V(\sigma)$. These in turn dualize to inclusions $T(\sigma)\to L(\sigma)$ into direct summands of $L(\sigma)^{(2)}$ in $L(\sigma)$. By construction each $\rho(\sigma) : T\to T(\sigma)$ is surjective and these define an isomorphism:
$$T\stackrel{\cong}{\longrightarrow} \varprojlim_\sigma T(\sigma).$$

On the other hand, let $T(\sigma, k)$ be the linear span of the iterated Lie brackets of length $k$ in elements of $T(\sigma)$. Since $L(\sigma)$ is the closure of the free Lie algebra $E(\sigma)$ generated by $T(\sigma)$, the fact that dim$\, T(\sigma)<\infty$ implies (Proposition \ref{p2}) that for each $n$,
$$L(\sigma)^{(n)} = \varprojlim_k E(\sigma)^n/E(\sigma)^{k+n} = T(\sigma, n) \oplus L(\sigma)^{(n+1)}.$$
Since $\rho_\sigma : T\to T_\sigma$ is surjective each $\rho_\sigma : T(k)\to T(\sigma, k)$ is also surjective. Thus Lemma \ref{l10} gives
$$\overline{T(n)} \stackrel{\cong}{\longrightarrow} \varprojlim_\sigma T(n,\sigma).$$

Finally, it follows from Lemma \ref{l6} that
$$L^{(n)}/L^{(n+1)} \to \varprojlim_\sigma L(\sigma)^{(n)}/L(\sigma)^{(n+1)}$$
is the dual of the isomorphism $V_{n+1}/V_n \cong \varinjlim_\sigma V(\sigma)_{n+1}/V(\sigma)_n$. Thus the commutative diagram
$$
\xymatrix{
\overline{T(n)} \ar[d]\ar[rr]^-\cong && \varprojlim_\sigma T(\sigma, n)\ar[d]^\cong\\
L^{(n)}/L^{(n+1)} \ar[rr]^-\cong && \varprojlim_\sigma L(\sigma)^{n)}/L(\sigma)^{(n+1)}
}$$
implies that $\overline{T(n)} \oplus L^{(n+1)} = L^{(n)}$.

(ii) This is immediate from the fact that $[T(k),T(\ell)]\subset T(k+\ell)$ and the fact that the Lie bracket preserves closures.

(iii) This is immediate from the relation
$$J(n) = \varprojlim_{\sigma, k\geq n} T(\sigma, k) = \varprojlim_{\sigma} L(\sigma)^{(n)} = L^{(n)}$$
where the first equality again follows from Lemma \ref{l6}. 

(iv) The first isomorphism is a consequence of (i). The second isomorphism follows from Lemma \ref{l10}.

\hfill$\square$

  \vspace{3mm} The next Proposition explains the relation between $L^2$ and $L^{(2)}$ for general profree Lie algebras.
  
 \begin{Prop}\label{ll2} If $L=L_0$ is a profree Lie algebra   and dim$\, L/[L,L]=\infty$, then $L^2\neq L^{(2)}$.
 \end{Prop}

  \vspace{3mm}\noindent {\sl proof.}
  Indeed, let   $(\land W,d)$ be the quadratic model of $L$. We can decompose $W$ as an union $W = \cup_n W_n$ with $W_0 = W \cap \mbox{ker}\, d$,  and for $n>1$, $W_n= d^{-1}(\land^2W_{n-1})$. We denote by $Z_n$ a direct summand of $W_{n-1}$ in  $W_n$. Then by hypothesis $W_0$ is infinite dimensional and there is a quasi-isomorphism $\varphi : (\land W,d)\to (\mathbb Q \oplus W_0,0)$ that is the identity on $W_0$ and that maps each $Z_n$ to $0$.

For sake of simplicity, we write $V = W_0$ and $Z = Z_1$. By construction the isomorphism $L= (sW)^\vee$    induces isomorphisms
$$L/L^{(2)} \cong (sV)^\vee ,\hspace{1cm}\mbox{and } L^{(2)}/L^{(3)} \cong (sZ)^\vee.$$

  Suppose first $W$ countable and denote by $w_1, w_2, \dots$ a basis of $  V$. Since $d=0$ on $V$, $d: Z\to \land^2V$ is an isomorphism. We denote by $w_{ij}$, $i<j$ the basis of $Z$ defined by $d(w_{ij})= w_i \land w_j$.  Denote by $E$ teh vector space of column matrices $X = (x_i)$ with only a finite number of nonzero $a_i$. Then the map $(a_i) \leadsto \sum a_i w_i$ defines an isomorphism $E\stackrel{\cong}{\to} V$. 
  
 Let represent an element $\varphi\in Z^\vee$ by the infinite dimensional antisymmetric matrix $M_\varphi$,$$(M_\varphi)_{ij}= \varphi (w_{ij}).$$
  In a similar process, an element $f\in V^\vee$ can be represented by a column matrix $A_f$, with $(A_f)_i = f(w_i)$.  The vector space ker$\, f$ can then be identified with the sub vector space of $E$ formed by the column matrices $X$ satisfying $A_f^t\cdot X = 0$. (Here  $A^t$ denotes the line matrix  transposition of a column matrix $A$.)
  
  Note that when we have two column matrices $A$ and $B$, we can form the antisymmetric matrix $A\cdot B^t - B\cdot A^t$. Now remark that for $f, g\in L^2$, we have
  $$M_{[f,g]} = A_f\cdot B_g^t-B_g\cdot A_f^t.$$
  
  Let $\varphi_0\in Z^\vee$ be the particular element defined for $i<j$ by
  $$\varphi_0 (w_{ij}) = \left\{
  \begin{array}{ll} 1 & \mbox{if $w_{ij}= w_{2k+1, 2k+2}$, for some $k$}\\0 & \mbox{otherwise}\end{array}
  \right.
  $$
  The element $\varphi_0$ can be extended to all of $W$, by $\varphi_0(V)=0$ and $ \varphi_0(Z_n)= 0$, for $n>1$. By construction $\varphi_0\in L^{(2)}$. 
  The associated matrix is
  $$M_0 = \left( \begin{array}{cccc}
  B_0 & 0 & 0 & \dots\\
  0 & B_0 & 0 & \dots\\
  0 & 0 & B_0 & \dots \\
  \dots & \dots &\dots &\dots\end{array}
  \right) \hspace{15mm}\mbox{with } B_0 = \left(\begin{array}{cc} 0 & 1\\ -1 & 0\end{array}\right).$$
   
   Now consider a finite sum $\sum_{i=1}^n [f_i, g_i]$, with $f_i$ and $g_i\in V^\vee$. The associated matrix is $\sum_{i=1}^n M_{[f_i, g_i]}$. Then $$K = (\cap_{i=1}^n \mbox{ker}\, f_i) \cap (\cap_{i=1}^n \mbox{ker}\, g_i)$$
   is infinite dimensional, and so for some non-zero $X\in Z^\vee$,   $(\sum M_{[f_i,g_i]})\cdot X = 0$. Since $M_0\cdot X\neq 0$, it follows that $L^2\subset_{\neq} L^{(2)}$.   
   
   In the general case, let $E\subset V$ be a countable subvector space and let $\land T$ be the minimal model of $(\mathbb Q \oplus E,0)$. Then its homotopy Lie algebra $L_T$ is a retract of $L$, 
   $$\xymatrix{ 
   L_T \ar@/^/[r]^j & L\ar@/^/[l]^\rho .}$$
   Let $\varphi_0\in L_T^{(2)}$, not in $L_T^2$.Then $j(\varphi_0)\in L^{(2)}$ and not in $L^2$ because otherwise $\varphi_0= \rho(\varphi)$ would belong to $L_T^{2}$.

   \hfill$\square$

 \newpage
 \part{Topological spaces, Sullivan models, and their homotopy Lie algebras}
 
 \section{The homotopy groups of a Sullivan completion}
 
 Enriched Lie algebras effectively describe the homotopy groups, $\pi_*(X_{\mathbb Q})$ and $\pi_*<\land V>$ for any connected space, $X$, and minimal Sullivan algebra, $\land V$. (Recall that if $\land V$ is the minimal Sullivan model of $X$, then $X_{\mathbb Q}= <\land V>$. Recall also that for simplicity we write $H(X)$ to mean $H^*(X;\mathbb Q)$)
 
 To make this description explicit, \emph{we fix a minimal Sullivan algebra, $(\land V,d)$.} Recall from \cite[Chap 2]{RHTII} that the \emph{homotopy Lie algebra} of $(\land V,d)$ is the graded Lie algebra, $L = L_{\geq 0}$, defined by:
 $$sL = V^\vee \hspace{5mm}\mbox{and } <v, s[x,y]> = (-1)^{deg\, y+1} <d_1v, sx,sy>,$$
 where $d_1v$ is the component of $dv$ in $\land^2V$. Now observe that $(\land V,d)$ determines a natural enriched structure in $L$: it is immediate from the defining condition for $d$ that $V$ is the union of the finite dimensional subspaces $V_\alpha$ for which $\land V_\alpha$ is preserved by $d$. The inclusions $V_\alpha \to V$ then dualize to surjections
 $$L\to L/I_\alpha = L_\alpha$$
 onto finite dimensional nilpotent Lie algebras, and this yields an enriched structure $(L, \{I_\alpha\})$ for $L$.
 
 On the other hand, the linear map $d_1: V\to \land^2V$ extends to a derivation in $\land V$, and $(\land V, d_1)$ is a quadratic Sullivan algebra: the \emph{quadratic Sullivan algebra associated with $(\land V,d)$}. It is immediate that $(L, \{I_\alpha\})$ is  the homotopy Lie algebra for $(\land V, d_1)$, introduced in \S 4. In particular, it is complete.
 
 \vspace{3mm}\noindent {\bf Remark.} The choice of generating space $V$ can be modified by a map, $v\mapsto v+ \sigma (v)$ in which $\sigma : V\to \land^{\geq 2}V$, without changing the associated quadratic Sullivan algebra $(\land V, d_1)$ or the   enriched Lie algebra $(L, \{I_\alpha\})$, although it may replace the ideals $I_\alpha$ by an equivalent set of ideals.
 
 \vspace{3mm}Now consider the adjoint bijections
 \begin{eqnarray}\label{SN1}
 \mbox{Cdga} (\land V, A_{PL}(X)) = \mbox{Simpl} (X, <\land V>)
 \end{eqnarray}
 for any connected space $X$.  In particular, adjoint to $id_{<\land V>}$ is a morphism
 $$\varphi_{\land V} : \land V\to A_{PL}<\land V>.$$
 Then adjoint to any morphism $\varphi : \land V\to A_{PL}(X)$ from a minimal Sullivan algebra $\land V$ is a simplicial map
 $$<\varphi> : X\to <\land V>.$$
 A straightforward check from the definitions then shows that $\varphi$ decomposes as the composite
 \begin{eqnarray}
 \label{S11}
 \varphi : \xymatrix{\land V \ar[rr]^{\varphi_{\land V}} && A_{PL}<\land V>\ar[rr]^{A_{PL}<\varphi>} && A_{PL}(X).
 }
 \end{eqnarray}
 There follows
 
\begin{lem}  \label{lS12} Suppose $\varphi : \land V\to A_{PL}(X)$ is a morphism from a minimal Sullivan algebra. If $A_{PL}<\varphi>$ is a quasi-isomorphism then $\varphi_{\land V}$ is a quasi-isomorphism if and only if $\varphi$ is the minimal Sullivan model of $X$.
 \end{lem}

\vspace{3mm}
Next observe that (\ref{SN1}) provides bijections, \emph{linear for $n\geq 2$},
 $$\pi_n<\land V> = [\land V, A_{PL}(S^n)] = (V^n)^\vee, \hspace{1cm} n\geq 1.$$
Since $\land V$ is minimal, the retraction $\xi : \land^{\geq 1}V\to V$ with kernel $\land^{\geq 2}V$ induces a linear map $H(\xi) : H^{\geq 1}(\land V)\to V$. It is immediate from the definitions that for $x\in V^\vee= \pi_*<\land V>$ and $[\Phi]\in H(\land V)$, we have
\begin{eqnarray}
\label{red}
<H(\xi)[\Phi],x>= <H(\varphi_{\land V})[\Phi], \mbox{hur}\, x>
\end{eqnarray}
where hur: $\pi_*<\land V>\to H_*<\land V>$ denotes the Hurewicz homomorphism.

Since $(V^n)^\vee = L_{n-1}$, the bijection $\pi_n<\land V> = (V^n)^\vee$ also identifies   $\pi_*<\land V>$ as the suspension $sL$ of $L$.
With this identification, the enriched Lie algebra structure of $L$ provides explicit formulas for the product in $\pi_1<\land V>$, the action of this group on $\pi_n<\land V>, n\geq 2$, and for the   Whitehead products.

 The formulas involving $\pi_1$ depend on a functor $L\leadsto G_L\subset \overline{UL_0}$ from complete enriched Lie algebras to groups, together with a natural bijection
 $$\exp : L_0\stackrel{\cong}{\longrightarrow} G_L.$$
 These are defined as follows. Since each $L_\alpha$ is finite dimensional, it is a classical fact (\cite{Serre}, \cite{Q}, \cite[Chap 2]{RHTII}) that the standard power series is a bijection,
 $$\exp : (L_\alpha)_0 \stackrel{\cong}{\longrightarrow} G_{L_\alpha}$$
 onto the group of group-like elements in $\widehat{UL_\alpha}$.

Passing to inverse limits yields   a functor $L\leadsto G_L\subset \overline{UL_0}$ from complete enriched Lie algebras to groups, and the natural bijection
 $$\exp : L_0\stackrel{\cong}{\longrightarrow} G_L.$$
 Moreover, the Baker-Campbell-Hausdorff series in $L_0$ is convergent, and provides an explicit formula for the product.
 
 Next recall (Corollary 2 to Proposition 2) that $L_0$ is the direct limit of pronilpotent sub Lie algebras $L(\sigma)$ satisfying dim$\, L(\sigma)/L(\sigma)^2<\infty$. It follows that $G_L$ is the direct limit of the corresponding groups $G(\sigma)$, and that $\exp$ is the direct limit of the bijection $\exp (\sigma)$.
 
 Now denote by $G^n$ the subgroup of a group $G$ generated by iterated commutators of length $n$. The direct limits above then induce bijections
 $$\varinjlim_\sigma L(\sigma)^n \stackrel{\cong}{\longrightarrow} L_0^n \hspace{5mm}\mbox{and } \varinjlim_\sigma G(\sigma)^n \stackrel{\cong}{\longrightarrow} G^n(\sigma).$$
 Moreover \cite[Sec. 2.4]{RHTII} shows that $\exp$ restricts to a family of bijections $L(\sigma)^n \stackrel{\cong}{\longrightarrow} G(\sigma)^n$. Therefore $\exp$ restricts to bijections
 $$\exp : L_0^n \stackrel{\cong}{\longrightarrow} G_L^n\,, \hspace{5mm} n\geq 1.$$
 
 Moreover, again by \cite[Sec. 2.4]{RHTII} for each $\sigma$ the bijections $L(\sigma)^n \stackrel{\cong}{\to} G(\sigma)^n$ factor to give bijections,
 $$L(\sigma)^n/L(\sigma)^{n+k} \stackrel{\cong}{\longrightarrow} G(\sigma)^n/G(\sigma)^{n+k}, \hspace{5mm} k\geq 1,$$
 which are linear isomorphisms when $k= 1$.
 
 Since the direct limit of quotients is the quotient of the direct limits these define bijections, which are linear when $k= 1$,
 $$L_0^n/L_0^{n+k} \stackrel{\cong}{\longrightarrow} G_L^n/G_L^{n+k}, \hspace{1cm} k\geq 1.$$
 In particular each $G_L^n/G_L^{n+1}$ is a rational vector space. Finally, and in the same way, the diagrams
 $$\xymatrix{
 L_0^q/L_0^{q+1} \otimes L_0^p/L_0^{p+1} \ar[rr]\ar[d]^\cong && L_0^{p+q}/L_0^{p+q+1}\ar[d]^\cong\\
 G_L^q/G_L^{q+1}\otimes G_L^p/G_L^{p+1} \ar[rr] && G_L^{p+q}/G_L^{p+q+1},}$$
 induced by the respective adjoint actions, commute.
 
 \vspace{3mm}

Moreover, all these relations translate immediately to $\pi_1<\land V>$ as follows: Theorem 2.4 in \cite{RHTII} provides natural isomorphisms $G_{L_\alpha}\cong \pi_1(\land V_\alpha)$ of groups. Passing to inverse limits yields the natural isomorphism
 $$G_L\cong \pi_1<\land V>.$$
 Then all the other relations translate simply by replacing $G_L$ by $\pi_1<\land V>$.
 
 Thus Theorem 2.5 in \cite{RHTII} shows that the right action of $\pi_1<\land V>$ on $\pi_n<\land V>$, $n\geq 2$, is given by 
 $$\beta\cdot \alpha = [\beta,\alpha], \hspace{1cm} \alpha \in L_0, \beta\in L_n, n\geq 1,$$
while, for $k, \ell\geq 2$, the Whitehead products
 $$\pi_k<\land V> \times \pi_\ell <\land V>\to \pi_{k+\ell -1}<\land V>$$
also translate the Lie bracket in $L$.
 
 \vspace{3mm}
Finally, suppose $\varphi : \land V \to \land W$ is a morphism of minimal Sullivan algebras. Filtering by wedge degree yields a morphism $\varphi_1 : (\land W, d_1)\to (\land V, d_1)$ of the associated quadratic algebras, whose restriction to $W$ dualizes to a morphism
 $$L_\varphi : L_V \leftarrow L_W$$
 between the homotopy Lie algebras. This is almost by definition a morphism of enriched Lie algebras. In particular it preserves the group structure in $\pi_1$, its action on $\pi_n$, $n\geq 1$, and the higher Whitehead products. Moreover, a standard argument establishes the
 
 \begin{lem}\label{l11}
 Suppose $\varphi, \psi : \land V\to \land W$ are two morphisms. Then the following conditions are equivalent:
 \begin{enumerate}
 \item[(i)] $\varphi$ and $\psi$ are based homotopic.
 \item[(ii)] $\varphi_1=\psi_1$.
 \item[(iii)] $L_\varphi = L_\psi$.
 \end{enumerate}
 \end{lem}
 
 \section{Sullivan rationalizations}
 
Let $\varphi_X : \land V\stackrel{\simeq}{\to} A_{PL}(X)$ be a minimal Sullivan model for a connected space, $X$, and let $L_X$ denote the homotopy Lie algebra of $\land V$. The adjoint map   
  $$<\varphi_X> : X\to  <\land V>:= X_{\mathbb Q}$$
is called   the \emph{Sullivan completion} of $X$. The discussion in \S 8 then identifies 
  $$\pi_*(X_{\mathbb Q}) = sL_X$$
  and shows how the Lie bracket and enriched structure in $L_X$ determine the product in $\pi_1(X_\mathbb Q)$, its action in $\pi_{\geq 2}(X_\mathbb Q)$, and the higher order Whitehead products.
  
  Thus while $H(X)$ and $\pi_*(X_{\mathbb Q})$ are directly computable from $\land V$, $H(X_{\mathbb Q})$ may not be suh a simple invariant of the model. For example, $H^2((S^1\vee S^1)_{\mathbb Q})$ is uncountably infinite (\cite{IM}). This leads to the 
  
  \vspace{3mm}\noindent {\bf Definition.} If a Sullivan completion $<\varphi_X> : X\to X_{\mathbb Q}$ satisfies
  $$H(<\varphi_X>) : H(X_{\mathbb Q}) \stackrel{\cong}{\longrightarrow} H(X)$$
  then $<\varphi_X>$ is a \emph{Sullivan rationalization} of $X$.
  
  If $<\varphi_X>$ is a homotopy equivalence, then $X$ is \emph{Sullivan rational}.
  
  \vspace{3mm}\noindent {\bf Remarks. 1.} The condition that $<\varphi_X>$ be a Sullivan rationalization implies that the minimal Sullivan model of $X$ directly computes both $H(X_{\mathbb Q})$ and $\pi_*(X_{\mathbb Q})$.
  
  {\bf 2.} If $X_{\mathbb Q}$ is Sullivan rational then clearly $A_{PL}<\varphi_X>$ is a quasi-isomorphism and $<\varphi_X>$ is a Sullivan rationalization. Conversely, if $<\varphi_X>$ is a Sullivan rationalization then by Lemma 11 it identifies the minimal Sullivan model $\land V$ of $X$ with a minimal model of $<\land V>=X_{\mathbb Q}$. Thus $X_{\mathbb Q}\simeq (X_{\mathbb Q})_{\mathbb Q}$ and $X_{\mathbb Q}$ is Sullivan rational.

 \vspace{3mm} Recall that(\cite{Q})   Quillen's rationalization   for simply connected spaces assigns to each such space a map $Y\to \mathbb Q (Y)$ which induces isomorphisms
  $$H_{\geq 1}(Y)\otimes \mathbb Q \stackrel{\cong}{\to} H_{\geq 1}(\mathbb Q(Y)) \hspace{5mm}\mbox{and } \pi_*(Y)\otimes \mathbb Q \stackrel{\cong}{\to} \pi_*(\mathbb Q (Y)).$$
  In particular, $\mathbb Q (\mathbb Q Y) = \mathbb Q (Y)$.

 Sullivan's completion, $X_{\mathbb Q}$, is analogous to Quillen's rationalization, and also to the Bousfield-Kan completion $\mathbb Q_{\infty}(X)$ - cf \cite{BK}. Moreover, if $H(X)$ has finite type then $\mathbb Q_\infty (X)$ and $X_{\mathbb Q}$ are homotopy equivalent (\cite{BG}). But it is not always true that every $X_{\mathbb Q}$ is Sullivan rational. For instance (\cite{IM}), $H^2((S^1\vee S^1)_{\mathbb Q})$ has uncountable dimension, where as if $X_{\mathbb Q}$ is Sullivan rational by the Corollary to Theorem 2 below, $H(X_{\mathbb Q})$ has finite type.

 \vspace{3mm}\noindent {\bf Remarks} {\bf 1.} If $\varphi_X : \land V\to A_{PL}(X)$ is a minimal Sullivan model then $X$ is Sullivan rational if and only if each $$\pi_k<\varphi_X> : (V^k)^\vee \stackrel{\cong}{\longrightarrow} \pi_k(X).$$

{\bf 2.}  Let $Y$ be a Sullivan rational space and let $X$ be a connected space. Since $f_{\mathbb Q}\circ \varphi_X= \varphi_Y\circ f$, and $\varphi_Y$ is a homotopy equivalence, it follows that $f\mapsto f_{\mathbb Q}$ yields an injection
$  [X,Y]\to   [X_{\mathbb Q}, Y_{\mathbb Q}].$
 
 \vspace{3mm} The next Theorem extends results in \cite[Chap.7]{RHTII}.
 
 \begin{Theorem}
 \label{ntt2}
 Suppose $X$ is a connected space with fundamental group $G_X$ and universal covering space $\widetilde{X}$. Then the following conditions are equivalent:
 \begin{enumerate}
 \item[(i)] $X$ is Sullivan rational.
 \item[(ii)] $BG_X$ and $\widetilde{X}$ are Sullivan rational, and $G_X$ acts locally nilpotently in $H(\widetilde{X})$.
 \end{enumerate}
 \end{Theorem}

\vspace{3mm}\noindent {\bf Corollary.}  {\sl   If a connected space $X$ is Sullivan rational, then
  \begin{enumerate}
  \item[(i)] Each $\pi_1(X)^n/\pi_1(X)^{n+1}$ is a finite dimensional rational vector space.
  \item[(ii)] $\pi_1(X)= \varprojlim_n \pi_1(X)/\pi_1(X)^n$.
  \item[(iii)] Each $\pi_k(X)$, $k\geq 2,$ is a rational vector space and a finite dimensional nilpotent $\pi_1(X)$-module.
  \item[(iv)] $H(X)$ is a graded vector space of finite type.
  \end{enumerate}}

\vspace{3mm}Before proceeding to the proof of Theorem \ref{ntt2} and of its Corollary, we establish two Lemmas.

 \begin{lem} 
 \label{lemmaA}Suppose for some minimal Sullivan algebra, $\land V$, that $\varphi_{\land V} : \land V\to A_{PL}<\land V>$ is a quasi-isomorphism. Then
 \begin{enumerate}
 \item[(i)] $<\land V>$ is Sullivan rational.
 \item[(ii)] A basis of $V$ is at most countable.
 \item[(iii)] $H(\land V)$ is a graded vector space of finite type.
 \end{enumerate}
 \end{lem}
 
 \vspace{3mm}\noindent {\sl proof.} (i) Since $<\varphi_{\land V}> = id_{<\land V>}$, $<\varphi_{\land V}>$ is a homotopy equivalence and $<\land V>$ is Sullivan rational by definition.
 
 (ii) Formula (\ref{red}) identifies the surjection $V^\vee \to (V\cap \mbox{ker}\, d)^\vee$ as a composite
 $$\xymatrix{V^\vee \ar[rr]^{\mbox{hur}} &&H_*<\land V> \ar[rr] && (V\cap \mbox{ker}\, d)^\vee.}$$
 Thus $H_*<\land V>\to (V\cap \mbox{ker}\, d)^\vee$ is surjective. Dualizing gives injections $(V^n\cap \mbox{ker}\, d)^{\vee\vee} \to H^n(\land V).$
 
 Now assume by induction that $V^{<n}$ has an at most countable basis and let $Z$ be the space of cycles in $(\land V)^n$. Then $Z\subset \land V^{<n}\oplus V^n$, and division by $V^n\cap \mbox{ker}\, d$ gives an exact sequence
 $$0\to V^n\cap \mbox{ker}\, d\to Z\to \land V^{<n}.$$
 Since $(V^n\cap \mbox{ker}\, d)^{\vee\vee}$ embeds in $H^n(\land V)$, it follows that dim$\, V^n\cap \mbox{ker}\, d<\infty$. 
 
 On the other hand, $V^n= \cup_{p\geq 0}V^n_p$ with $V^n_0= V^n \cap d^{-1}(\land V^{<n})$ and $V_{p+1}^n = V^n \cap d^{-1}(\land V^{<n}\otimes \land V^n_p)$. Since dim$\,V^n\cap \mbox{ker}\, d<\infty$, it follows by induction on $p$ that each $V^n_p$, and hence $V^n $ itself, has an at most countable basis.
 
 (iii) It follows from (i) that each $H^n(\land V)$ has an at most countable basis. But since $H^n(\land V)= H^n(<\land V>)= H_n<\land V>^\vee$, this implies that each dim$\, H^n(\land V)<\infty$. \hfill$\square$

\begin{lem}
\label{lemmaB} Suppose $\land V$ is a minimal Sullivan algebra. If $V^1=0$ then the following  conditions are equivalent:
\begin{enumerate}
\item[(i)] $\varphi_{\land V} : \land V\to A_{PL}<\land V>$ is a quasi-isomorphism.
\item[(ii)] Each dim$\, H^k(\land V)<\infty.$
\item[(iii)] Each dim$\, V^k<\infty$.
\item[(iv)] Each dim$\, \pi_k<\land V><\infty.$
\end{enumerate}
\end{lem}

\vspace{3mm}\noindent {\sl proof.} (i) $\Rightarrow$ (ii). This is established in Lemma \ref{lemmaA}.

(ii) $\Rightarrow $ (iii). This is immediate because, since $V^1= 0$, $d : V^k\to \land^{\geq 2}V^{<k}$.

(iii) $\Rightarrow$ (iv) This follows from the isomorphism $(V^k)^\vee \cong \pi_k<\land V>$ established in \S 8. 

(iv) $\Rightarrow$ (i). The morphisms $\land V^n\to A_{PL}<\land V^n>$ map $V^n$ isomorphically to the dual of $(V^n)^\vee = \pi_n<\land V^n>$, as shown in (10). Since $<\land V^n>$ is the Eilenberg-MacLane space $K(\pi_n<\land V^n>,n)$ it follows that $\land V^n \stackrel{\simeq}{\longrightarrow} A_{PL}<\land V^n>.$ 
 
 Now suppose by induction that $\land V^{<n}\to A_{PL}<\land V^{<n}>$ is a quasi-isomorphism. Then in the commutative diagram,
 $$\xymatrix{
 \land V^{<n} \ar[d] \ar[rr] && \land V^{<n}\otimes \land V^n \ar[d] \ar[rr] && \land V^n\ar[d]\\
 A_{PL}<\land V^{<n}> \ar[rr] && A_{PL}<\land V^{\leq n}> \ar[rr] && A_{PL}<\land V^n>,}$$
 the left and right morphisms are quasi-isomorphisms. But, \cite[Proposition 17.9]{FHTI}
 $$\xymatrix{
 <\land V^{\leq n-1}> && <\land V^{\leq n}>\ar[ll] && <\land V^n>\ar[ll]}$$
 is a fibration. A standard Serre spectral sequence shows that $\land V^{\leq n}\to A_{PL}<\land V^{\leq n}>$ is a quasi-isomorphism. 
 \hfill$\square$

 \vspace{3mm}\noindent {\sl proof of Theorem 2.}  Denote by $G_X$, $BG_X$, and $\widetilde{X}$ the fundamental group of $X$, its classifying space, and the homotopy fibre of $X\to BG_X$. Then let $\varphi_X : \land V\stackrel{\simeq}{\to} A_{PL}(X)$ be the minimal Sullivan model of $X$. Now \cite[formula (7.4)]{RHTII} gives the commutative diagram
 \begin{eqnarray}
 \label{xa}
 \xymatrix{
 \land V^1\ar[rr]\ar[d]_{\varphi_B} && \land V\ar[rr]\ar[d]^\simeq_{\varphi_X} && \land V^{\geq 2}\ar[d]_{\widetilde{\varphi}}\\
 A_{PL}(BG_X) \ar[rr] && A_{PL}(X) \ar[rr] &&A_{PL}(\widetilde{X}).}
 \end{eqnarray}
 The corresponding adjoint diagram is then
 \begin{eqnarray}
 \label{xb}
\xymatrix{ <\land V^1>  &&<\land V>\ar[ll] && <\land V^{\geq 2}>\ar[ll]\\
 BG_X \ar[u]_{<\varphi_B>} && X\ar[ll] \ar[u]_{<\varphi_X>} && \widetilde{X}.\ar[ll]\ar[u]_{<\widetilde{\varphi}>}
 }
 \end{eqnarray}
 As observed in \cite[Prop. 17.9]{FHTI}, (\ref{xb}) is a map of fibrations. In particular, $\pi_*(<\varphi_B>)= \pi_1<\varphi_X>$ and $\pi_*<\widetilde{\varphi}>= \pi_{\geq 2}<\varphi_X>$. Thus  $<\varphi_X>$ is a homotopy equivalence if and only if  both $<\varphi_B>$ and $<\widetilde{\varphi}>$ are homotopy equivalences. In particular, if either (i) or (ii) holds, then each of $A_{PL}(<\varphi_B>)$, $A_{PL}(<\varphi_X>)$ and $A_{PL}(\widetilde{\varphi})$ are quasi-isomorphisms. 
 Moreover, in this case
 $$\pi_1<\varphi_B> : \pi_1<\land V^1>\stackrel{\cong}{\longrightarrow} G_X.$$

 \vspace{3mm} Now we show separately that (i) $\Rightarrow $ (ii) and (ii) $\Rightarrow$ (i).
 
 \vspace{2mm}\noindent (i) $\Rightarrow$ (ii). Since $\varphi_X$ is a quasi-isomorphism and $<\varphi_X>$ is a homotopy equivalence, Lemma 12 implies that $\varphi_{\land V} : \land V\to A_{PL}<\land V>$ is a quasi-isomorphism. Therefore (Lemma \ref{lemmaA}) $V$ has a countable basis and (Lemma \ref{lemmaB})
 $$\varphi_{\land V^{\geq 2} } : \land V^{\geq 2} \to A_{PL}<\land V^{\geq 2}>$$
 is a quasi-isomorphism. Since, as observed above, $A_{PL}<\varphi_{\land V^{\geq 2}}>$ is also a quasi-isomorphism it follows (Lemma 12) that $\widetilde{\varphi}$ is a minimal Sullivan model of $\widetilde{X}$. As noted above, $<\widetilde{\varphi}>$ is a homotopy equivalence. Therefore by definition, $\widetilde{X}$ is Sullivan rational.
 
 On the other hand, it follows from Lemma \ref{lemmaA} that dim$\, H^1(\land V)<\infty$, and from Lemma \ref{lemmaB} that dim$\, V^k<\infty$, $k\geq 2$. Thus the hypotheses of \cite[Theorem 7.8(i)]{RHTII} are satisfied. But the proof of that Theorem shows that
 $$\varphi_B : \land V^1\to A_{PL}<\land V^1>$$
 is a quasi-isomorphism. Since $<\varphi_B>$ is also a homotopy equivalence it follows that $BG_X$ is Sullivan rational.

 Finally \cite[Corollary 4.3]{RHTII} implies that the holonomy representation of $BG_X$ in $H(\widetilde{X})$ is nilpotent.
 
 \vspace{3mm}\noindent (ii) $\Rightarrow $ (i). Again consider diagram (12). As observed in the proof of \cite[Theorem 7.1]{RHTII}, $\varphi_B$ extends to a minimal Sullivan model of the form $\land V^1\otimes \land Z^{\geq 2}\stackrel{\simeq}{\longrightarrow} A_{PL}(BG_X)$. Since $BG_X$ is Sullivan rational it follows that $(V^1\otimes Z^{\geq 2})^\vee\cong \pi_1(BG_X)$, and so $Z^{\geq 2}= 0$. Therefore $\varphi_B : \land V^1\stackrel{\cong}{\longrightarrow } A_{PL}(BG_X)$. This, together with the hypothesis on the action of $G_X$ on $H(\widetilde{X})$ allows us to apply \cite[Theorem 5.1]{RHTII} and conclude that $\widetilde{\varphi}$ is a quasi-isomorphism. But by hypothesis, $\widetilde{X}$ is Sullivan rational. Therefore $<\widetilde{\varphi}>$ is a homotopy equivalence. This in turn implies that $<\varphi_X>$ is a homotopy equivalence and so $X$ is Sullivan rational. 
 
 \hfill$\square$
 
 \vspace{3mm}\noindent {\sl proof of the Corollary to Theorem 2.} 
 
Let $L_0$ be the homotopy Lie algebra of $\land V^1$. Then \S 8 gives isomorphisms of groups,
 $$L_0/[L_0,L_0] \cong G_X/[G_X, G_X] = H_1(BG_X;\mathbb Z).$$
 In particular, $H_1(G_X;\mathbb Z)$ is a rational vector space. 
 
 On the other hand, Lemma \ref{lemmaA} implies that $H^1(\land V^1)$ is a  finite dimensional vector space and therefore that dim$\, H^1(BG_X)=$ dim$\, L_0/[L_0,L_0]<\infty$.

 In particular, we obtain that dim$\, L_0/L_0^{(2)}<\infty$. It follows that $L_0$ is pronilpotent and $L_0^{(n)}= L_0^n$, $n\geq 1$. Thus we may   combine Lemma 5 with \S 8 to conclude that 
 $$G_X/G_X^{n+2} \cong L_0/L_0^{(n+2)} \cong (V_n^1)^\vee$$
 where $V_0= V\cap \mbox{ker}\, d$ and $V_{n+1} = d^{-1}(\land^2V_n)$. It follows that
 $$G_X = (V^1)^\vee = \varprojlim_n (V_n^1)^\vee = \varprojlim_n G_X/G_X^{n+2},$$
 and that each $G_X^k/G_X^{k+1}$ is a finite dimensional rational vector space.

 This proves (i) and (ii), while (iii) is established directly in the Theorem. Moreover, if $X$ is Sullivan rational then by definition $X \simeq <\land V>$, where $\land V$ is its minimal Sullivan model. Thus (iv) follows from Lemma 13.
 \hfill$\square$

\section{Towers   of Lie algebras}
 
 \begin{Prop}
 \label{ptower}
 Let $L$ be the homotopy Lie algebra of a minimal algebra $\land V$. Then the following are equivalent :
 \begin{enumerate}
 \item[(i)] $L$ is the possibly finite inverse limit of a tower of finite dimensional nilpotent Lie algebras
 $$L = \varprojlim \,\,  \dots L(n)\to L(n-1)\to \dots L(0)\to 0 $$
 \item[(ii)] dim$\,V$ is finite or countably infinite
 \item[(iii)] the dimension of $H(\land V)$ is finite   or countably infinite.
 \end{enumerate}
 \end{Prop}
 
 \vspace{3mm}\noindent {\sl proof.} Denote by $d_1$ the quadratic part of the differential $d$ in $\land V$, $d_1: V\to \land^2V$. Then,   $(\land V, d_1)$ is the quadratic model of $L$. 
 
 (i) $\Longrightarrow$ (ii). Denote by $\land V_n$ the quadratic model of $L(n)$. Then the quadratic model $(\land V, d_1)$ of $L$ satisfies $\land V= \cup_n \land V_n$. Since each $V_n$ is finite dimensional, dim$\,V$ is finite or countably infinite.
 
 (ii) $\Longrightarrow$ (iii). This follows because in that case $\land V$ is countably infinite and, as vector spaces we have an injection  $H(\land V) \subset\land V$. 
 
 (iii) $\Longrightarrow$ (i). Write $V$ as the increasing union of finite dimensional subspaces $W_n$. Then define subspaces $V_n\subset W_n$ inductively by setting
 $$V_n = \{v\in W_n \, \vert\, dv\in \land V_{n-1}\}.$$
 A straightforward argument shows that $V= \cup V_n$. Evidently the homotopy Lie algebras $L(n)$ of $\land V_n$ satisfy (i).
 
 \hfill$\square$

 \vspace{3mm}\noindent {\bf Corollary.} Let $X$ be a path connected space. Then $L_X$ is the projective limit of a tower of finite dimensional nilpotent Lie algebras if and only if $H^*(X)$ is a finite type graded vector space.
 
 \vspace{3mm}\noindent {\sl proof.} Recall that  for any integer $n$, either $H^n(X)$ is finite dimensional or else uncountably infinite.   \hfill$\square$

\section{Wedge of spheres and rationally wedge-like spaces}

Suppose $\{S^{n_\alpha}\}_{\alpha \in \mathcal S}$ is a collection of spheres in which $\mathcal S$ is a linearly ordered set and each $n_\alpha \geq 1$. For each finite subset $\sigma \subset {\mathcal S}$ we write $\sigma = \{\sigma_1, \dots , \sigma_r\}$ with $\sigma_1<\dots <\sigma_r$, and set $\vert \sigma\vert = r$. If $\sigma \subset \tau$, so that $\vert \sigma\vert \leq \vert\tau \vert=q$, the inclusion defines an inclusion
$$j_{\sigma, \tau} : S^{n_{\sigma_1}}\vee \dots \vee S^{n_{\sigma_r}} \longrightarrow S^{n_{\tau_1}}\vee \dots \vee S^{n_{\tau_q}},$$
and
$$\vee_{\alpha} S^{n_\alpha} = \varinjlim_{\sigma_1<\dots <\sigma_r }\,\, S^{n_{\sigma_1}}\vee \dots \vee S^{n_{\sigma_r}}.$$

On the other hand, if $\sigma \subset \tau$, collapsing the remaining spheres to the basepoint defines a retraction
$$p_{\tau, \sigma} : S^{n_{\tau_1}}\vee\dots \vee S^{n_{\tau_q}} \longrightarrow S^{n_{\sigma_1}}\vee \dots \vee S^{n_{\sigma_r}}.$$

\vspace{3mm}\noindent {\bf Definition.} A connected space, $Y$, is \emph{rationally wedge-like} if for some ${\mathcal S}$ as above
$$Y \stackrel{\simeq}{\longrightarrow} \varprojlim_{\sigma_1<\dots <\sigma_r}\, \, (S^{n_{\sigma_1}}\vee \dots \vee S^{n_{\sigma_r}})_{\mathbb Q},$$
where the maps $\rho_{\tau, \sigma}$ defining the inverse system satisfy $\rho_{\tau, \sigma} = (p_{\tau, \sigma})_{\mathbb Q}$.

\vspace{3mm}\noindent {\bf Remark.} If each sphere has dimension $\geq 2$, and if there are only finitely many in each dimension, then $Y$ is Sullivan rational.

\begin{Prop}\label{p13}
\begin{enumerate}
\item[(i)] A minimal Sullivan algebra $\land V$ with homotopy Lie algebra, $L$, is the Sullivan model of a wedge of spheres if and only if $L$ is profree and $H(\land V)$ is the dual of a graded vector space, $H_*$.
\item[(ii)] A connected space $Y$ is rationally wedge-like if and only if for some minimal Sullivan algebra, $\land V$, 
$Y \simeq <\land V>$
and the homotopy Lie algebra of the minimal Sullivan algebra, $\land V$, is profree.
\end{enumerate}
\end{Prop}

\vspace{3mm} The Proposition depends on the next Lemma, which converts some results in Part II for quadratic Sullivan algebras to results for general minimal Sullivan algebras.

\begin{lem}\label{l12}
 The following conditions on a minimal Sullivan algebra, $(\land V,d)$, with   $V\neq 0$, are equivalent:
\begin{enumerate}
\item[(i)] There is a quasi-isomorphism $(\land V,d) \stackrel{\simeq}{\to} \mathbb Q \oplus S$ with zero differential in $S$ and $S\cdot S = 0$.
\item[(ii)] The generating space $V$ can be chosen so that $d : V\to \land^2V$ and 
$$\mathbb Q \oplus (V\cap \mbox{ker}\, d) \stackrel{\cong}{\longrightarrow} H(\land V,d).$$
\item[(iii)] $\mathbb Q \oplus V\cap \mbox{ker}\, d_1 \stackrel{\cong}{\to} H(\land V, d_1)$.
\item[(iv)] The homotopy Lie algebra, $L$, is profree.
\end{enumerate}
If they hold then $(\land V,d) \cong (\land V, d_1)$  and, if $V$ is chosen to satisfy (iii), then $V\cap \mbox{ker}\, d = V\cap \mbox{ker}\, d_1$.
 \end{lem}

\vspace{3mm}\noindent {\sl proof.} (i) $\Leftrightarrow$ (ii). First, let $S = S^{\geq 1}$ be any graded vector space. Then $\mathbb Q \oplus S$, with   zero differential and $S\cdot S = 0$ is a cdga, and its minimal model has the form 
$$(\land V,d) \stackrel{\simeq}{\longrightarrow} (\mathbb Q \oplus S,0).$$

On the other hand, a successive adjoining of variables $w$ to $S$ constructs a quadratic Sullivan algebra $(\land W, d_1)$ such that $H^{[1]}(\land W, d_1) = S$ and $H^{[2]}(\land W) = 0$. (Here as in Part II, $H^{[k]}(\land W, d_1)$ denotes the homology classes represented by cycles in $\land^kW$.) Now Lemma \ref{l9} asserts that $H^{[k]}(\land W)= 0$, $k\geq 2$. Thus division by $\land^{\geq 2}W$ and by a direct summand of $S$ in $W$ defines a quasi-isomorphism
$$(\land W, d_1) \stackrel{\simeq}{\to} \mathbb Q \oplus S.$$
Since minimal models are unique it follows that there is an isomorphism
$$(\land W, d_1) \cong (\land V,d).$$
Now choose $V$ so this isomorphism takes $W\stackrel{\cong}{\to} V$. 

On the other hand, if (ii) holds,  then with the given choice of $V$, division by $\land^{\geq 2}V$ and by a direct summand of $V\cap \mbox{ker}\, d$ in $V$ defines a quasi-isomorphism 
$$(\land V,d) \stackrel{\simeq}{\to} \mathbb Q \oplus (V\cap \mbox{ker}\, d).$$
This gives (i).

(i) $\Leftrightarrow$ (iii). If (i) holds, 
 the observations above yield an isomorphism $(\land W, d_1) \cong (\land V,d)$. This induces an isomorphism $(\land W, d_1)\stackrel{\cong}{\to} (\land V, d_1)$, which for appropriate choice of $V$  preserves wedge degree.   It follows that with this choice $(\land V,d)$ is quadratic and $V\cap \mbox{ker}\, d= V\cap \mbox{ker}\, d_1$.  This proves (iii), and the same argument shows that (iii) $\Rightarrow$ (ii).
 
Finally, the assertion (iiii) $\Leftrightarrow $ (iv) is Theorem \ref{t1}. \hfill$\square$

\vspace{3mm} Recall that a connected space, $X$, is \emph{formal} if its minimal Sullivan model is also a Sullivan model for $(H(X),0)$. Lemma 15 has th following two Corollaries :

\vspace{3mm}\noindent {\bf Corollary 1.} If $X$ is a connected space and $H^{\geq 1}(X)\neq 0$ then the homotopy Lie algebra $L_X$ is profree if and only if $X$ is formal and $H^{\geq 1}(X)\cdot H^{\geq 1}(X) \neq 0$. 

\vspace{3mm}\noindent {\bf Corollary 2.} If a connected co-H-space, $X$, satisfies $H^{\geq 1}(X)\neq 0$ (in particular, if $X$ is a wedge of spheres), then $L_X$ is profree.

\vspace{3mm}\noindent {\sl proof.} Let $U$ be a contractible open neighbourhood of a base point in $X$. Then $X\vee U$ and $U\vee X$ form an open cover for $X\vee X$
 
Since  X is a co-H-space, the diagonal map $\Delta: X\to X\times X$ factors up to homotopy through $X\vee X$,
$$\xymatrix{
X \ar[rrd]^\Delta \ar[rr]^f && X\vee X\ar[d]\\
&& X\times X.}$$
The open sets $f^{-1}(X\vee U)$ and $ f^{-1}(U\vee X)$ make then a covering of $X$ by contractible open sets, and 
 cat$\, X= 1$. It follows \cite[Theorem 9.2]{RHTII} that the minimal Sullivan model, $\land V,$ of $X$ satisfies
$$\mbox{cat}\, (\land V)= 1.$$
By Proposition 12, this   implies the result. \hfill $\square$

\vspace{3mm}\noindent {\sl proof of Proposition \ref{p13}}.  (i) Suppose first that $\land V$ is the minimal Sullivan model of a wedge of spheres, $X$. Then Corollary 2 asserts that the homotopy Lie algebra, $L$, of $\land V$ is profree. Moreover, $H(\land V)= H(X)= H_*(X;\mathbb Q)^\vee$ is the dual of a graded vector space.

In the reverse direction, since $L$ is profree, by Lemma \ref{l12}, there is a quasi-isomorphism $\land V \stackrel{\simeq}{\to} \mathbb Q \oplus S$. Moreover, by hypothesis $S= H^{\geq 1}(\land V)$ is the dual of a graded vector space $E$. Let $X= \vee_\alpha S^{n_\alpha} $ be a wedge of spheres satisfying $H_{\geq 1}(X;\mathbb Q)= E$. By Corollary 2  for the minimal Sullivan model, $\land W$, of $X$, there is a quasi-isomorphism
$$\land W \stackrel{\simeq}{\to} \mathbb Q \oplus H^{\geq 1}(X) = \mathbb Q \oplus E^\vee = \mathbb Q \oplus S.$$
It follows that $\land W \cong \land V$.

(ii) Suppose first that the homotopy Lie algebra, $L$, of a minimal Sullivan algebra, $\land V$, is profree. By Lemma \ref{l12} there is a quasi-isomorphism
$$\land V \stackrel{\simeq}{\longrightarrow} \mathbb Q \oplus S,$$
in which the differential in $S$ is zero and $S\cdot S= 0$. Fix a linearly ordered basis, $\{z_\alpha\}_{\alpha \in \mathcal S}$ of $S$ and denote deg$\, z_\alpha= n_\alpha$. Then for each finite subset $\sigma \subset \mathcal S$ write $\sigma = \{\sigma_1, \dots , \sigma_r\}$ with $\sigma_1<\dots <\sigma_r$ and set $\vert \sigma\vert = r$.

Now by induction on $\vert \sigma\vert$ we construct quasi-isomorphisms
$$\lambda_\sigma : \land V(\sigma) \stackrel{\simeq}{\longrightarrow} \mathbb Q \oplus (\oplus_i \,\mathbb Q z_{\sigma_i}),$$
and for $\gamma\subset \tau \subset \mathcal S$, wedge degree preserving Sullivan representatives
$$\lambda_{\tau, \gamma}: \land V(\gamma) \to \land V(\tau)$$
for the inclusions $\oplus_j \,\mathbb Q z_{\gamma_j} \to \oplus_i \,\mathbb Q z_{\tau_i}$.  These will be constructed so that $$\lambda_{\tau, \gamma}\circ \lambda_{\gamma, \omega} = \lambda_{\tau, \omega}\hspace{5mm}\mbox{and } \lambda_\tau \circ \lambda_{\tau, \gamma} = \lambda_\gamma.$$

In fact, suppose these have been constructed for all $\tau, \gamma$ such that $\vert \tau\vert, \vert\gamma\vert \leq r$. If $\vert \sigma\vert = r+1$ we set
$$\land V(\sigma) =\varinjlim_{\tau
 \renewcommand{\arraystretch}{0.3}
\begin{array}{c} {\scriptstyle \subset}\\{\scriptstyle \neq} \end{array}
 \renewcommand{\arraystretch}{1} 
\sigma} \, \land V(\tau) \hspace{6mm}\mbox{and } \lambda_\sigma = \varinjlim_{\tau \renewcommand{\arraystretch}{0.3}
\begin{array}{c}{\scriptstyle \subset} \\
{\scriptstyle \neq} \end{array}\renewcommand{\arraystretch}{1} \sigma} \, \lambda_\tau.$$
Then set $\lambda_{\sigma, \tau}$ to be the corresponding inclusion $\land V(\tau) \to \land V(\sigma).$

Finally, for any $\sigma_1<\dots <\sigma_r$, $\land V(\sigma)$ is the minimal Sullivan model of $S^{n_{\sigma_1}}\vee\dots \vee S^{n_{\sigma_r}}$ and $\lambda_{\tau, \sigma}$ is a Sullivan representative for the retraction $p_{\tau, \sigma}$ corresponding to the inclusion $\sigma \subset \tau$. All together this identifies $\land V = \varinjlim_\sigma \land V(\sigma)$ and
$$<\land V > = \varprojlim <\land V(\sigma)> = \varprojlim <S^{n_{\sigma_1}}\vee \dots \vee S^{n_{\sigma_r}}>_{\mathbb Q}.$$
It follows that $<\land V>$ is rationally wedge-like.

 In the reverse direction, suppose $X$ is rationally wedge like, so that
 $$X \simeq \varprojlim_{\sigma_1<\dots < \sigma_k} (S^{n_{\sigma_1}} \vee \dots \vee S^{n_{\sigma_k}})_{\mathbb Q}.$$

The explicit construction above then identifies $X = <\land V>$, where $\land V$ is the minimal Sullivan model of $\vee_\sigma S^{n_\sigma}$.

 \hfill$\square$

\vspace{3mm}\noindent {\bf Corollary.}  \emph{If $X = \vee_\sigma S^{n_\sigma}$ is a wedge of spheres, then $X_{\mathbb Q}$ is rationally wedge-like. If all the spheres are circles then $X_{\mathbb Q}$ is aspherical.}

\vspace{3mm}\noindent {\bf Remark.} Rationally wedge-like spaces provide examples of minimal Sullivan algebras $\land Z$ for which $\langle \land Z\rangle$ is not the Sullivan completion of a space. For example, suppose $Z= Z^3$ has a countably infinite basis, so that $\pi_*\langle \land Z\rangle = \pi_3\langle \land Z\rangle= (Z^3)^\vee$.
 But if $\land V$ were the minimal model of a space $X$ then we would have $V^3\cong H^3(X)= H_3(X)^\vee$ and so either dim$\, V^3<\infty$ or card$\, (V^3)\geq $ card$\,  \mathbb R$. 

\section{Wedges and free products of enriched Lie algebras}

The category of enriched Lie algebras has free products (to be constructed immediately below). Here we consider the inclusion $i : X\vee Y \to X\times Y$, and establish

\begin{Prop}\label{p14}
Suppose $X$ and $Y$ are connected spaces. Then
\begin{enumerate}
\item[(i)] The homotopy Lie algebra $L_{X\vee Y}$ of the wedge $X\vee Y$ is the free product
$$L_{X\vee Y} = L_X \,\widehat{\amalg}\, L_Y$$
of the homotopy Lie algebras of $X$ and $Y$. In particular, the correspondence $X\mapsto L_X$ preserves coproducts.
\item[(ii)] If one of $H(X)$, $H(Y)$ is a graded vector space of finite type, then the  homotopy fibre, $F$, of the map $i_{\mathbb Q} : (X\vee Y)_{\mathbb Q}\to (X\times Y)_{\mathbb Q}$ is rationally wedge-like.
\item[(iii)] If $X_{\mathbb Q}$ and $Y_{\mathbb Q}$ are aspherical then so are $F$ and $(X\vee Y)_{\mathbb Q}$.
\end{enumerate}
\end{Prop}

\vspace{3mm}\noindent {\bf Remark.}  This result is analogous to the fact that the usual fibre of the injection $X\vee Y\to X\times Y$ is the join of $\Omega X$ and $\Omega Y$ and thus a suspension. (But note that $(X\vee Y)_{\mathbb Q}$ may be different from $X_{\mathbb Q}\vee Y_{\mathbb Q}$.)

\vspace{3mm} The main step in the proof of Proposition \ref{p14} is the explicit description of the homotopy Lie algebra of a fibre product $\land W\times_{\mathbb Q} \land Q$ of any two minimal Sullivan algebras. After defining $L_X\, \widehat{\amalg}\, L_Y$ we establish this description in Proposition 22, and  return to Proposition \ref{p14} and its application to Sullivan completions.

First, note that the classical construction of the free product, $L \amalg L'$, of two graded Lie algebras  extends naturally to  enriched Lie algebras, $(L, \{I_\alpha\})$ and $(L', \{I'_\beta\})$, in which the enriched structure is given by the surjections
$$\xi_{\alpha, \beta, n} : L\,\amalg\, L'\to L_\alpha\, \amalg\, L_\beta \,/\, (L_\alpha \, \amalg\, L_\beta)^n.$$
Its completion will be denoted by $L\, \widehat{\amalg}\, L'$, and almost by definition,
$$L\, \widehat{\amalg}\, L' \stackrel{\cong}{\longrightarrow} \overline{L}\, \widehat{\amalg}\, \overline{L'}.$$

\begin{lem}\label{l13}  With the notation above, any two morphisms $f : (L, \{I_\alpha\}) \to E$ and $g : (L', \{I_\beta '\})\to E$ into a complete enriched Lie algebra  extend  uniquely to a morphism
$$L\,\widehat{\amalg}\, L' \to E.$$
This characterizes $L\,\widehat{\amalg}\,L'$ up to natural isomorphism.
%\item[(ii)] If $L= \varinjlim_\sigma L(\sigma)$ and $L'= \varinjlim_\tau L'(\tau)$ are the direct limits of closed sub %Lie algebras then $L(\sigma)\,\widehat{\amalg}\, L'(\tau)$ is a closed sub Lie algebra of $L\,\widehat{\amalg}\,L'$ and
%$$L\,\widehat{\amalg}\, L' = \varinjlim_{\sigma, \tau} L(\sigma)\, \widehat{\amalg}\, L'(\tau).$$
%\end{enumerate}
\end{lem}

\vspace{3mm}\noindent {\sl proof.} Let $\{J_\gamma\}$ denote the enriched structure for $E$. Then for each $\gamma$ there are indices $\alpha (\gamma), \beta(\gamma)$ such that $f$ and $g$ factor to yield morphisms
$$f_\gamma : L_{\alpha (\gamma)} \to E_\gamma \hspace{5mm}\mbox{and } g_\gamma : L'_{\beta (\gamma)} \to E_\gamma.$$
Moreover, since $E_\gamma^{n(\gamma)} = 0$, some $n(\gamma)$, $f_\gamma \amalg g_\gamma$ factors to yield a morphism
$$h_\gamma : (L_{\alpha (\gamma)} \amalg L'_{\beta(\gamma)}) / (L_{\alpha (\gamma)} \amalg L'_{\beta (\gamma)})^n \longrightarrow E_\gamma.$$
Since $E$ is complete these extend to 
$$h : = \varprojlim_\gamma h_\gamma : L\,\widehat{\amalg}\, L'\to E.$$
  The uniqueness is immediate. \hfill $\square$

\vspace{3mm}
Now, consider minimal Sullivan algebras,   $\land W$ and $\land Q$. The natural surjection $\land W\otimes \land Q \to \land W\times_{\mathbb Q}\land Q$ is surjective in homology, and so the corresponding $\Lambda$-extension, 
$$\varphi : \land T:= \land W\otimes \land Q\otimes \land R \stackrel{\simeq}{\longrightarrow} \land W\times_{\mathbb Q}\land Q,$$ in which $\varphi (R)= 0$, is a minimal Sullivan model for $\land W\oplus_{\mathbb Q} \land U$.  
The sequence,  
$$\land W\otimes \land Q \to \land W\otimes \land Q\otimes \land R \to \land R$$
then (\S 8)  induces the short exact sequence
$$0 \leftarrow L_W\times L_Q \leftarrow L_T\leftarrow L_R\leftarrow 0$$
of complete enriched Lie algebras.  

\begin{Prop}
\label{p15} With the hypotheses and notation above, 
\begin{enumerate}
\item[(i)] The surjection $L_T\to L_W\times L_Q$ factors as 
$$L_T\stackrel{\cong}{\to} L_W\, \widehat{\amalg}\, L_Q \to L_W\times L_Q.$$
\item[(ii)] $L_R$ is a profree Lie algebra. In particular   $<\land R>$ is rationally wedge-like.
\item[(iii)] When $\land W$ and $\land Q$ are quadratic Sullivan algebras, then, for appropriate choice of $R$,  $\land T$ is the quadratic model of $L_W\, \widehat{\amalg}\, L_Q$.
\item[(iv)] If $f : L\to L'$ and $g : E\to E'$ are surjections of complete enriched Lie algebras, then
$$f\, \widehat{\amalg}\, g: L\, \widehat{\amalg}\, E\to L'\, \widehat{\amalg}\, E'$$
is also surjective.
\end{enumerate}
\end{Prop}

\vspace{3mm}\noindent {\bf Corollary.} If $L$ and $L'$ are complete enriched Lie algebras then the kernel of 
$$L\, \widehat{\amalg}\, L' \to L\times L'$$
is profree. Moreover, if $L$ and $L'$ are both profree then $L\,\widehat{\amalg}\, L'$ is also profree.

\vspace{3mm}\noindent {\sl proof of the Corollary.} For the first assertion apply Proposition \ref{p15} to the quadratic models $\land W$ and $\land W'$ of $L$ and $L'$. If now $L$ and $L'$ are profree then by Theorem \ref{t1} there are quasi-isomorphisms
$$\land W \stackrel{\simeq}{\longrightarrow} \mathbb Q \oplus S \hspace{5mm}\mbox{and } \land W' \stackrel{\simeq}{\longrightarrow} \mathbb Q \oplus S'$$
with $S\cdot S = 0= S'\cdot S'$, and vanishing differentials in $S$ and $S'$. Thus $\land W\times_{\mathbb Q}\land W' \simeq \mathbb Q\oplus S\oplus S'$, and its homotopy Lie algebra, $L\, \widehat{\amalg}\, L'$, is profree by Theorem \ref{t1}. \hfill$\square$

\vspace{3mm} \vspace{3mm}\noindent {\sl proof of Proposition \ref{p15}.} (i) Denote by $[\land M, \land N]_*$ the set of based homotopy classes of morphisms. Then for any minimal Sullivan algebra, $\land V$, composition with $\varphi$ induces a bijection (\cite[Proposition 1.10]{RHTII}),
$$[\land V, \land T]_* \to [\land V, \land W\times_{\mathbb Q} \land Q]_*,$$
of based homotopy classes of morphisms. The surjections $\land W\times_{\mathbb Q}\land Q\to \land W, \land Q$ then identify $[\land V, \land W\times_{\mathbb Q}\land Q]_* = [\land V, \land W]_* \times [\land V, \land Q]_*$. Now apply Lemma \ref{l13} to convert these bijections to a bijection,
$${\mathcal C}(L_T, L_V) \stackrel{\cong}{\longrightarrow} {\mathcal C}(L_W, L_V) \times {\mathcal C}(L_Q, L_V),$$
where as usual, ${\mathcal C}(-,-)$ denotes the set of morphisms in the category of enriched Lie algebras. This identifies $L_T$ as $L_W\, \widehat{\amalg}\, L_Q$ and establishes (i).

 (ii) Let $\land W\otimes \land U_W$ and $\land Q\otimes \land U_Q$ denote the respective acyclic closures. 
Then  $\land R$ is quasi-isomorphic to
 $$\land T \otimes_{\land W\otimes \land Q} (\land W\otimes \land U_W\otimes \land Q\otimes \land U_Q) \simeq A:= (\land W\times_{\mathbb Q}\land Q) \otimes \land U_W\otimes \land U_Q.$$
Now divide $A$ by the ideal generated by $W$ to obtain the short exact sequence
 $$0\to \land^{\geq 1}W\otimes \land U_W\otimes \land U_Q \to A \to  \land Q\otimes \land U_W\otimes \land U_Q \to 0.$$
 
 Next, decompose the differential in $\land W\otimes \land U_W$ in the form $d= d_1+d'$ with $d_1(W)\subset \land^2W$, $d_1(U_W)\subset W\otimes \land U_W$, $d'(W)\subset \land^{\geq 3}W$ and $d'(U_W)\subset \land^{\geq 2}W\otimes \land U_W$. Then $d_1$ is a differential and $(\land W\otimes \land U_W,d_1)$ is the acyclic closure of $(\land W,d_1)$. Choose a direct summand,  $S$, of $d_1(\land^{\geq 1}U_W) $ in $W\otimes \land U_W$. Then $I = (\land^{\geq 2}W\otimes \land U_W)\oplus S$ is acyclic for the differential $d_1$ and therefore also for the differential $d$. 
Thus $J = I\otimes \land U_Q$ is an acyclic ideal in $A$ and $A\stackrel{\cong}{\to} A/J$.  

On the other hand, consider the short exact sequence
$$0\to d(\land U_W)\otimes \land U_Q\to A/J\to \land Q\otimes \land U_W\otimes \land U_Q\to 0$$
in which $d(\land U_W)\otimes \land U_Q$ is an ideal with trivial multiplication and trivial differential. It follows from the long homology sequence that
$$\mathbb Q \oplus d(\land U_W)\otimes \land^{\geq 1}U_Q \stackrel{\simeq}{\longrightarrow} A/J.$$
Thus (ii) follows from Lemma \ref{l13}.

(iii) Assign  $\land W$ and $\land Q$ wedge degree as a second degree and assign $U_W$ and $U_Q$   second degree $0$. Then $(\land W\otimes \land U_W,d_1)$ and $(\land Q\otimes \land U_Q,d_1)$ are the respective acyclic closures of $(\land W,d_1)$ and $(\land Q, d_1)$, and $d_1$ increases the second degree by 1. Now $ {\varphi}$ and $ {T}$ 
 may be constructed so that $ {R}$ is equipped with a second gradation for which $d$ increases the second degree by one and $ {\varphi}$ is bihomogeneous of degree zero. 

The argument in the proof of (ii) now yields a sequence of bihomogeneous quasi-isomorphisms connecting
$$\mathbb Q \oplus \left( d_1(\land^+U_W)\otimes \land^+U_Q\right)\simeq \land  {R}.$$
Thus $H^{\geq 1}(\land {R})$ is concentrated in second degree 1. Therefore   $ {R}$ is concentrated in second degree 1 and so the second degree in $R$ (and therefore in $T$) is just the wedge degree. Thus $\land T$ is quadratic and therefore is the quadratic model of $L_W\, \widehat{\amalg}\, L_Q$.

(iv) It is sufficient to consider the case $g = id_E : E\stackrel{=}{\longrightarrow} E$. In this case we have the row-exact commutative diagram
$$\xymatrix{
0\ar[r] &  K\ar[rr]\ar[d]^h && L\, \widehat{\amalg}\, E\ar[d]^{f\, \widehat{\amalg}\, id_E} \ar[rr] &&   L\times E\ar[r]\ar[d]^{f\times id_E} & 0\\
0 \ar[r] & K' \ar[rr] && L'\, \widehat{\amalg}\, E \ar[rr] && L'\times E \ar[r] & 0,}$$
and so it is also sufficient to show that $h$ is surjective.

Now let $\land W$, $\land W'$ and $\land Q$ be, respectively the quadratic models of $L$, $L'$ and $E$, and let $\varphi : \land W'\to \land W$ be the morphism corresponding to $f$. Then the restriction $\varphi : W'\to W$ dualizes to $f$, and so $\varphi$ is injective. In particular, $\varphi$ extends to an inclusion
$$\varphi : \land W'\otimes \land U_{W'}\to \land W\otimes \land U_W$$
of the acyclic closures, restricting to an inclusion $U_{W'}\to U_W$.

On the other hand, the diagram above is the diagram of homotopy Lie algebras associated with a commutative diagram,
$$\xymatrix{\land W\otimes \land Q \ar[rr] && \land W\otimes \land Q\otimes \land R \ar[rr] && \land R\\
\land W'\otimes \land Q\ar[u]^{\varphi\otimes id} \ar[rr] && \land W'\otimes \land Q\otimes \land R' \ar[u]\ar[rr] && \land R'\ar[u]^\sigma,}$$
which identifies $L_\sigma$ with $h$.

Moreover, the computations in the proof of (ii) identify
$$H(\sigma) = \varphi : \mathbb Q \oplus d(\land U_{W'}) \to \mathbb Q \oplus d(\land U_W).$$
It follows that $H(\sigma)$ is injective. But since $L_R$ and $L_{R'}$ are profree, Lemma \ref{l5} identifies
$$L_R/L_R^{(2)} = H^{\geq 1}(\land R)^\vee \hspace{5mm}\mbox{and } L_{R'}/L_{R'}^{(2)} = H^{\geq 1}(\land R')^\vee.$$
It follows that $L_\sigma : L_R\to L_{R'}$ induces a surjection $L_R/L_R^{(2)}\to L_{R'}/L_{R'}^{(2)}$. Now Lemma \ref{l8} (i) implies that $L_\sigma$ is surjective.
\hfill$\square$

\vspace{3mm}\noindent {\sl proof of Proposition \ref{p14}.} Let $\land W$ and $\land Q$ be the minimal Sullivan models of $X$ and $Y$. Then recall from Proposition \ref{p15} the quasi-isomorphism
$$\land T := \land W\otimes \land Q\otimes \land R \stackrel{\simeq}{\longrightarrow} \land W\times_{\mathbb Q}\land Q.$$
Since $H(X\vee Y) =H(X)\times_{\mathbb Q} H(Y)$ it follows that $A_{PL}(X\vee Y) \stackrel{\simeq}{\to} A_{PL}(X)\times_{\mathbb Q} A_{PL}(Y)$. Therefore $\land T$ is a minimal Sullivan model of $X\vee Y$, and so (i) follows from Proposition \ref{p15}(i).

(ii) If one of $H(X)$, $H(Y)$ has finite type then $H(X\times Y)= H(X)\otimes H(Y)$. It follows that the inclusion
$$\land W\otimes \land Q\to \land T$$
is a Sullivan representative for the inclusion $i : X\vee Y\to X\times Y$. This identifies $i_{\mathbb Q}$ as the map $<\land T>\to <\land W\otimes \land Q>$. But by \cite[Proposition 17.9]{FHTI} the sequence
$$<\land R>\to <\land T>\stackrel{i_{\mathbb Q}}{\longrightarrow} <\land W\otimes \land Q>$$
is a fibration. Thus (Proposition \ref{p15}) the fibre, $<\land R>$ of $i_{\mathbb Q}$ is rationally wedge-like.

(iii) When $X_{\mathbb Q}$ and $Y_{\mathbb Q}$ are aspherical then $W = W^1$ and $Q= Q^1$; Now it follows from Proposition \ref{p15}(iii) that $R= R^1$ and so $F = <\land R>$ is aspherical. But since $R = R^1$, $T= T^1$ and so $<\land T>= (X\vee Y)_{\mathbb Q}$ is also aspherical. \hfill$\square$

\subsection{Free products of completions of weighted Lie algebras}
 
\vspace{1mm} Let $E$ and $E'$ be weighted enriched Lie algebras with enriched structures given by morphisms $E\to E_\alpha$ and $E'\to E_\beta'$.  
 Denote by $L$ and $L'$ the completions of $E$ and $E'$.  The weight decomposition in $E$ and $E'$ then extends in the standard way to $E\,\amalg\, E'$. The corresponding decomposition of the quadratic model for $E\, \amalg\, E'$ is given as follows. The decomposition of the models $\land V$ and $\land V'$ induce one in $\land V\oplus_{\mathbb Q}\land V'$ which then in turn defines the decomposition in its minimal quadratic model. 
 
 \begin{Prop}
 \label{pXfree}  The completion of $E\, \amalg\, E'$ is $L\, \widehat{\amalg}\, L'.$
 \end{Prop}
 
 \vspace{3mm}\noindent {\sl proof.}   This follows directly from Proposition 22(iii).

\hfill$\square$

\newpage\noindent {\bf \LARGE Appendix : Connected spaces and minimal Sullivan models}

\vspace{1cm} For the convenience of the reader we recall briefly some of the definitions and main elements of Sullivan's approach (\cite{S}) to the rational homotopy theory of a connected space. A more detailed survey and  proofs are provided respectively in \cite{Sur} and \cite{RHTII}. (Here by space we mean either a CW complex or a simplicial set;  these categories are identified ny the inverse homotopy equivalences provided by the singular simplex and Milnor realization functor (\cite[\S 16]{May}).

The central concept is the Sullivan algebra, $(\land V,d)$. These are the   commutative graded differential algebras (cdga's), for which 
\begin{enumerate}
\item[$\bullet$] $\land V$ is the free graded commutative algebra on a graded vector space $V = V^{\geq 1}$, and
\item[$\bullet$] the differential is constrained by the condition $V = \cup_{n\geq 0} V_n$, where $V_0 = V\cap \mbox{ker}\, d$ and $V_{n+1}= V\cap d^{-1}(\land V_n)$.
\end{enumerate}
Here $\land V= \oplus_{p\geq 0} \land^pV$, where $\land^pV$ is the linear span of the monomials in $V$ of length $p$.
In particular, $(\land V,d)$ is \emph{minimal} if $d: V\to \land^{\geq 2}V$. We frequently suppress the differential from the notation, and write $\land V$ for $(\land V,d)$. Then two cdga morphisms $\varphi_0, \varphi_1 : \land V\to A$ are \emph{homotopic} if for some morphism $\varphi : \land V\to A\otimes \land (t,dt)$ with deg$\, t= 0$ the augmentations $t\mapsto 0,1$ convert $\varphi$ to $\varphi_0, \varphi_1$; $\varphi_0$ and $\varphi_1$ are \emph{based homotopic} with respect to an augmentation $\varepsilon_A : A\to \mathbb Q$ if $\varphi : V\to \mbox{ker}\, \varepsilon_A \otimes \land (t,dt)$.  

We also recall  \cite[\S 6]{FHTI} that a $\land V$-module, $P$, is a graded vector space equipped with a differential, $d$ and a multiplication $\land V\otimes P\to P$ compatible with the differential. The $\land V$-module is \emph{semifree} if it has the form $P= \land V \otimes M$ in which the module structure is simply multiplication by $\land V$ and, in addition $P = \cup_{k\geq 0} \land V\otimes M_k$ is the increasing union of submodules in which the quotient differentials in $\land V\otimes M_{k+1}/M_k$ are just $d\otimes id$.

Now, suppose $(A, \varepsilon_A)$ is an augmented cdga satisfying $H^0(A) = \mathbb Q$. Then a $\Lambda$-extension of $(A, \varepsilon_A)$ is a sequence of morphisms
$$\xymatrix{A \ar[rr]^\lambda && A\otimes \land Z \ar[rr]^\rho && \land Z}$$
in which $\lambda (a) = a\otimes 1$, $\rho = \mathbb Q\otimes_A -$, and the differential is constrained by the conditions: $Z = \cup_{n\geq 0} Z_n$ with $Z_0 = Z\cap d^{-1}(A\otimes 1)$ and $Z_{n+1}= Z\cap d^{-1}(A\otimes \land Z_n)$. It is \emph{minimal} if the quotient cdga, $(\land Z, \overline{d})$ satisfies $\overline{d}: Z\to \land^{\geq 2}Z$. 

Every cdga, $A$, with $H^0(A) = \mathbb Q$ admits a quasi-isomorphism $\sigma : \land V \stackrel{\simeq}{\to} A$ from a minimal Sullivan algebra; this is a \emph{minimal Sullivan model} for $A$, and $\land V$ is unique up to isomorphism. Similarly, if $\varphi : A\to B$ is a morphism of augmented cdga's satisfying $H^0= \mathbb Q$ then $\varphi$ extends to a quasi-isomorphism
$$A\otimes \land Z\stackrel{\simeq}{\to} B$$
from a minimal $\Lambda$-extension. Again, $A\otimes \land Z$ is unique up to isomorphism.

Now consider the case in which $A$ is a minimal Sullivan algebra $\land V$, and $H^1(\varphi)$ is injective. In this case the quotient $(\land Z, \overline{d}) = \land V\otimes \land Z\otimes_{\land V}\mathbb Q$ will also be a minimal Sullivan algebra. Write $\land V\otimes \and Z = \land (V\oplus Z)$ and note that the differential is the sum
$$d= \sum_{n\geq 0} d_n$$
of derivations in which $d_k$ raises the wedge degree by $k$. In particular, $\land V\otimes \land Z$ is a minimal Sullivan algebra if and only if $d_0= 0$.

In general, $d_0$ is a linear map $Z\to V$. Its dual is a linear map $\partial : L_Z\leftarrow L_V$ between the respective homotopy Lie algebras $ L_V = (sV)^\vee$ and $L_Z= (sZ)^\vee$ defined in \cite[Chap 2]{RHTII} and recalled in \S 8. Here we have the

\vspace{3mm}\noindent {\bf Lemma A.1} {\sl With the hypotheses and notation above, the image of $\partial$ is contained in the centre of $L_Z$.}

\vspace{3mm}\noindent {\sl proof.} We may identify $\land^2Z$ as a subspace of $Z\otimes Z$. Then a single calculation gives
$$<(d_0\otimes id)\overline{d}_1z, s\partial x\otimes sy> = \pm <\overline{d}_1z, s\partial x\otimes sy>.$$
Since $(d_0\otimes id)\overline{d}_1z$ is the component in $V\otimes Z$ of $d^2(1\otimes z)$ it follows that $<\overline{d}_1z, s\partial x\otimes sy>= 0$. On the other hand, the Lie bracket in $L_Z$ is given by
$$<\overline{d}_1z, sy_1\otimes sy_2> = (-1)^{1+ deg\, y_1} <z_1, s[y_1, y_2]>.$$
It follows that for $x\in L_V$, $y\in L_S$ we have $[\partial x,y]= 0$. 
\hfill$\square$

\vspace{3mm}For the connection between Sullivan algebras and topology we identify spaces with simplicial sets via the singular simplex functor, and refer to the objects in either category as 'spaces'. 
The connection between Sullivan algebras and topology is provided by
 the simplicial commutative cochain algebra, $A_{PL}$, defined as follows:
 $$(A_{PL})_n =\frac{\land (t_0, \dots , t_n, dt_0, \dots , dt_n)}{\sum t_i-1, \sum dt_i}.$$
 The simplicial faces $\partial_i$ and degeneracies $s_j$ are defind by their restrictions to the $t_k$ 
 $$\partial_it_k= \left\{\begin{array}{ll} t_k & \mbox{if } k<i\\
 0 & \mbox{if } k=i\\
 t_{k-1} & \mbox{if } k>i\end{array}\right. \hspace{15mm}\mbox{and } s_jt_k= \left\{\begin{array}{ll} t_k & \mbox{if } k<j\\
 t_k+t_{k+1} & \mbox{if } k=j\\
 t_{k+1} & \mbox{if } k>j.
 \end{array}\right.$$
  The \emph{algebra of polynomial forms} on $X$ is then defined by
 $$A_{PL}^n(X)= \mbox{Simpl} (X, A_{PL}^n),$$
 and there is a natural isomorphism $H(A_{PL}(X)) \cong H^*(X;\mathbb Q)$. 
 
 A \emph{Sullivan minimal model} of $X$ is then a minimal model for $A_{PL}(X)$, i.e., a minimal Sullivan algebra equipped with a quasi-isomorphism
 $$\varphi : \land V \stackrel{\simeq}{\to}   A_{PL}(X).$$ The minimal model exists always and is unique up to isomorphism. 
 
 Similarly, if $f : X\to Y$ is a map of pointed connected spaces then $A_{PL}(f)$ extends to a quasi-isomorphism
 $$A_{PL}(Y)\otimes \land Z\stackrel{\simeq}{\to} A_{PL}(X)$$
 from a minimal $\Lambda$-extension; here the quotient $(\land Z, \overline{d})$ is the analogue of the homotopy fibre, $F$, of $f$. In fact, there is a natural morphism $(\land Z, \overline{d})\to A_{PL}(F)$ which, with certain hypotheses is the minimal Sullivan model (\cite[Theorem 5.1]{RHTII}).

 Sullivan (\cite{S},\cite{RHTII}) also introduces the \emph{realization functor}, $A \mapsto <A>$, from cdga's to simplicial sets  adjoint to $A_{PL}(-)$:    
$$<A>_p = \mbox{{Cdga}}(A, (A_{PL})_p),$$
Since the realization functor converts tensor products to products, it follows that $<A\otimes B> = <A> \times <B>$. 

For a connected space $X$ this yields a natural bijection,
 $$\mbox{Cdga} (A, A_{PL}(X) \cong \mbox{Simpl}(X, <A>),$$ denoted $\varphi \mapsto <\varphi>$. This    extends, for Sullivan algebras,  to a natural bijection of homotopy classes
$$[\land V, A_{PL}(X)] \cong [X, <\land V>].$$

Applied to a minimal model $\varphi : \land V\stackrel{\simeq}{\to} A_{PL}(X)$, this yields   a natural map
$$<\varphi> : X \to X_{\mathbb Q}:= <\land V>,$$
 the \emph{Sullivan rationalization} or \emph{completion} of $X$.   When $X$ is simply connected and $H(X)$ is a   graded vector space of finite type, then (\cite{S}) $<\varphi>$ is a rational homotopy equivalence, and $<\varphi>$ induces isomorphisms of graded vector spaces
$$\pi_n(X)\otimes \mathbb Q\to \pi_n(X_{\mathbb Q}).$$ 
 The hypothesis on the $H(X)$ is important as we see in the next Proposition

\vspace{3mm}\noindent {\bf Proposition.} When $X$ is simply connected and $H(X)$ is not finite type, then $<\varphi>$ is not a rational homotopy equivalence.

\vspace{3mm}\noindent {\sl proof}. Since some Betti numbers are infinite, there is a smallest integer $p$ with dim$\, \pi_p(X)\otimes \mathbb Q=\infty$. Suppose first that $X$ is $(p-1)$-connected and write $E= \pi_p(X)\otimes \mathbb Q$. Then by Hurewicz theorem, $H^p(X;\mathbb Q) \cong E^\vee$. Let $(\land W,d)\to A_{PL}(X)$ be the minimal model of $X$, then $W^p\cong H^p(X;\mathbb Q)\cong E^p$. It follows that $\pi_p<\land W>$ is the bidual of $E$. Now since $E$ is infinite, $E$ and its bidual are not isomorphic, and this proves the result in this case.

In the general situation, consider the Postnikov fibration
$$X^{(p)}\stackrel{i}{\to} X\stackrel{p}{\to} Y$$
where $\pi_q(Y)= 0$ for $q\geq p$, $\pi_q(p)$ is an isomorphism for $q<p$ and $\pi_q(i)$ is an isomorphism for $i\geq p$. By \cite[Theorem 5.1]{RHTII}, we have a commutative diagram   in which $(\land V\otimes \land W,d)$ is a relative minimal model and the vertical arrows are quasi-isomorphisms,
$$\xymatrix{
A_{PL}(Y) \ar[r] & A_{PL}(X) \ar[r] & A_{PL}(X^{(p)})\\
\land V \ar[r]\ar[u]^\simeq & \land V\otimes \land W\ar[r]\ar[u]_\simeq & (\land W, \overline{d})\ar[u]_\simeq}$$
By \cite[Proposition 17.9]{FHTI} we deduce a commutative diagram of fibrations
$$\xymatrix{
<\land W> \ar[r] & <\land V\otimes \land W> \ar[r] & <\land V>\\
X^{(p)}\ar[r] \ar[u]^{r_{X^{(p)}}} & X\ar[r]\ar[u]^{r_X} & Y\ar[u]^{r_Y}
}$$

Since $H(Y)$ is a finite type graded vector space, $r_Y$ is a rational homotopy equivalence. If $r_X$ is a rational homotopy equivalence, then the same is true for $r_{X^{(p)}}$ which is not the case. Therefore $r_X$ is not a rational homotopy equivalence.
\hfill$\square$

 \vspace{5mm}\noindent Institut de Math\'ematique et de Physique, Universit\'e Catholique de Louvain, 2, Chemin du cyclotron, 1348 Louvain-La-Neuve, Belgium, yves.felix@uclouvain.be

 \vspace{1mm}\noindent Department of Mathematics, Mathematics  Building, University of Maryland, College Park, MD 20742, United States, shalper@umd.edu
\end{document}